\documentclass[a4paper,12]{article}
\usepackage{amsmath, amsfonts, amssymb, amsthm, bm, graphics, bbm, color}
\usepackage{graphicx}
\usepackage{subfig}
\usepackage{pgfplots}
\usepackage{mathtools}
\usepackage{xcolor,colortbl}
\usepackage{xurl}
\usepackage{mathabx}
\usepackage{cancel}
\usepackage{multicol}
\usepackage{float}
\usepackage{caption}
\usepackage{mathalfa}
\usepackage{hyperref}
\usepackage{afterpage}
\usepackage{multirow}
\usepackage{amssymb}
\usepackage{tikz} 
\usepackage{lscape} 
\usepackage{longtable}
\usepackage{algorithm}
\usepackage{algpseudocode}
\usepackage{graphbox} 

\usepackage{tikz}
\usetikzlibrary{patterns}

\DeclareMathAlphabet\mathbfcal{OMS}{cmsy}{b}{n}
\newtheorem{theorem}{Theorem}[section]
\newtheorem{lemma}[theorem]{Lemma}

\theoremstyle{definition}

\newtheorem{remark}[theorem]{Remark}

\makeatletter
\newcommand\makebig[2]{%
 \@xp\newcommand\@xp*\csname#1\endcsname{\bBigg@{#2}}%
 \@xp\newcommand\@xp*\csname#1l\endcsname{\@xp\mathopen\csname#1\endcsname}%
 \@xp\newcommand\@xp*\csname#1r\endcsname{\@xp\mathclose\csname#1\endcsname}%
}
\makeatother

\makebig{biggg} {3.0}
\makebig{Biggg} {3.5}
\makebig{bigggg}{4.0}
\makebig{Bigggg}{4.5}
\makebig{biggggg}{5.0}
\makebig{Biggggg}{5.5}
\makebig{bigggggg}{6.0}
\makebig{Bigggggg}{6.5}

\newcommand{\dif}{\mathrm{d}}
\newcommand{\im}{\mathrm{i}}

\pgfplotsset{compat=1.15}

\begin{document}
\title{Efficient Computation of Magnetic Polarizability Tensor Spectral Signatures for Object Characterisation in Metal Detection}

\author{J. Elgy and P.D. Ledger\\
School of Computer Science and Mathematics, Keele University\\
Keele, Staffordshire,  ST5 5BG.  U.K\\
corresponding author: j.elgy@keele.ac.uk}
\maketitle

\section*{Abstract}

%
%

\textbf{Purpose:} Magnetic polarizability tensors (MPTs) provide an economical characterisation of conducting magnetic metallic objects and their spectral signature can aid in the solution of metal detection inverse problems, such as scrap metal sorting, searching for unexploded ordnance in areas of former conflict, and security screening at event venues and transport hubs. In this work, the authors discuss methods for efficiently building large dictionaries for classification approaches. 

\noindent\textbf{Design/methodology/approach:} Previous work has established explicit formulae for MPT coefficients, underpinned by a rigorous mathematical theory. To assist with the efficient computation of MPTs at differing parameters and objects of interest this work applies new observations about the way the MPT coefficients can be computed. Furthermore, the authors discuss discretisation strategies for $hp$--finite elements on meshes of unstructured tetrahedra combined with prismatic boundary layer elements for resolving thin skin depths and using an adaptive proper orthogonal decomposition (POD) reduced order modelling methodology to accelerate computations for varying parameters.  

\noindent\textbf{Findings:} The success of the proposed methodologies is demonstrated using a series of examples. A significant reduction in computational effort is observed across all examples. The authors identify and recommend a simple discretisation strategy, and improved accuracy is obtained using adaptive POD. 

\noindent\textbf{Originality:} The authors present novel computations, timings, and error certificates of MPT characterisations of realistic objects made of magnetic materials. A novel postprocessing implementation is introduced, and an adaptive POD algorithm is demonstrated. 

\noindent\textbf{Keywords:} Metal detection, magnetic polarizability tensor, eddy current problems, reduced order models, thin skin depth effects. 

\noindent\textbf{Paper type:} Research paper

\section{Introduction}\label{sect:intro}

Metal detection uses low frequency magnetic induction to locate and identify highly conducting objects within the imaging region, where objects are typically made up of one or more metallic materials. Most metals can be regarded as being highly conducting, having a conductivity in excess of $10^6$ S/m, and many metals are also magnetic with a permeability much larger than that of free space (given by $\mu_0 = 4 \pi \times10^{-7} $ H/m). Traditional metal detectors rely on simple thresholding of measured field perturbations, {which} are mapped to an audible tone to aid with detection. For this reason, hobbyist metal detectors  wear headphones, listen for changes in pitch and volume, and have become trained in recognising the different signals that a detector receives for targets made of different materials (such as coins and other buried treasure)  buried at different depths. 
Here, small objects close to the surface may give rise to similar signals to those of larger objects buried at greater depth and false positives are common. There are other safety critical applications of metal detection technology, such as in the location and identification of hidden unexploded ordnance (UXO) in areas of former conflict, where minimising the number of false positives and false negatives is vital.
This work presents new computational procedures and algorithmic developments that can aid with accurately and efficiently characterising highly conducting and magnetic objects  for use in classification approaches with the aim of reducing the number of both false positives  and false negatives.

\begin{figure}[h!]
\centering
$\begin{array}{c c}
\raisebox{-0.5\height}{
 
\tikzset{
pattern size/.store in=\mcSize, 
pattern size = 5pt,
pattern thickness/.store in=\mcThickness, 
pattern thickness = 0.3pt,
pattern radius/.store in=\mcRadius, 
pattern radius = 1pt}
\makeatletter
\pgfutil@ifundefined{pgf@pattern@name@_axzpthoyo}{
\pgfdeclarepatternformonly[\mcThickness,\mcSize]{_axzpthoyo}
{\pgfqpoint{0pt}{0pt}}
{\pgfpoint{\mcSize+\mcThickness}{\mcSize+\mcThickness}}
{\pgfpoint{\mcSize}{\mcSize}}
{
\pgfsetcolor{\tikz@pattern@color}
\pgfsetlinewidth{\mcThickness}
\pgfpathmoveto{\pgfqpoint{0pt}{0pt}}
\pgfpathlineto{\pgfpoint{\mcSize+\mcThickness}{\mcSize+\mcThickness}}
\pgfusepath{stroke}
}}
\makeatother

 
\tikzset{
pattern size/.store in=\mcSize, 
pattern size = 5pt,
pattern thickness/.store in=\mcThickness, 
pattern thickness = 0.3pt,
pattern radius/.store in=\mcRadius, 
pattern radius = 1pt}
\makeatletter
\pgfutil@ifundefined{pgf@pattern@name@_9b2vxhdkv}{
\pgfdeclarepatternformonly[\mcThickness,\mcSize]{_9b2vxhdkv}
{\pgfqpoint{0pt}{0pt}}
{\pgfpoint{\mcSize+\mcThickness}{\mcSize+\mcThickness}}
{\pgfpoint{\mcSize}{\mcSize}}
{
\pgfsetcolor{\tikz@pattern@color}
\pgfsetlinewidth{\mcThickness}
\pgfpathmoveto{\pgfqpoint{0pt}{0pt}}
\pgfpathlineto{\pgfpoint{\mcSize+\mcThickness}{\mcSize+\mcThickness}}
\pgfusepath{stroke}
}}
\makeatother
\tikzset{every picture/.style={line width=0.75pt}} 

\begin{tikzpicture}[x=0.75pt,y=0.75pt,yscale=-0.6,xscale=0.6]

\draw  [draw opacity=0][pattern=_axzpthoyo,pattern size=8.774999999999999pt,pattern thickness=0.75pt,pattern radius=0pt, pattern color={rgb, 255:red, 0; green, 0; blue, 0}] (80,250) -- (420,250) -- (420,260) -- (80,260) -- cycle ;
\draw  [fill={rgb, 255:red, 255; green, 255; blue, 255 }  ,fill opacity=1 ] (216,281.33) .. controls (236,271.33) and (326,261.33) .. (306,281.33) .. controls (287.21,300.12) and (286.07,310.09) .. (302.59,336.11) .. controls (303.66,337.78) and (304.79,339.52) .. (306,341.33) .. controls (326,371.33) and (236,371.33) .. (216,341.33) .. controls (196,311.33) and (196,291.33) .. (216,281.33) -- cycle ;

\draw    (80,250) -- (420,250) ;
\draw    (406.33,189.83) -- (576,189.83) ;
\draw    (399,178.5) -- (406,185.25) ;
\draw    (406,185.25) -- (406.5,213.75) ;
\draw    (357,229.67) -- (406.5,213.75) ;
\draw    (357,198.73) -- (406,185.25) ;
\draw    (349.27,191) -- (399,178.5) ;
\draw    (262.6,179.4) -- (399,178.5) ;
\draw    (262.6,179.4) -- (207,191) ;
\draw    (206.6,219) -- (80,250) ;
\draw    (576,189.83) -- (420,250) ;
\draw  [fill={rgb, 255:red, 255; green, 255; blue, 255 }  ,fill opacity=1 ] (141.16,117.69) -- (309.2,184.97) .. controls (311.08,185.72) and (309.99,186.33) .. (306.77,186.33) .. controls (303.55,186.33) and (299.42,185.72) .. (297.54,184.97) -- (129.5,117.69) .. controls (127.63,116.94) and (128.71,116.33) .. (131.93,116.33) .. controls (135.15,116.33) and (139.29,116.94) .. (141.16,117.69) .. controls (143.04,118.44) and (141.95,119.05) .. (138.73,119.05) .. controls (135.51,119.05) and (131.38,118.44) .. (129.5,117.69) ;
\draw  [fill={rgb, 255:red, 208; green, 2; blue, 33 }  ,fill opacity=1 ] (325.71,214.33) .. controls (329.72,213.23) and (333.86,212.33) .. (334.97,212.33) -- (372.52,212.33) .. controls (373.62,212.33) and (371.27,213.23) .. (367.26,214.33) -- (345.47,220.33) .. controls (341.46,221.44) and (337.32,222.33) .. (336.21,222.33) -- (298.66,222.33) .. controls (297.56,222.33) and (299.92,221.44) .. (303.93,220.33) -- cycle ;
\draw  [fill={rgb, 255:red, 255; green, 255; blue, 255 }  ,fill opacity=1 ] (329.36,215.76) .. controls (331.47,215.17) and (333.65,214.7) .. (334.23,214.7) -- (356.05,214.7) .. controls (356.63,214.7) and (355.39,215.17) .. (353.28,215.76) -- (341.82,218.91) .. controls (339.71,219.49) and (337.53,219.96) .. (336.95,219.96) -- (315.13,219.96) .. controls (314.55,219.96) and (315.79,219.49) .. (317.9,218.91) -- cycle ;
\draw   (207,191) .. controls (207,191) and (207,191) .. (207,191) -- (349.5,191) -- (357,198.5) -- (357,229.67) -- (207,229.67) -- cycle ;
\draw  [fill={rgb, 255:red, 2; green, 208; blue, 27 }  ,fill opacity=1 ] (264.71,213.33) .. controls (268.72,212.23) and (272.86,211.33) .. (273.97,211.33) -- (311.52,211.33) .. controls (312.62,211.33) and (310.27,212.23) .. (306.26,213.33) -- (284.47,219.33) .. controls (280.46,220.44) and (276.32,221.33) .. (275.21,221.33) -- (237.66,221.33) .. controls (236.56,221.33) and (238.92,220.44) .. (242.93,219.33) -- cycle ;
\draw  [fill={rgb, 255:red, 255; green, 255; blue, 255 }  ,fill opacity=1 ] (268.36,214.76) .. controls (270.47,214.17) and (272.65,213.7) .. (273.23,213.7) -- (295.05,213.7) .. controls (295.63,213.7) and (294.39,214.17) .. (292.28,214.76) -- (280.82,217.91) .. controls (278.71,218.49) and (276.53,218.96) .. (275.95,218.96) -- (254.13,218.96) .. controls (253.55,218.96) and (254.79,218.49) .. (256.9,217.91) -- cycle ;

\draw  [draw opacity=0][fill={rgb, 255:red, 255; green, 255; blue, 255 }  ,fill opacity=1 ] (125.5,113.25) -- (187,113.25) -- (187,131.5) -- (125.5,131.5) -- cycle ;
\draw  [pattern=_9b2vxhdkv,pattern size=3pt,pattern thickness=0.75pt,pattern radius=0pt, pattern color={rgb, 255:red, 0; green, 0; blue, 0}] (329.36,215.76) .. controls (331.47,215.17) and (333.65,214.7) .. (334.23,214.7) -- (356.05,214.7) .. controls (356.63,214.7) and (355.39,215.17) .. (353.28,215.76) -- (341.82,218.91) .. controls (339.71,219.49) and (337.53,219.96) .. (336.95,219.96) -- (315.13,219.96) .. controls (314.55,219.96) and (315.79,219.49) .. (317.9,218.91) -- cycle ;

\draw    (386.2,154.2) -- (339,213.77) ;
\draw [shift={(339,213)}, rotate = 308.72] [color={rgb, 255:red, 0; green, 0; blue, 0 }  ][line width=0.75]    (10.93,-3.29) .. controls (6.95,-1.4) and (3.31,-0.3) .. (0,0) .. controls (3.31,0.3) and (6.95,1.4) .. (10.93,3.29)   ;
\draw    (346.2,155) -- (300,211.26) ;
\draw [shift={(300,211.5)}, rotate = 309.83] [color={rgb, 255:red, 0; green, 0; blue, 0 }  ][line width=0.75]    (10.93,-3.29) .. controls (6.95,-1.4) and (3.31,-0.3) .. (0,0) .. controls (3.31,0.3) and (6.95,1.4) .. (10.93,3.29)   ;

\draw (451,85.67) node [anchor=north west][inner sep=0.75pt]   [align=left] {\tiny{\underline{Air}}};
\draw (410,116.07) node [anchor=north west][inner sep=0.75pt]    {\tiny{$\mu _{\alpha } = \mu _{0} ,\ \sigma _{\alpha } = 0$}};

\draw (381,285.67) node [anchor=north west][inner sep=0.75pt]   [align=left] {\tiny{\underline{Soil}}};
\draw (340,316.07) node [anchor=north west][inner sep=0.75pt]    {\tiny{$\mu _{\alpha } \approx \mu _{0} ,\ \sigma _{\alpha } \approx 0$}};
\draw (214.87,279.9) node [anchor=north west][inner sep=0.75pt]   [align=left] {\tiny{\underline{Object $\displaystyle B_{\alpha }$}}};
\draw (209,303.73) node [anchor=north west][inner sep=0.75pt]    {\tiny{$ \begin{array}{l}
\mu _{\alpha } =\mu _{*} ,\\
\sigma _{\alpha } =\sigma _{*}
\end{array}$}};
\draw (387.6,141.8) node [anchor=north west][inner sep=0.75pt]    {\tiny{$S$}};
\draw (348.4,140.2) node [anchor=north west][inner sep=0.75pt]    {\tiny{$\boldsymbol{J}_0$}};
%

\end{tikzpicture}} &
\raisebox{-0.5\height}{\input{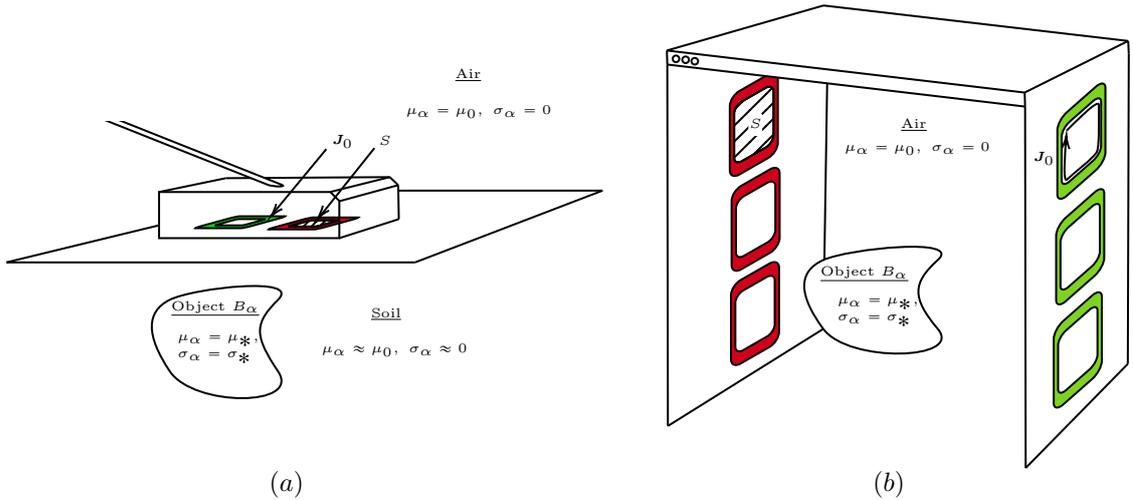}} \\
(a) & (b)
\end{array}$
\caption{Illustrative diagram of a ($a$) handheld and ($b$) walk through metal detector indicating the coil arrays for the transmitting (green) and receiving (red) coils. The object $B_\alpha$ and its material properties are also indicated.}
\label{fig:metal_detector_diagrams}
\end{figure}

When a low frequency sinusoidal time varying current with volume current density ${\boldsymbol J}_0$ is passed through an exciting coil it generates a background (exciting) magnetic field with complex amplitude $\boldsymbol{H}_0$ and the presence of a highly conducting magnetic object, with (possibly inhomogeneous) electrical conductivity $0\ll \sigma_*< \infty $ and magnetic permeability $0< \mu_*< \infty$, leads to the generation of eddy currents that perturb this field. A second set of measurement coils sense this perturbation in the form of an induced voltage of complex amplitude $V$ in each coil. Note that pulsed excitation is also used in metal detection~\cite{ambrus} and has recently been successfully employed for both measurement~\cite{simic2023} and classification between large complex metallic objects and electronic components~\cite{Thomas2024} and smaller landmines~\cite{Simic2023b}.
The voltages contain considerably more information about the size, shape, materials and location of hidden target than the simple {aforementioned} audible alarm may suggest. As illustrated in Figure~\ref{fig:metal_detector_diagrams}, for a walk-through metal detector configuration and a hand-held metal detector, $V$ is  related to the perturbed magnetic field $(\boldsymbol{H}_\alpha- \boldsymbol{H}_0)(\boldsymbol{x})$ at position ${\boldsymbol x}$ by
 \begin{equation}
 V = \im \omega \mu_0  \int_S (\boldsymbol{H}_\alpha- {\boldsymbol{H}}_0)(\boldsymbol{x}) \cdot \boldsymbol{n} \dif \boldsymbol{x},
\label{eqn:voltage}
\end{equation}
where $\im := \sqrt{-1}$, $S$ is a closed surface defined by the geometry of the coil, $\omega$ is the angular frequency of excitation, and ${\boldsymbol n}$ is the unit outward normal with respect to $S$. In the case of a buried object, we assume  that the soil is non-conducting and non-magnetic, this is a reasonable approximation given the high value of $\sigma_*$, but the accuracy of the approximation deteriorates in the case of wet soils. The position dependent electrical conductivity, $\sigma_\alpha$, and magnetic permeability, $\mu_\alpha$, describing the materials in the object, $B_\alpha$, and the surrounding non-conducting region $B_\alpha^c= {\mathbb R}^3 \setminus \overline{B_\alpha}$, where the overbar denotes the closure (later an overbar is also used to denote the complex conjugate, but it should clear from the context as to which applies) are
 \begin{align}
  \sigma_\alpha = \left \{ \begin{array}{cl}
                            \sigma_* & \text{in $B_\alpha$} \nonumber                                       \\
                            0        & \text{in $B_\alpha^c$}
  \end{array} \right ., \qquad
  \mu_\alpha = \left \{ \begin{array}{cc}
                         \mu_* & \text{in $B_\alpha$} \nonumber \\
                         \mu_0 & \text{in $B_\alpha^c $}
  \end{array} \right . .
 \end{align}
It is also convenient to define $\mu_r: = \mu_* /\mu_0$ as the (position dependent) relative magnetic permeability in the object. 


The metal detection inverse problem consists of determining information about the shape, size, material properties, and location of the hidden object given voltage measurements $V$ in the measurement coils and the prescribed ${\boldsymbol J}_0$ in the exciting coils.

In the case of a  highly conducting spherical magnetic object, analytical solutions are available for the perturbation in magnetic field
$(\boldsymbol{H}_\alpha- \boldsymbol{H}_0)(\boldsymbol{x})$ at positions $\boldsymbol{x}$ away from the object for different forms of exciting low-frequency time harmonic fields in the eddy current regime. This includes results for {the case} where $\boldsymbol{H}_0$ is a uniform background magnetic field~\cite{Wait1951} and $\boldsymbol{H}_0$ is the background field due a dipole source~\cite{Wait1960}.
Related analytical and semi-analytical solutions have been obtained for highly  conducting and magnetic spheroids~\cite{barrowes2004} and ellipsoids~\cite{barrowes2008}. More generally, approximate dipole models suggest that {the} field perturbation due to the presence of a highly conducting magnetic object has a magnetic dipole moment that can be expressed in terms of a complex symmetric magnetic polarizbility tensor (MPT) and the background magnetic field at the position of the object~\cite{landaubook}, which offers possibilities for identifying hidden metallic objects using classification. Using this model, MPTs have been approximately computed for thin shells~\cite{Gabbay2019} and for other objects by using commercial codes~\cite{rehim2015}. 
However, the dipole model only provides an approximation to $(\boldsymbol{H}_\alpha- \boldsymbol{H}_0)(\boldsymbol{x})$. Its accuracy can be assessed by  comparing it to a rigorously described asymptotic expansion of $(\boldsymbol{H}_\alpha- \boldsymbol{H}_0)(\boldsymbol{x})$ due to the presence of a small highly conducting magnetic object $B_\alpha$ as its size tends to zero~\cite{Ammari2014} for the eddy current regime and ${\boldsymbol x}$ away from $B_\alpha$. In this expansion, the object $B_\alpha:= \alpha B + \boldsymbol{z}$  is described in terms of a size parameter $\alpha \ll 1$ m, a dimensionless object $B$ having the same shape as $B_\alpha$ and placed at the origin, and a translation $\boldsymbol{z}$. By using tensor manipulations and integration by parts, the expansion has been shown to reduce to~\cite{LedgerLionheart2019}
\begin{equation}
(\boldsymbol{H}_\alpha - \boldsymbol{H}_0)(\boldsymbol{x})_i = (\boldsymbol{D}_x^2 G(\boldsymbol{x},\boldsymbol{z}))_{ij} ({\mathcal M})_{jk} (\boldsymbol{H}_0 (\boldsymbol{z}))_k + (\boldsymbol{R}(\boldsymbol{x}))_i , \label{eqnasymp}
\end{equation}
for $\boldsymbol{x}$ away from $B_\alpha$. In the above, $| \boldsymbol{R}(\boldsymbol{x})| \le C \alpha^4 \| \boldsymbol{H}_0 \|_{W^{2,\infty}(B_\alpha)}$ is a residual, with $C$ being a constant independent of $\alpha$, $G(\boldsymbol{x},\boldsymbol{z}):= 1/ ( 4 \pi | \boldsymbol{x}-\boldsymbol{z}|)$ is the free space Green's Laplace function, {$\boldsymbol{D}_x^2 G(\boldsymbol{x},\boldsymbol{z})$ denotes the Hessian of $G({\boldsymbol x},{\boldsymbol z})$}, ${\mathcal M} = ({\mathcal M})_{ij} {\boldsymbol e}_i \otimes {\boldsymbol e}_j$ is the complex symmetric rank--2 MPT with six complex coefficients, ${\boldsymbol e}_i$ is the $i$th orthonormal unit vector and Einstein summation convention has been applied. The MPT coefficients $({\mathcal M})_{ij} $ are independent of position~\cite{Ammari2015} and 
so the leading order term in~(\ref{eqnasymp}) separates the object position from the characterisation of the object.
{The advantages of the asymptotic expansion {are} that it provides a measure of accuracy of the approximation and has explicit  expressions for computing the MPT coefficients. Furthermore, these explicit expressions hold for  inhomogeneous objects and multiple objects~\cite{LedgerLionheart2018},
where { $ \sigma_*$   and  $ \mu_* $} are isotropic and independent of  $\omega$, but not necessarily homogeneous.

In the case of a single object, substituting (\ref{eqnasymp}) into  (\ref{eqn:voltage}), and considering a single exciter-measurement coil pair, it can be shown that 
\begin{equation}\label{eqn:meas_voltage}
V \approx ({\boldsymbol H}_0^{ms}({\boldsymbol z}) )_i ( {\mathcal M})_{ij} ({\boldsymbol H}_0({\boldsymbol z}))_j, 
\end{equation}
where ${\boldsymbol H}_0^{ms}({\boldsymbol z})$ is the background magnetic field at the position of the object if the measurement coil acts as an exciter~\cite{LedgerLionheart2018}. 
The background magnetic fields
$ {\boldsymbol H}_0^{ms}({\boldsymbol z}) $ and $ {\boldsymbol H}_0({\boldsymbol z})$ do not depend on the shape, size or materials of the object.
 To find ${\boldsymbol z}$, a multiple signal classification (MUSIC) algorithm (e.g.~\cite{Ammari2014}) can be applied, and, once ${\boldsymbol z}$ is known, the coefficients $({\mathcal M})_{ij}$ for the hidden object follow from the solution of a linear least squares problem. 
 The obtained $({\mathcal M})_{ij}$ depend on the object orientation and transform in the standard way for a rank--2 tensor as $({\mathcal M}')_{ij} = ({\bf R})_{ip} ({\mathbf R})_{jq} ({\mathcal M})_{pq} $ where ${\mathbf R}$ is an orthogonal rotation matrix. As the object orientation is unknown, an appropriate feature for classification and dictionary matching is the eigenvalues of the real and imaginary parts of the MPT (or equivalently other tensor invariants) as these are independent of ${\mathbf R}$~\cite{ledgerwilsonamadlion2021}.
 
Exploiting the MPT's spectral signature~\cite{LedgerLionheart2019} (the variation of the MPT coefficients $({\mathcal M}[\omega])_{ij}$ as a function of exciting frequency $\omega$) {provides additional information} compared to an MPT at a single fixed frequency, which only characterises the object's shape, size, and materials up to the best fitting ellipsoid. {Indeed, the measured spectral signature $({\mathcal M}[\omega])_{ij}$ can be determined from $V[\omega]$ in an analogous way as the fixed frequency case by repeating the process for each exciting frequency considered. To identify the target, the dictionary matching is performed using eigenvalues (or tensor invariants) of the real and imaginary parts of ${\mathcal M}[\omega]$ against a dictionary of MPTs spectral signature obtained for other objects.
The prediction of $({\mathcal M}[\omega])_{ij}$ for a known object requires the approximation of an eddy current transmission problem using the finite element method (FEM)~\cite{monkbook} for each $\omega$ of interest, and, to accurately capture the MPT coefficients, a discretisation with small mesh spacing $h$ and elements of high order $p$ is required. To reduce the cost of repeated FEM simulations needed for the spectral signature, 
a reduced order methodology using proper orthogonal decomposition (POD) has been developed~\cite{Wilson2021}. This separates into two stages. Firstly, in the off-line stage, a small number of expensive full FEM simulations at selected $\omega$ snapshots are obtained. Then, in the on-line stage, the MPT spectral signature is rapidly obtained for $\omega_{min} \le \omega\le \omega_{max}$. The approach has been applied to compute dictionaries of non-magnetic objects~\cite{ledgerwilsonamadlion2021} that have, in turn, been used for training machine learning classifiers to identify hidden objects~\cite{ledgerwilsonlion2022}.

Practical applications of MPTs and related technology include at transport hubs using single frequencies~\cite{marsh2013}, spectral signatures in the detection of multiple objects ~\cite{marsh2014}, using machine learning classifiers~\cite{marsh2014b} and improving their reliability~\cite{marsh2015}. MPTs have also been used for identification of landmines~\cite{dekdouk}, landmine surrogates~\cite{abdel}, other UXOs~\cite{barrowes2008}, objects buried in  underwater sediment~\cite{Song2016} and ice~\cite{Wilson2020},  ensuring food safety~\cite{zhao2016}, and scrap metal sorting of shredded metallic scrap~\cite{karimian2017} and vapes~\cite{Williams2024}.
 
While a reduced order modelling approach has been developed for computing MPT spectral signatures, our previous approach still requires large computational demands for complex geometries, which impede its application to the characterisation of complex inhomogeneous magnetic real-world targets. Furthermore, it was not clear how best to choose the number and thicknesses of prismatic boundary layers in FEM discretisations in order to achieve accurate results and faster convergence of the solution, nor was it clear how best to choose the number of full order solution snapshots in order to achieve a reduced order model that is in close agreement to the underlying full order model. This work addresses these important shortcomings through the following novelties:

\begin{itemize}
\item {An improved algorithm  for the numerical computation of the MPT coefficients for varying parameters using a reduced order model POD approach that significantly improves the efficiency and scalability compared to a naive approach which does not account of the special structure of the MPT coefficients and their relationship to the underlying FEM matrices.}

\item {The presentation of a greedy algorithm~\cite{Hesthaven2014} for adaptively choosing new POD snapshot parameters and the performance benefit compared to a non--adaptive strategy for computing the MPT coefficients for varying parameters.}

\item {A simple recipe for designing prismatic boundary layer discretisations in FEM, which, in combination with $p$--refinement, lead to exponential convergence of the MPT coefficients.}

\end{itemize}

The paper is organised as follows: Brief comments on notation are provided in Section~\ref{sect:notation}. Section~\ref{sect:mpt} reviews the explicit formulae for computing the MPT coefficients, which includes, in Section~\ref{sect:efften}, new discrete formulae for computing the MPT coefficients in terms of FEM matrices. Then, Section~\ref{sect:rom}, briefly reviews the off-line and on-line stages of the POD {with projection} (PODP) approach and explains how the calculation of the MPT coefficients can be accelerated by combining the approach described in Section~\ref{sect:efften} with a PODP reduced order description. Section~\ref{sect:rom} also outlines an adaptive algorithm for computation of new snapshot frequencies to improve the accuracy of the POD approach. In Section~\ref{sect:compres}, the computational resources used for the numerical experiments conducted in this work are described, details of the specific versions of the libraries used for simulations are included as well as the details of the open-access software accompanying this work. Then, in Section~\ref{sect:boundarylayers}, recipes for choosing the number and thicknesses of prismatic boundary layer discretisations in order to resolve the thin skin depths associated with highly conducting and highly magnetic objects are described and compared. Section~\ref{sect:numexp} presents a series of  examples to demonstrate the success of the proposed methodologies and the paper closes with some concluding remarks in Section~\ref{sect:Conclusions}.

\section{Notation} \label{sect:notation}
{Calligraphic symbols e.g. ${\mathcal M}$ are used for rank--2 tensors and their coefficients are denoted by $({\mathcal M})_{ij}$. The $i$th orthonormal basis vector is ${\boldsymbol e}_i$  and bold face italics, e.g. ${\boldsymbol H}$,  are used for vector fields. Components of vector fields are denoted by $({\boldsymbol H})_i = {\boldsymbol e}_i \cdot {\boldsymbol H}$, which should be distinguished from indexed vectors e.g. ${\boldsymbol e}_i$. Bold face upper case Roman font is used for linear algebra matrices e.g. ${\mathbf A}$ and bold face lower case Roman font for linear algebra vectors e.g. ${\mathbf b}$, and their coefficients are denoted by ${(\mathbf A})_{ij}$ and $({\mathbf b})_i$, respectively. Finally, position dependence is sometimes highlighted using parentheses e.g. $\boldsymbol{H}_0(\boldsymbol{z})$ and parameter dependence using square brackets, e.g. $\mathbf{A}[\boldsymbol{\omega}]$.}

\section{Magnetic Polarizability Tensor Object Characterisation} \label{sect:mpt}

{Following~\cite{LedgerLionheart2019}, in order to compute the MPT coefficents, the auxiliary vector field  ${\boldsymbol \theta}_i^{(1)}\in {\mathbb C}^3$ is introduced, for $i=1,2,3$, as the solution of the vectorial transmission problem
\begin{subequations}\label{eqn:theta1trans}
\begin{align}
\nabla \times \mu_r^{-1} \nabla \times {\boldsymbol \theta}_i^{(1)} -\im \nu {\boldsymbol \theta}_i^{(1)} &= \im \nu {\boldsymbol \theta}_i^{(0)} && \text{in $B$}, \\
\nabla \times\nabla \times {\boldsymbol \theta}_i^{(1)}  &=  {\boldsymbol 0} && \text{in $B^c:= {\mathbb R}^3 \setminus \overline{B}$} ,\\
\nabla \cdot {\boldsymbol \theta}_i^{(1)} & = 0 && \text{in $B^c$}, \\
[{\boldsymbol n} \times {\boldsymbol \theta}_i^{(1)}]_\Gamma ={\boldsymbol 0}, \ [{\boldsymbol n} \times \tilde{\mu}_r^{-1} \nabla \times {\boldsymbol \theta}_i^{(1)}]_\Gamma &  ={\boldsymbol 0} && \text{on $\Gamma:=\partial B$},  \\
{\boldsymbol \theta}_i^{(1)} & = O ( |{\boldsymbol \xi}|^{-1} ) && \text{as $| {\boldsymbol \xi} | \to \infty$},
\end{align}
\end{subequations}
where ${\boldsymbol \xi}$ is measured from the origin, which lies inside $B$, $[\cdot ]_\Gamma$ denotes the jump over $\Gamma$, ${\boldsymbol n}$ is the unit outward normal, $\nu:= \alpha^2  \mu_0 \omega \sigma_*$, 
 $\tilde{\mu}_r ({\boldsymbol \xi}) = \mu_r({\boldsymbol \xi})$ for ${\boldsymbol \xi}\in B$ and $\tilde{\mu}_r ({\boldsymbol \xi})=1 $ for ${\boldsymbol \xi} \in B^c$.
Furthermore, $ {\boldsymbol \theta}_i^{(0)}\in {\mathbb R}^3$ is the solution of a related transmission problem with $\nu=0$ in $B$ and $\tilde{\boldsymbol \theta}_i^{(0)} := {\boldsymbol \theta}_i^{(0)} - {\boldsymbol e}_i \times {\boldsymbol \xi} \in {\mathbb R}^3$.  Given the solution to these problems,
 the six coefficients $({\mathcal M})_{ij} $ of the complex symmetric MPT ${\mathcal M} =({\mathcal M})_{ij} {\boldsymbol e}_i \otimes {\boldsymbol e}_j$ are obtained by post-processing using the additive decomposition $(\mathcal{M})_{ij}:=(\tilde{\mathcal{R}})_{ij}+\im(\mathcal{I})_{ij}=(\mathcal{N}^0)_{ij}+(\mathcal{R})_{ij}+\im(  \mathcal{I})_{ij}$ where}
\begin{subequations}
\label{eqn:NRI}
\begin{align}
(\mathcal{N}^0[ \alpha B,\mu_r] )_{ij}&:=\alpha^3\delta_{ij}\int_{B}(1-\tilde{\mu}_r^{-1})\dif \boldsymbol{\xi}+\frac{\alpha^3}{4}\int_{B\cup B^c}\tilde{\mu}_r^{-1}\nabla\times\tilde{\boldsymbol{\theta}}_i^{(0)}\cdot\nabla\times\tilde{\boldsymbol{\theta}}_j^{(0)}\dif \boldsymbol{\xi},\\
(\mathcal{R}[\alpha B, \omega,\sigma_*,\mu_r])_{ij}&:=-\frac{\alpha^3}{4}\int_{B\cup B^c}\tilde{\mu}_r^{-1}\nabla\times\overline{\boldsymbol{\theta}_i^{(1)}}\cdot\nabla\times{\boldsymbol{\theta}_j^{(1)}}\dif \boldsymbol{\xi}, \label{eqn:Rtensor}\\
(\mathcal{I}[\alpha B, \omega,\sigma_*,\mu_r])_{ij}&:=\frac{\alpha^3}{4}\int_B\nu\Big(\overline{\boldsymbol{\theta}_i^{(1)}+\tilde{\boldsymbol{\theta}}_i^{(0)}+\boldsymbol{e}_i\times\boldsymbol{\xi}}\Big)\cdot\Big({\boldsymbol{\theta}_j^{(1)}+\tilde{\boldsymbol{\theta}}_j^{(0)}+\boldsymbol{e}_j\times\boldsymbol{\xi}}\Big)\dif \boldsymbol{\xi}, \label{eqn:Itensor}
\end{align}
\end{subequations}
are each the coefficients of real symmetric {rank--2} tensors and $\delta_{ij}$ is the Kronecker delta}.

The problem (\ref{eqn:theta1trans}) is set on an unbounded domain and, for the purposes of approximate computation, it is replaced by a problem on a bounded domain $\Omega$ with truncation in the form of  a convex outer boundary $\partial \Omega$ placed sufficiently far from the object of interest $B$ {with} ${\boldsymbol n} \times {\boldsymbol \theta}_i^{(1)}={\boldsymbol 0}$ applied on $\partial \Omega$ as an approximation to (\ref{eqn:theta1trans}e). {In addition, the Coulomb gauge (\ref{eqn:theta1trans}c) is circumvented by regularisation~\cite{ledgerzaglmayr2010}, where  $\varepsilon$ is a small regularisation parameter, leading to} 
\begin{subequations}\label{eqn:theta1transtrun}
\begin{align}
\nabla \times \mu_r^{-1} \nabla \times {\boldsymbol \theta}_i^{(1)} -\im \nu {\boldsymbol \theta}_i^{(1)} &= \im \nu {\boldsymbol \theta}_i^{(0)} && \text{in $B$}, \\
\nabla \times\nabla \times {\boldsymbol \theta}_i^{(1)} + \varepsilon {\boldsymbol \theta}_i^{(1)}  &=  {\boldsymbol 0} && \text{in $\Omega \setminus \overline{B}$} ,\\
[{\boldsymbol n} \times {\boldsymbol \theta}_i^{(1)}]_\Gamma ={\boldsymbol 0}, \ [{\boldsymbol n} \times \tilde{\mu}_r^{-1} \nabla \times {\boldsymbol \theta}_i^{(1)}]_\Gamma &  ={\boldsymbol 0} && \text{on $\Gamma$},  \\
{\boldsymbol n} \times {\boldsymbol \theta}_i^{(1)} & = {\boldsymbol 0} && \text{on $\partial \Omega$}.
\end{align}
\end{subequations}

\subsection{Finite Element Discretisation}
Following the approach in~\cite{Wilson2021}, the weak solution: find $ {\boldsymbol \theta}_i^{(1)} \in Y^\epsilon $ such that
\begin{align}
a({\boldsymbol \theta}_i^{(1)}, {\boldsymbol w}) := k ({\boldsymbol \theta}_i^{(1)}, {\boldsymbol w})  - \im \omega c({\boldsymbol \theta}_i^{(1)}, {\boldsymbol w}) + \epsilon m({\boldsymbol \theta}_i^{(1)}, {\boldsymbol w})= r ({\boldsymbol w}) \qquad \forall {\boldsymbol w}\in Y^\epsilon\nonumber
\end{align}
where  $Y^\epsilon := \{  {\boldsymbol u} \in {\boldsymbol H}(\text{curl}) : {\boldsymbol n} \times {\boldsymbol u} = {\boldsymbol 0} \text{ on $\partial \Omega$} \}$ and
\begin{align*}
k({\boldsymbol u}, {\boldsymbol v}) : = \int_\Omega \tilde{\mu}_r^{-1} \nabla \times {\boldsymbol u} \cdot  \nabla \times \overline{\boldsymbol v} \, \dif {\boldsymbol \xi} , \qquad
c({\boldsymbol u}, {\boldsymbol v}) : = \alpha^2 \mu_0  \int_B \sigma_* {\boldsymbol u} \cdot  \overline{\boldsymbol v} \, \dif {\boldsymbol \xi} , \qquad
m({\boldsymbol u}, {\boldsymbol v}) : =   \int_{\Omega \setminus \overline{B}} {\boldsymbol u} \cdot  \overline{\boldsymbol v} \, \dif {\boldsymbol \xi} ,
\end{align*}
is approximated by a $hp$-version ${\boldsymbol H}(\hbox{curl})$ conforming FEM approximation.
The discrete approximation  on unstructured tetrahedral meshes of variable size $h$  with elements of uniform order $p$ takes  the form}
\begin{align}
\boldsymbol{\theta}_i^{(1,hp)}(\boldsymbol{\xi})[\boldsymbol{\omega}]&\vcentcolon=\sum_{n=1}^{N_d}\boldsymbol{N}^{(n)}(\boldsymbol{\xi})({\mathbf q}_{i}[\boldsymbol{\omega}])_n = {\mathbf N}\boldsymbol{q}_i[\mathbf{\omega}]   \in {Y^\epsilon \cap W^{(hp)}},  \label{eqn:theta1apprx}
\end{align}
where  $W^{(hp)} \subset {\boldsymbol H}(\text{curl})$ is an appropriate set of ${\boldsymbol H}(\text{curl})$ conforming FEM basis functions, ${\boldsymbol N}^{(n)}$ is a typical basis function with $\mathbf{N}:=[\boldsymbol{N}^{(1)}, \boldsymbol{N}^{(2)},\cdots, \boldsymbol{N}^{(N_d)}]\in \mathbb{R}^{N_d \times N_d}$ , $({\mathbf q}_{i}[\boldsymbol{\omega}])_n$ is the $n$th solution coefficient for the $i$th solution direction, which depends on parameters of interest ${\boldsymbol \omega} = \left[\alpha B, \omega, \sigma_*, \mu_r \right]$,  and $N_d$ {is} the number of degrees of freedom. 
  The FEM approximation for the $i$th direction then corresponds to the solution of the  linear system of equations 
\begin{equation}\label{eqn:Linear}
\mathbf{A}[\boldsymbol{\omega}]\mathbf{q}_i[\boldsymbol{\omega}]=\mathbf{r}(\boldsymbol{\theta}_i^{(0,hp)})[ \boldsymbol{\omega}],
\end{equation}
for $\mathbf{q}_i[\boldsymbol{\omega}]\in {\mathbb C}^{N_d}$ where ${\mathbf r}({\boldsymbol \theta}_i^{(0,hp)})[ \boldsymbol{\omega}]\in{\mathbb C}^{N_d}$ is a known source term~\cite{Wilson2021}[Equation (17)]. In the above, {${\mathbf A}[\boldsymbol{\omega}]= {\mathbf K}[ \mu_r ] - \im \omega {\mathbf C}[\tilde{\nu}] +\varepsilon{\mathbf M}  \in{\mathbb C}^{N_d\times N_d}$}  is a large parameter dependent complex  symmetric  sparse matrix, where 
\begin{subequations}
\begin{align}
({\mathbf K}\left[\mu_r\right])_{mn}:=& k({\boldsymbol N}^{(m)}, {\boldsymbol N}^{(n)}) , \label{eqn:kmatrix} \\
  ({\mathbf C}\left[\tilde{\nu}\right])_{mn}:=&({\mathbf C}\left[\alpha^2 \mu_0 \sigma_*\right])_{mn}=c({\boldsymbol N}^{(m)}, {\boldsymbol N}^{(n)}) ,
\label{eqn:cmatrix}\\
 ({\mathbf M})_{mn}:=& m({\boldsymbol N}^{(m)}, {\boldsymbol N}^{(n)}) . \end{align}
\end{subequations}
In the following, the parameter dependence for these matrices is omitted and ${\mathbf K} \in {\mathbb R}^{N_d \times N_d}$  is identified as a curl-curl stiffness matrix, ${\mathbf C} \in {\mathbb R}^{N_d \times N_d}$ as a damping matrix associated with the conducting region $B$, and ${\mathbf M} \in {\mathbb R}^{N_d \times N_d}$ as a mass matrix  associated with regularisation in $\Omega \setminus \overline{B}$.

Unlike our earlier work, FEM discretisations with unstructured tetrahedral elements with prismatic layers, where appropriate, are included in the present contribution. These prismatic boundary layers are used, in conjunction with $p$--refinement, to model thin skin depths, which is discussed further in Section~\ref{sect:boundarylayers}. Following~\cite{cairlet}, and noting the convention that $p=0$ elements have constant tangential components on edges, but consist of vector valued linear basis functions, 
the order of the Gaussian quadrature used for computing $({\mathbf K})_{mn}$, 
 $({\mathbf C})_{mn}$ and $({\mathbf M})_{mn}$ is chosen to be sufficiently high so that it can integrate products of appropriate polynomials exactly. The approximation error associated with under-integration of the non-linear rational functions arising from the interpolating polynomial functions of degree $g$ used to approximate the geometry of curved boundary and transmission faces is absorbed into other approximation errors in the FEM discretisation. The linear system (\ref{eqn:Linear}) is solved to a relative tolerance $TOL$ using a conjugate gradient solver and a balancing domain decomposition by constraints (BDDC) preconditioner~\cite{Dohrmann2003}.
 {While the order of Gaussian quadrature for computing $({\mathbf K})_{mn}$, 
 $({\mathbf C})_{mn}$ and $({\mathbf M})_{mn}$ can be determined directly from $p$,  the same is not true when computing the MPT coefficients by post-processing as the approximation errors can no longer be absorbed. Hence, in this case, attention must also be paid to using a Gaussian quadrature scheme of sufficient order, not only to account for the order of the FEM approximation, but also  for the curvature of elements placed next to  a curved boundary or interface.}
 
A first approach to computing MPT coefficients is to apply Gaussian quadrature to (\ref{eqn:Rtensor}) and (\ref{eqn:Itensor}) with $\boldsymbol{\theta}_i^{(1)}[ {\bm \omega} ] \approx \boldsymbol{\theta}_i^{(1,hp)} [{\bm \omega}] \approx {\mathbf N} {\mathbf q}_i[{\bm \omega}] $. We call this approach the Integral Method (IM) and this requires  $O(N_{elem}N_{nip} N_d^2 )$ operations to compute the MPT coefficients for each parameter set ${\bm \omega}$ once (\ref{eqn:Linear}) has been solved. Here, $N_{elem}$ and $N_{nip}$ are the number of elements and number of integration points, respectively. However, for a large problems, with many output parameter sets, this is costly and scales poorly. We address this in the following sections.

\subsection{Improved Efficiency for the Calculation of the MPT Spectral Signature} \label{sect:efften}

{This section focuses on the improved efficiency for the computation of $(\mathcal{R}[\alpha B, \omega,\sigma_*,\mu_r])_{ij}$ and \\$(\mathcal{I}[\alpha B,\omega,\sigma_*,\mu_r])_{ij}$, and, in particular, for the case where $\omega$ is varied to produce the MPT spectral signature.} Note that a similar approach could also be applied to variations of other parameters of interest. For $(\mathcal{R}[\alpha B, \omega,\sigma_*,\mu_r])_{ij}$, using (\ref{eqn:Rtensor}), (\ref{eqn:theta1apprx}) {and (\ref{eqn:kmatrix}),} {it follows that}
\begin{align}
(\mathcal{R}[\alpha B, \omega,\sigma_*,\mu_r])_{ij} = & - \frac{\alpha^3}{4}  \sum_{m=1}^{N_d} \sum_{n=1}^{N_d} \overline{({\mathbf q}_i)_m } \int_{\Omega} \tilde{\mu}_r^{-1} \nabla \times {\boldsymbol N}^{(m)} \cdot 
 \nabla \times {\boldsymbol N}^{(n)} \dif {\boldsymbol \xi}   {({\mathbf q}_{j})_n } \nonumber \\
 = & - \frac{\alpha^3}{4} \overline{{\mathbf q}_i}^T {\mathbf K} {\mathbf q}_j , \label{eqn:Rmatrixmeth}
\end{align}
where {the coefficients of ${\mathbf K}$  are  independent of $\omega$ and $T$ denotes the transpose}. 

Next, for $(\mathcal{I}[\alpha B, \omega,\sigma_*,\mu_r])_{ij}$, expanding (\ref{eqn:Itensor}), it follows that
\begin{align}
(\mathcal{I}[\alpha B, \omega,\sigma_*,\mu_r])_{ij} = &\frac{\alpha^3}{4} \left ( \int_B \nu {\boldsymbol \theta}_j^{(1)} \cdot \overline{{\boldsymbol \theta}_i ^{(1)}} \dif {\boldsymbol \xi} +  \int_B \nu \tilde{\boldsymbol \theta}_j^{(0)}\cdot \tilde{\boldsymbol \theta}_i^{(0)} \dif {\boldsymbol \xi} + \int_B \nu {\boldsymbol e}_i \times {\boldsymbol \xi} \cdot {\boldsymbol e}_j \times {\boldsymbol \xi} \dif {\boldsymbol \xi} +  \int_B \nu {\boldsymbol \theta}_j^{(1)} \cdot \tilde{\boldsymbol \theta}_i^{(0)} \dif {\boldsymbol \xi} \right . \nonumber \\
& \left . 
 + \int_B \nu {\boldsymbol \theta}_j^{(1)} \cdot  {\boldsymbol e}_ i \times {\boldsymbol \xi} \dif {\boldsymbol \xi} + \int_B \nu  \tilde{\boldsymbol \theta}_j^{(0)} \cdot \overline{{\boldsymbol \theta}_i^{(1)}} \dif {\boldsymbol \xi} +     \int_B \nu \tilde{\boldsymbol \theta}_j^{(0)} \cdot {\boldsymbol e}_i \times {\boldsymbol \xi} \dif {\boldsymbol \xi} \right . \nonumber \\
& \left .
+ \int_B \nu {\boldsymbol e}_j \times {\boldsymbol \xi} \cdot   \overline{{\boldsymbol \theta}_i^{(1)}}  \dif {\boldsymbol \xi} + \int_B \nu {\boldsymbol e}_j \times {\boldsymbol \xi} \cdot  \tilde{\boldsymbol \theta}_i^{(0)} \dif {\boldsymbol \xi}
\right ), \nonumber
\end{align}
and, since $(\mathcal{I}[\alpha B, \omega,\sigma_*,\mu_r])_{ij} = (\mathcal{I}[\alpha B, \omega,\sigma_*,\mu_r])_{ji} \in {\mathbb R}$, then

 {
\begin{align}
(\mathcal{I}[\alpha B, \omega,\sigma_*,\mu_r])_{ij} = &\frac{\alpha^3}{4} \left ( \int_B \nu \tilde{\boldsymbol \theta}_i^{(0)}\cdot \tilde{\boldsymbol \theta}_j^{(0)} \dif {\boldsymbol \xi} + \int_B \nu {\boldsymbol e}_i \times {\boldsymbol \xi} \cdot {\boldsymbol e}_j \times {\boldsymbol \xi} \dif {\boldsymbol \xi} +  \int_B \nu \tilde{\boldsymbol \theta}_j^{(0)} \cdot {\boldsymbol e}_i \times {\boldsymbol \xi} \dif {\boldsymbol \xi} \right . \nonumber \\
& +\int_B \nu \tilde{\boldsymbol \theta}_i^{(0)} \cdot {\boldsymbol e}_j \times {\boldsymbol \xi} \dif {\boldsymbol \xi}+ \text{Re} \left ( 
 \int_B \nu {\boldsymbol \theta}_j^{(1)} \cdot \overline{{\boldsymbol \theta}_i ^{(1)}} \dif {\boldsymbol \xi}  +   \int_B \nu {\boldsymbol \theta}_j^{(1)} \cdot \tilde{\boldsymbol \theta}_i^{(0)} \dif {\boldsymbol \xi} \right .  \nonumber  \\
& \left .  \left . + \int_B \nu \overline{{\boldsymbol \theta}_i^{(1)}} \cdot \tilde{\boldsymbol \theta}_j^{(0)} \dif {\boldsymbol \xi}
 + \int_B \nu {\boldsymbol \theta}_j^{(1)} \cdot  {\boldsymbol e}_ i \times {\boldsymbol \xi} \dif {\boldsymbol \xi} + \int_B \nu \overline{{\boldsymbol \theta}_i^{(1)}} \cdot  {\boldsymbol e}_ j \times {\boldsymbol \xi} \dif {\boldsymbol \xi}  \right )
\right ) \nonumber .
\end{align}}

Writing {$\tilde{\boldsymbol{\theta}}_i^{(0,hp)}(\boldsymbol{\xi})[\boldsymbol{\omega}] \vcentcolon=\tilde{\boldsymbol{\theta}}_i^{(0,hp)}(\boldsymbol{\xi}) [\mu_r] \vcentcolon=\sum_{m=1}^{M_d}\tilde{\boldsymbol{N}}^{(m)}(\boldsymbol{\xi})   ( \mathbf{o}_{i}[ \mu_r ] )_m $, where $( \mathbf{o}_{i}[\mu_r] )_m$ is the $m$th solution coefficient for the $i$th direction  in a similar way to (\ref{eqn:theta1apprx}),} and following a similar procedure to above{,} then
 {\begin{align}\label{eqn:I_tensor_coeffs}
(\mathcal{I}[\alpha B, \omega,\sigma_*,\mu_r])_{ij} = \frac{\omega \alpha^3}{4} &\left ( {\mathbf o}_i^T {\mathbf C}^{(1)} {\mathbf o}_j + {\rm c}_{ij} + {\mathbf s}_i^T {\mathbf o}_j +
 {\mathbf s}_j^T {\mathbf o}_i +
 \text{Re} \left (
\overline{\mathbf q}_i^T {\mathbf C} {\mathbf q}_j +  {\mathbf o}_i^T {\mathbf C}^{(2)} {\mathbf q}_j+
 {\mathbf o}_j^T {\mathbf C}^{(2)} \overline{\mathbf q}_i \right . \right . \nonumber \\
 & \left . \left .  +
 {\mathbf t}_i^T {\mathbf q}_j +  {\mathbf t}_j^T \overline{\mathbf q}_i 
\right ) \right ). 
\end{align}}
{Recall that ${\mathbf C} \in {\mathbb R}^{N_d \times N_d}$ is independent of $\omega$, and its coefficients are defined in (\ref{eqn:cmatrix}).} {The other terms are defined as}
\begin{align}
{\rm c}_{ij} :=& \int_B \tilde{\nu} {\boldsymbol e}_i \times {\boldsymbol \xi} \cdot  {\boldsymbol e}_j \times {\boldsymbol \xi} \dif  {\boldsymbol \xi}, \qquad
({\mathbf C}^{(1)})_{mn} := \int_B \tilde{\nu} \tilde{\boldsymbol N}^{(m)} \cdot  \tilde{\boldsymbol N}^{(n )} \dif {\boldsymbol \xi}, \qquad ({\mathbf C}^{(2)})_{mn } := \int_B \tilde{\nu} \tilde{\boldsymbol N}^{(m)} \cdot  {\boldsymbol N}^{(n)} \dif {\boldsymbol \xi}, \nonumber \\
({\mathbf s}_i)_n &:= \int_B \tilde{\nu} {\boldsymbol e}_i \times {\boldsymbol \xi} \cdot \tilde{\boldsymbol N}^{(n)} \dif {\boldsymbol \xi}, \qquad ({\mathbf t}_i)_n: = \int_B \tilde{\nu} {\boldsymbol e}_i \times {\boldsymbol \xi} \cdot {\boldsymbol N}^{(n)} \dif {\boldsymbol \xi}, \nonumber 
\end{align}
{and have dimensions}
 ${\mathbf C}^{(1)} \in {\mathbb R}^{{M}_d \times {M_d}}$, ${\mathbf C}^{(2)} \in {\mathbb R}^{{M}_d \times {N_d}}$, ${\mathbf s}_i \in {\mathbb R}^{M_d}$, ${\mathbf t}_i \in {\mathbb R}^{N_d}$ and ${\rm c}_{ij} \in {\mathbb R}$, respectively, and are also independent of $\omega$. Hence, the computation  of $(\mathcal{R}[\alpha B, \omega,\sigma_*,\mu_r])_{ij}$ and $(\mathcal{I}[\alpha B, \omega,\sigma_*,\mu_r])_{ij}$ has been reduced to matrix vector products, where the matrices can be precomputed and stored as they are all independent of $\omega$.

\begin{remark} \label{rem:basisfun}
The  number of matrices required for computing (\ref{eqn:I_tensor_coeffs}) can be reduced from six to three by using the same discretisation for both the $\tilde{\boldsymbol{\theta}}^{(0)}_i$ and $\boldsymbol{\theta}^{(1)}_i$ problems. In this case, $M_d=N_d$, $\tilde{\boldsymbol N}^{(n)} = {{\boldsymbol N}}^{(n)}$ and $\mathbf{C} = \mathbf{C}^{(1)} = \mathbf{C}^{(2)} \in {\mathbb R}^{N_d \times N_d}$ and ${\mathbf t}_i = {\mathbf s}_i \in {\mathbb R}^{N_d}$. However, vector basis functions that are gradients of scalar polynomials~\cite{SchoberlZaglmayr2005} must be removed via post-processing \cite[pg 145-150]{zaglmayrphd} for the $\tilde{\boldsymbol \theta}_i^{(0)}$ problem.
\end{remark}

Once, the solutions ${\mathbf q}_i[{\bm \omega}] $ to (\ref{eqn:Linear}) have been obtained, and aforementioned matrices and vectors precomputed,  then the computation of the MPT coefficients in this case requires $O(nz N_d)$ operations for each ${\bm \omega}$ where $nz$ is the number of non-zeros in ${\mathbf K}, {\mathbf C}, {\mathbf C}^{(1)}$ and ${\mathbf C}^{(2)}$.  We call this approach the Full Matrix Method (FMM).

\section{PODP Approach} \label{sect:rom}
To further accelerate the computation of the MPT spectral signature, a POD approach is employed. The off-line and on-line stages of projected POD (known as PODP) are reviewed, and the details from Section~\ref{sect:efften} are combined with a PODP reduced order description.
 Then, an a-posteriori error estimate is  used to drive an adaptive procedure for improving the accuracy of the approximated MPT spectral signature.

\subsection{Off-line Stage} \label{sect:romoffline}

Following the solution of (\ref{eqn:Linear}) for   {$\mathbf{q}_i[\boldsymbol{\omega}]$   for different sets  of  snapshot parameters}, $\boldsymbol{\omega}$, {the matrices} $\mathbf{D}_i\in\mathbb{C}^{N_d\times N}$, $i=1,2,3$, whose columns are vectors  of solution coefficients, in the form
\begin{equation}\label{D}
\mathbf{D}_i \vcentcolon=\big[\mathbf{q}_i[ \boldsymbol{\omega}_1],\mathbf{q}_i[\boldsymbol{\omega}_2],...,\mathbf{q}_i[\boldsymbol{\omega}_{N}]\big],
\end{equation}
are constructed, where  $N\ll N_d$ denotes the number of such snapshots. {In the above, it is assumed that the same FEM discretisation is employed for all snapshots, which offers computational savings over repeatably computing the matrices ${\mathbf K}$, ${\mathbf C}$ and ${\mathbf M}$. This also ensures high levels of accuracy are maintained compared to using different discretisations for each snapshot and using interpolation to the finest grid.}
 Application of a singular value decomposition (SVD)~e.g.~\cite{bjorck} gives
\begin{equation}\label{eq:SVD}
\mathbf{D}_i=\mathbf{U}_i \mathbf{\Sigma}_i \mathbf{V}_i^{{H}},\end{equation}
where $\mathbf{U}_i\in\mathbb{C}^{N_d\times N_d}$ and $\mathbf{V}_i \in\mathbb{C}^{N\times N}$ are unitary matrices and $\mathbf{\Sigma}_i\in\mathbb{R}^{N_d\times N}$ is a diagonal matrix  {padded by zeros such} that it becomes rectangular. In the above, $\mathbf{V}_i^{{H}}= \overline{\mathbf{V}_i}^T$ is the Hermitian of $\mathbf{V}_i$.

The diagonal entries $(\mathbf{\Sigma}_i)_{nn}=s_n$ are the singular values of $\mathbf{D}_i$ and they are arranged as $s_1>s_2>...>s_{N}$. {The singular values $s_n$ decay rapidly towards zero,} motivating the introduction of the truncated SVD (TSVD)~e.g.~\cite{bjorck}
\begin{equation}\label{eq:truncSVD}
\mathbf{D}_i\approx
\mathbf{D}_i^M = \mathbf{U}_i^M\mathbf{\Sigma}_i^M\left(\mathbf{V}_i^M\right)^{{H}},
\end{equation}
where $\mathbf{U}_i^M\in\mathbb{C}^{N_d\times M} $ are the first $M$ columns of $\mathbf{U}_i$, $\mathbf{\Sigma}_i^M\in{\mathbb R}^{M\times M}$ is a diagonal matrix containing the first $M$ singular values and $(\mathbf{V}_i^M)^{{H}} \in{\mathbb C}^{M\times N}$ are the first $M$ rows of $\mathbf{V}_i^{{H}}$. Furthermore, $M$ follows from retaining singular values $s_1,\ldots,s_M$ where $s_M$ is the largest singular value such that $s_M /s_1 < TOL_\Sigma$.
\subsection{On-line Stage} \label{sect:podp}
 In the on-line stage of PODP, {$\mathbf{q}_i^{PODP} [ {\boldsymbol \omega}] \approx\mathbf{q}_i[ {\boldsymbol \omega}]$  is obtained as a linear combination of the columns of ${\mathbf U}_i^M$, where the coefficients of this projection are contained in the vector ${\mathbf p}_i^M [ \boldsymbol{\omega}] \in {\mathbb C}^M$, so that
 \begin{align}
\boldsymbol{\theta}_i^{(1,hp)}\left(\boldsymbol{\xi} \right )[\boldsymbol{\omega}]  \approx\left(\boldsymbol{\theta}_i^{(1,hp)}\right)^{\text{PODP}} ({\boldsymbol \xi})[ \boldsymbol{\omega]} :=& \textbf{N}(\boldsymbol{\xi})
\mathbf{q}_i^{PODP} [ {\boldsymbol \omega}] = 
\mathbf{N}(\boldsymbol{\xi}) \mathbf{U}_i^M\mathbf{p}_i^M[ \boldsymbol{\omega}]   \in Y^{(PODP)} , \label{eqn:solution:rebreakdown}
\end{align}
with} $Y^{(PODP)}\subset Y^\epsilon \cap W^{(hp)}$~\cite{Wilson2021}. To obtain ${\mathbf p}_i^M [ \boldsymbol{\omega}] $, {the reduced linear system}
\begin{equation}\label{eqn:ReducedA}
\mathbf{A}_i^M[ \boldsymbol{\omega}]\mathbf{p}_i^M[\boldsymbol{\omega}]=\mathbf{r}_i^M\left(\boldsymbol{\theta}_i^{(0,hp)}\right )[  \boldsymbol{\omega}],
\end{equation}
is solved. This system is obtained by a Galerkin projection, where $\mathbf{A}_i^M[\boldsymbol{\omega}]\vcentcolon=(\mathbf{U}_i^M)^{{H}} \mathbf{A} [\boldsymbol{\omega}]\mathbf{U}_i^M\in {\mathbb C}^{M\times M}$ and $\mathbf{r}_i^M(\boldsymbol{\theta}^{(0,hp)} )[ \boldsymbol{\omega}] \vcentcolon=(\mathbf{U}_i^M)^{{H}} \mathbf{r} ( \boldsymbol{\theta}_i^{(0,hp)})[ \boldsymbol{\omega}] \in {\mathbb C}^M$. Since $M \le N \ll N_d$, the size of (\ref{eqn:ReducedA}) is smaller than (\ref{eqn:Linear}) and, therefore, substantially computationally cheaper to solve. After solving this reduced system, and obtaining $\mathbf{p}_i^M[ \boldsymbol{\omega}]$, an approximate solution for $\boldsymbol{\theta}_i^{(1,hp)}(\boldsymbol{\xi})[ \boldsymbol{\omega}]$ {is obtained} using (\ref{eqn:solution:rebreakdown}). Moreover, the matrix  $\mathbf{A}_i^M[\boldsymbol{\omega}]$ and right hand side $\mathbf{r}_i^M(\boldsymbol{\theta}_i^{(0,hp)})[ \boldsymbol{\omega}]$ can be computed efficiently for new ${\boldsymbol \omega}$, which leads to further {reductions in computational cost for varied parameters}.

\subsection{Improved Computational Efficiency for the PODP Prediction of the MPT Spectral Signature} \label{sect:podpefficent}
{An approach for improving the computational efficiency of the PODP prediction of the MPT spectral signature, for the case where ${\boldsymbol \omega} = \omega$ is varied and $\alpha B$, $\sigma_*$ and $\mu_r$ are fixed, is now described. Since $\mathbf{q}_i^{PODP} [ { \omega}]  = \mathbf{U}_i^M\mathbf{p}_i^M[ {\omega}]  $,} the PODP prediction of the MPT coefficients {for new $\omega$ is}
\begin{align}
 \left(\mathcal{R}^{PODP} [\alpha B, \omega,\sigma_*,\mu_r]\right)_{ij} = & - \frac{\alpha^3}{4} {\overline{{\mathbf p}_i^M}}^T {\mathbf K}_{ij}^M {\mathbf p}_j^M , \label{eqn:podprealmpt}
\end{align}
where ${\mathbf K}_{ij}^M :=( {\mathbf U}_i^M)^H {\mathbf K} {\mathbf U}_j^M \in {\mathbb C}^{M \times M}$ can be precomputed once {for all output $\omega$ of interest. Thus,} for each new value of $\omega$, the coefficients $(\mathcal{R}^{PODP} [\alpha B, \omega,\sigma_*,\mu_r ])_{ij}$ can be obtained from vector--matrix--vector products of small dimension $M$.

Similarly,
\begin{align}\label{eqn:podpimagmpt}
(\mathcal{I}^{PODP} [\alpha B, \omega,\sigma_*,\mu_r])_{ij} = &\frac{\omega \alpha^3}{4} \left ( {\mathbf o}_i^T {\mathbf C}^{(1)} {\mathbf o}_j + {\rm c}_{ij} +
{\mathbf s}_i^T {\mathbf o}_j +
 {\mathbf s}_j^T {\mathbf o}_i+
  \text{Re} \left (
{\overline{\mathbf p_i^M}}^T {\mathbf C}_{{ij}}^M {\mathbf p}_j^M +  {\mathbf o}_i^T {\mathbf C}_{{j}}^{(2),M} {\mathbf p}_j^M+\right . \right . \nonumber \\
 & \left . \left .  
 {\mathbf o}_j^T \overline{{\mathbf C}_{{i}}^{(2),M}}\, \overline{\mathbf p_i^M} +
 ({\mathbf t}_i^M)^T {\mathbf p}_j^M +  (\overline{\mathbf{t}_j^M})^T \overline{\mathbf p_i^M} 
\right ) \right ), 
\end{align}
where $ {\mathbf C}^{M}_{ij} : = ( {\mathbf U}_i^M)^H {\mathbf C} {\mathbf U}_j^M \in {\mathbb C}^{M \times M}$ ,  ${\mathbf C}^{(2),M}_{j} :=  {\mathbf C}^{(2)} {\mathbf U}_j^M \in {\mathbb C}^{M_d \times M}$,   ${\mathbf t}_i^M ={\mathbf t}_i  {\mathbf U}_j^M \in {\mathbb C}^M$ and further efficiencies are made by precomputing  {$ {\mathbf o}_i^T {\mathbf C}^{(1)} {\mathbf o}_j + {\rm c}_{ij} +
{\mathbf s}_i^T {\mathbf o}_j +
 {\mathbf s}_j^T {\mathbf o}_i$ and ${\mathbf o}_i^T {\mathbf C}_{{j}}^{(2),M}$.}

\begin{remark}
Once ${\mathbf K}_{ij}^M$, $ {\mathbf C}^{M}_{ij} $,  ${\mathbf t}_i^M$ and ${\mathbf o}_i^T {\mathbf C}^{(1)} {\mathbf o}_j + {\rm c}_{ij} +2 {\mathbf s}_i^T {\mathbf o}_j$ and ${\mathbf o}_i^T {\mathbf C}^{(2),M}_j$ have been precomputed, the cost of computing $(\mathcal{R}^{PODP} [\alpha B, \omega,\sigma_*,\mu_r])_{ij}$ and $(\mathcal{I}^{PODP} [\alpha B, \omega,\sigma_*,\mu_r])_{ij}$ is at most that of computing two vector-matrix-vector products involving ${\mathbf p}_i^M$ and ${\mathbf p}_j^M$ and four vector  scalar products. The small square $M \times M$ matrices and short vectors of length $M$ involved are independent of the mesh size and element order,  and so the cost of evaluation of (\ref{eqn:podprealmpt}) and (\ref{eqn:podpimagmpt}) is $O(M^2)$ and the cost of solving (\ref{eqn:ReducedA}) is at most $O(M^3)$. As $M$ is small, this is considerably cheaper than the repeated solution of (\ref{eqn:Linear}) for new parameters, {where each iteration of the iterative solve} involves a matrix vector product requiring $O(nz)$ operations {and} $nz$ is the number of non-zeros {in} ${\mathbf A}$. Furthermore, the aforementioned matrices and vectors needed for the PODP prediction can be computed once for all frequencies of interest. Further efficiencies can also be made by choosing $\tilde{\boldsymbol N}^{(n) } = {\boldsymbol N}^{(n)}$ as per Remark~\ref{rem:basisfun}. Similar efficiencies can also be gained when using PODP for other problem parameters. 
 \end{remark}

Once the aforementioned matrices and vectors have been precomputed, and the system (\ref{eqn:ReducedA}) solved for ${\mathbf p}_i^M[ { \omega}]$, the MPT coefficients obtained from (\ref{eqn:podprealmpt}) and (\ref{eqn:podpimagmpt}) are independent of the mesh and discretisation and the cost depends only on $M$ and is of order $O(M^2)$ for each $\omega$. We call this approach the Matrix Method (MM).

\subsection{Adaptive Selection of Frequency Snapshots {for the Off-line Stage}}
{
When the PODP approach is applied to the prediction of the MPT spectral signature, the off-line stage requires the selection of snapshot frequencies $\omega_n$, $n=1,\ldots,N$. Previous work~\cite{Wilson2021} has obtained accurate predictions by choosing the snapshots to be logarithmically spaced. However, in general, in order to achieve a desired accuracy in the on-line stage, the number of snapshots and their values will be problem dependent. To assist with this, an adaptive greedy algorithm along the lines of that suggested in~\cite{Hesthaven2014} and, based on our groups earlier a--posteriori error estimate, is proposed to choose new snapshots.} The estimate is restated below{:}
\begin{lemma}[Wilson, Ledger~\cite{Wilson2021}] \label{lemma:certificate}
An {a-posteriori error estimate} for the tensor coefficients computed using PODP is
\begin{subequations}
\label{eqn:certifcate}
\begin{align}
\left | (\mathcal{R}[\alpha B, \omega,\sigma_*,\mu_r ])_{ij} - (\mathcal{R}^{PODP}[\alpha B, \omega,\sigma_*,\mu_r ])_{ij}  \right |\le & (\Delta[\omega])_{ij} ,\\
\left | (\mathcal{I}[\alpha B, \omega,\sigma_*,\mu_r ])_{ij} - (\mathcal{I}^{PODP}[\alpha B, \omega,\sigma_*,\mu_r ])_{ij}  \right | \le &(\Delta[\omega])_{ij},
\end{align}
\end{subequations}
where 
\begin{align}
(\Delta[\omega])_{ij}: =  \frac{\alpha^3}{8\alpha_{LB}}
\left (  \| \hat{\boldsymbol{r}}_i [\omega] \|_{Y^{(hp)}}^2 + \| \hat{\boldsymbol{r}}_j [\omega] \|_{Y^{(hp)}}^2 + \| \hat{\boldsymbol{r}}_i [\omega] - \hat{\boldsymbol{r}}_j [\omega] \|_{Y^{(hp)}}^2 
\right ) ,  \nonumber
\end{align}
{$Y^{(hp)}: = Y^\varepsilon \cap W^{(hp)}$} and $\alpha_{LB}$ is a lower bound on a stability constant.
\end{lemma}

Note that the above error estimate applies to both $\left | (\mathcal{R}[\alpha B, \omega,\sigma_*,\mu_r ])_{ij} - (\mathcal{R}^{PODP}[\alpha B, \omega,\sigma_*,\mu_r ])_{ij}  \right |$ and
$\left | (\tilde{\mathcal{R}}[\alpha B, \omega,\sigma_*,\mu_r ])_{ij} - (\tilde{\mathcal{R}}^{PODP}[\alpha B, \omega,\sigma_*,\mu_r ])_{ij}  \right |$, since $({\mathcal N}^{0,PODP}[\alpha B,\mu_r])_{ij}= ({\mathcal N}^{0}[\alpha B,\mu_r])_{ij}$ is computed once and independently of {$\omega$}. {The  practical computation of $\| \hat{\boldsymbol{r}}_i [\omega] \|_{Y^{(hp)}} $  and $\| \hat{\boldsymbol{r}}_i [\omega] - \hat{\boldsymbol{r}}_j [\omega] \|_{Y^{(hp)}}$  is achieved using
\begin{align}
\| \hat{\boldsymbol{r}}_i [\omega] \|_{Y^{(hp)}} = &  \left ( 
\left( {\mathbf w}^{(i)}[\omega]  \right )^H {\mathbf G}^{(i,i)} {\mathbf w}^{(i)}[\omega] \right)^{1/2} ,\nonumber \\
\| \hat{\boldsymbol{r}}_i [\omega] - \hat{\boldsymbol{r}}_j [\omega]
 \|_{Y^{(hp)}} = & \left ( 
 \| \hat{\boldsymbol{r}}_i [\omega] \|_{Y^{(hp)}}^2 +  \| \hat{\boldsymbol{r}}_j [\omega] \|_{Y^{(hp)}}^2 - 2 \text{Re}\left (\left (
  {\mathbf w}^{(i)}[\omega]  \right )^H {\mathbf G}^{(i,j)} {\mathbf w}^{(j)}[\omega] 
  \right ) \right )^{1/2}, \nonumber
\end{align}
and definitions of ${\mathbf w}^{(i)}[\omega]$ and ${\mathbf G}^{(i,j)}$  that allow for a fast computation (see~\cite{Wilson2021}[pg. 1951]).
}

While the error estimate provides an estimate of the accuracy of the on--line prediction \\$\mathcal{M}^{PODP}[\alpha B, \omega,\sigma_*,\mu_r ]$ for new $\omega$ with respect to the full order model, it can also be used to adaptively choose the snapshots for the off-line stage, which is the approach followed here. The goal being to choose snapshots in the off-line stage so that the predictions in the on-line stage result in a small  $(\Delta[\omega])_{ij}$ for all $\omega$ of interest leading to the adaptive procedure presented in Algorithm~\ref{alg:adapt}.

\begin{algorithm}[h]
\begin{algorithmic}[1]
\State Prescribe $TOL_\Sigma$, $TOL_\Delta$, {$N_*\ge 1$, $k=1$, $MAX_{k}>1$, $N_{output}\gg N>1$}.
\State {\bf Initial Off-line Stage}
\State Choose initial log-spaced frequency snapshots $\omega_n$, $n=1, \ldots,N$. 
\State Compute ${\mathbf q}_i[\omega_n]$, $i=1,2,3$, $n=1,\ldots,N$ by solving (\ref{eqn:Linear}).
\State Form $\mathbf{D}_i \vcentcolon=\big[\mathbf{q}_i[{\omega}_1],\mathbf{q}_i[{\omega}_2],...,\mathbf{q}_i[  {\omega}_{N}]\big]$ for $i=1,2,3$.
\State Based on $TOL_\Sigma$,  obtain  the TSVDs $ \mathbf{D}_i^M = \mathbf{U}_i^M\mathbf{\Sigma}_i^M(\mathbf{V}_i^M)^{{H}}$ for $i=1,2,3$.
\State {\bf Adaptive Off-line Computation}
\While{$\Lambda / |B_\alpha|  > TOL_\Delta$ and $k  < MAX_{k}$}
\State $k \to k+ 1$ 
{\For{$n=1,\ldots,N-1$ snapshot intervals}
\State Obtain certificates $(\Delta[\omega])_{ij}$, $i,j=1,2,3$ for $\omega_{n } \le \omega < \omega_{n+1}$.
\State Set $\displaystyle \Lambda_{n}[\omega]:= \max_{\substack{\omega_{n} \le \omega < \omega_{n+1} \\ i,j=1,2,3}} (\Delta[\omega])_{ij}$ as the maximum error per interval.
\EndFor
\State Set $\displaystyle \Lambda := \max_{n=1,\ldots, N-1} \Lambda_{n}[\omega]$ as the global maximum.
\State Determine extra snapshots $\omega_{n_*}$, $n_*=1,\ldots, N_*$  as those $\omega_{n_*}$ corresponding to the $N_*$ largest $\Lambda_{n}[\omega]$
\State Compute ${\mathbf q}_i[\omega_{n_*}]$, $i=1,2,3$, $n_*=1,\ldots,N_*$ by solving (\ref{eqn:Linear}).}
\State $N \to N+N_*$.
\State Update $\mathbf{D}_i \vcentcolon=\big[\mathbf{q}_i[{\omega}_1],\mathbf{q}_i[ {\omega}_2],...,\mathbf{q}_i[ {\omega}_{N}]\big]$ for $i=1,2,3$.
\State Based on $TOL_\Sigma$,  update the TSVDs $ \mathbf{D}_i^M = \mathbf{U}_i^M\mathbf{\Sigma}_i^M(\mathbf{V}_i^M)^{{H}}$ for $i=1,2,3$.
\EndWhile
\State {\bf On-line Stage}
{\State{Choose output frequencies $\omega_n^o$, $n=1, \ldots, {N_{output}}$}
\For{$n=1,\ldots,N_{output}$}
\State Solve $\mathbf{A}_i^M[{\omega_n^o}]\mathbf{p}_i^M[  {\omega}_n^o]=\mathbf{r}_i^M(\boldsymbol{\theta}_i^{(0,hp)})[\omega_n^o]$ for $\mathbf{p}_i^M[  {\omega_n^o}]$, $i=1,2,3$.
\State Compute  $(\mathcal{R}^{PODP} [\alpha B, \omega_n^o,\sigma_*,\mu_r])_{ij}$ and $(\mathcal{I}^{PODP} [\alpha B, \omega_n^o,\sigma_*,\mu_r])_{ij} $, $i,j=1,2,3$, using (\ref{eqn:podprealmpt}) and (\ref{eqn:podpimagmpt}).
\EndFor}
\end{algorithmic}
\caption{{Adaptive algorithm for PODP prediction of the MPT spectral signature}} \label{alg:adapt}
\end{algorithm}

Here, $N_*$ controls how many additional snapshots are generated in each adaption.  In this work, a fixed value of $N_* = 2$ is used, with at most one additional snapshot being introduced per interval. A common, and equally valid implementation, is to set a threshold parameter $0 < \vartheta \le 1$ and choose new snapshot parameters as those $\omega_{n_*}$ such that $ (\Delta[\omega_{n_*}])_{ij}\ge \vartheta \Lambda$~\cite{Hesthaven2014}. This has the benefit that different numbers of snapshots can be included at different stages of the adaption, however, for the present application, this was found to lead to a more variable computation time and memory requirements.

\section{Computational Resources and Software}\label{sect:compres}
The computational resources used to perform the simulations reported in this work are the following workstation:

{\bf Workstation}: Intel Xeon W-2265 CPU (12 cores 24 threads)  with a clock speed of 3.50 GHz and 128GB of DDR4 RAM with a speed of 3200 MT/s.

{All timings were} performed using wall clock times with the \texttt{Memory-Profiler} package (version 0.61.0) and \texttt{Python} version 3.10. FEM simulations were performed using {the} {\texttt{NGSolve} library~\cite{NGSolve}, based on the basis functions proposed in~\cite{zaglmayrphd}, and the associated \texttt{Netgen} mesh generator~\cite{netgendet} using versions 6.2204 and 6.2203, respectively,} and using \texttt{NumPy} (version 1.23.3). These libraries are called from  our group's open source \texttt{MPT-Calculator}\footnote{\texttt{MPT-Calculator} is publicly available at \url{https://github.com/MPT-Calculator/MPT-Calculator/}} (April 2024 release) software.

In the off-line stage, two different forms of parallelism are {applied:}  The assembly of the matrices and the underlying iterative solution of (\ref{eqn:Linear}), which requires repeated matrix-vector products in the conjugate gradient solver, {are} accelerated by using shared memory parallelism across multiple computational threads {without requiring additional memory}. These operations are trivially parallelisable within the \texttt{NGSolve} {library.} Provided sufficient memory resources are available, the computation of the full order solutions {${\mathbf q}_i[\omega_n]$}  for different $\omega_n$ is further accelerated by using multi-processing with different {frequencies} being considered simultaneously, which leads to higher memory demands. {For the parallel computation of ${\mathbf q}_i[\omega_n]$ across $P$ computational processes (cores), stored memory is duplicated. Thus, memory usage for obtaining these solutions for multiple snapshot parameters  will be at least $P$ times greater.}

In the on-line stage, the computation of the solution ${\mathbf p}_i^M$ to (\ref{eqn:ReducedA}) and the computation of\\ $(\mathcal{R}^{PODP} [\alpha B, \omega,\sigma_*,\mu_r])_{ij}$ and $(\mathcal{I}^{PODP} [\alpha B, \omega,\sigma_*,\mu_r])_{ij}$ using (\ref{eqn:podprealmpt})  and (\ref{eqn:podpimagmpt}), respectively, which has already been reduced to small matrix vector products,  is further accelerated where appropriate.  We give examples of the actual timings in the following sections.

\section{{A Recipe for Choosing the} Number and Thicknesses of Prismatic Boundary Layers}\label{sect:boundarylayers}

The depth at which the amplitude of the {Ohmic currents in the conductor} decays to $1/e$ of its surface value is known as the skin-depth and, for a homogeneous isotropic conductor, is commonly approximated by~\cite{balanis2012}
\begin{equation}\label{eqn:delta}
\delta[\omega, \sigma_*, \mu_r] \approx \sqrt{\frac{2}{\omega \sigma_* \mu_0\mu_r}},
\end{equation}
which, for high $\sigma_*$, $\omega$ and $\mu_r$, can become very small compared to the object dimensions. This measure applies to ${\boldsymbol J}_i^e = \im \omega \sigma_*{\boldsymbol \theta}_i^{(1)}$, since (\ref{eqn:theta1transtrun}) is an eddy current type transmission problem and ensuring that the skin is accurately computed is also important for the computation of the MPT coefficients.

Combining prismatic boundary layer elements with unstructured tetrahedral meshes with $p$--refinement of the elements has previously been shown to {result in faster rates of convergence than} $hp$-refinement of purely tetrahedral meshes for {the MPT coefficients~\cite{elgy2023}. They have also been used for modelling  thin material coatings in plated objects.}
 Similar {increased rates of convergence} have also been reported in other applications, such as the Maxwell eigenvalue problem~\cite{zaglmayrphd} and singularly perturbed elliptic boundary value problems~\cite{schawbsuri}. Here, the prismatic layers allow $h$--refinement to be achieved  in a direction normal to the surface of the conductor, while leaving the tangential element spacing unchanged, which is ideal for addressing the high field gradients in the normal direction, but without resulting in a large increase in the number of degrees of freedom. On the other hand, using $p$--refinement alone on traditional unstructured meshes is sub--optimal, {converging} only at an algebraic rate due to the small skin depths, while {$h$--refinement of the unstructured tetrahedral mesh leads to an excessive $N_d$}. {As pointed out in Section~\ref{sect:romoffline},  the same FEM discretisation is used for all $\omega_n$  off--line snapshots and, in order  to ensure that this is be accurate for the complete signature,
a maximum target $\omega$ of interest  is fixed and, for a given $\mu_r$,  $\tau:=\sqrt{2/ (\mu_r\nu)}= \delta/\alpha$ is the smallest non-dimensional skin-depth 
that is to be resolved.}

While offering considerable benefits, the inclusion of prismatic layers must {never}theless also be weighed up against the increase in computational resources (including both run time and memory usage). Our goal in this section is to determine a simple recipe for choosing the number and thicknesses of prismatic boundary layers that can be used to achieve a high level of accuracy at a reasonable computational cost. {Three schemes for defining the structure of $L$ prismatic layers are considered:}

\begin{figure} 
\centering
$\begin{array}[t]{c c c}
\includegraphics[width=0.3\textwidth]{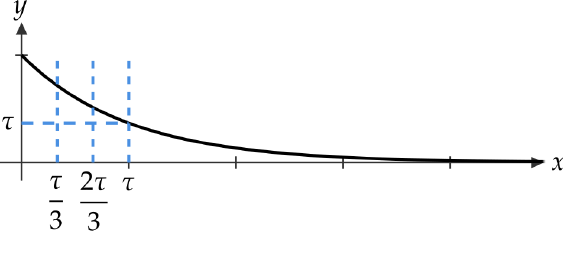} & \includegraphics[width=0.3\textwidth]{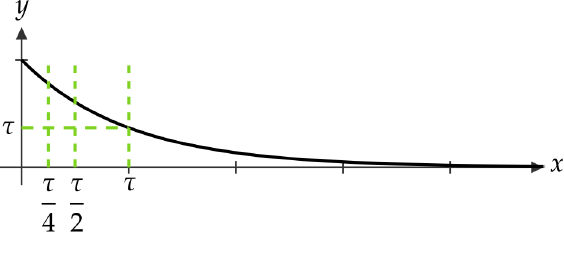} & \includegraphics[width=0.3\textwidth]{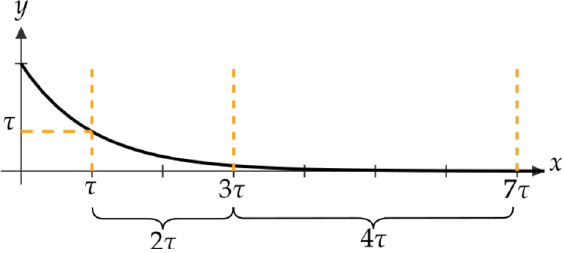} \\
(a) & (b) & (c)
\end{array}$
\caption{Proposals for boundary layer thickness for the simple function $y=e^{-x}$  showing the thickness{es} for $L=3$ layers of elements in terms
 of the non-dimensional skin depth $\tau$ being for  $(a)$ ``uniform'' distribution, $(b)$ ``geometric decreasing'' distribution, and $(c)$ ``geometric increasing'' distribution. }
\label{fig:boundary_layer_discretisation}
\end{figure}

\begin{enumerate}
\item ``Uniform'' refinement, where the total thickness of the layers is equal to $\tau$ and each layer of prismatic elements has thickness $t_\ell=\tau/L, $ $\ell=1,\ldots,L$ with $t_1$ being the closest to the conductor's surface and the layers numbered consecutively towards the inside of the conductor.
\item ``Geometric decreasing'' refinement, where the total thickness is limited to $\tau$ and the thickness of each layer is defined by the geometric series $t_{\ell+1} = 2t_\ell$ with $\sum_{l=1}^L t_\ell = {t_1 (2^L-1)} =  \tau$. 
\item ``Geometric increasing'' refinement, where the thickness of the layers are defined as $t_{\ell+1} = 2t_\ell$ with $t_1 = \tau$. We also ensure that $L$ is chosen such that no difficulties arise in the case of thin conducting objects.
\end{enumerate}
An illustration of each refinement strategy is shown in Figure~\ref{fig:boundary_layer_discretisation}, showing the thicknesses for $L=3$ layers of elements in terms of the non-dimensional skin depth $\tau$ for the three strategies.
In each of these schemes, the total number of prismatic elements remains constant and, in the case of $L=1$, the uniform, geometric decreasing and geometric increasing strategies all lead to identical discretisations. This means that for a given mesh, order of elements{,} and number of layers {$N_d$ and $M_d$ remain} the same for all three strategies.

A sphere  of radius $\alpha=1\times 10^{-3}$ m is considered, and to illustrate different skin depths, the cases corresponding to $\sigma_* = 1\times 10^6$ S/m, $\mu_r=1,16,64$ and
a fixed target frequency $\omega = 1\times 10^8$ rad/s are compared for the three strategies. Solutions are obtained using a a computational domain $\Omega$ consisting of a unit radius sphere $B$ placed centrally  in a $[-1000, 1000]^3$ units box and a mesh of 21\,151 unstructured tetrahedra. Here, the curved surface of the sphere is approximated by polynomials of order $g=5$ independently of the order of the elements, $p$. 
The mesh is augmented by the addition of prismatic layers placed just inside $\Gamma$ according to the three strategies above.

{Figure~\ref{fig:boundary_layer_tests} shows the relative error between the approximated MPT coefficents and the exact solution for the sphere~\cite{Wait1951}, $E = \lVert \mathcal{M}^{hp} - \mathcal{M}\rVert_F / \lVert\mathcal{M}\rVert_F$, where $\lVert \cdot \rVert_F$ denotes the Frobenius norm, under $p$--refinement for each of the different approaches for $L=1,2,3,4$ layers and $\mu_r= 1, 16, 64$, in turn.
Here, and unless otherwise stated, the regularisation was set as $\varepsilon = 1\times 10^{-10}$ and the iterative solver relative tolerance as $TOL = 1\times 10^{-8}$.
For each value of $\mu_r$, the convergence of $E$ with respect to $N_d$ and $E$ with respect to computational time (showing overall time using $P=2$ cores on the workstation described in Section~\ref{sect:compres} and the IM method described in Section~\ref{sect:mpt}) are shown.
In each case, the points on each curve correspond to increasing element order.

\begin{figure}
\centering
$\begin{array}{c c}
\includegraphics[width=0.45\textwidth]{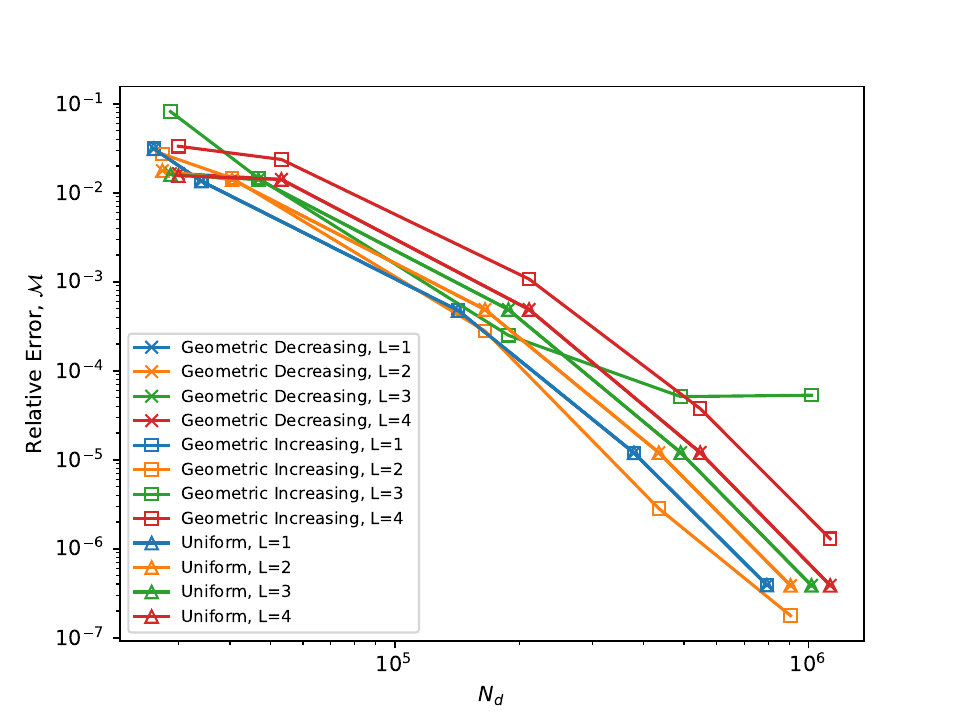} &
\includegraphics[width=0.45\textwidth]{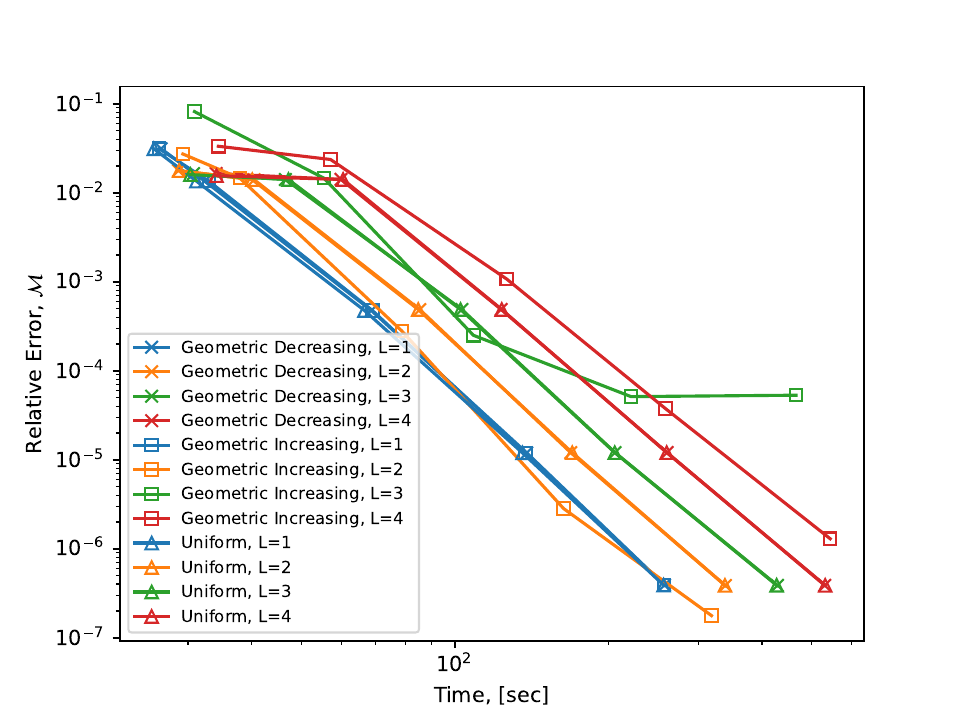} \\
(a) & (b)\\
\includegraphics[width=0.45\textwidth]{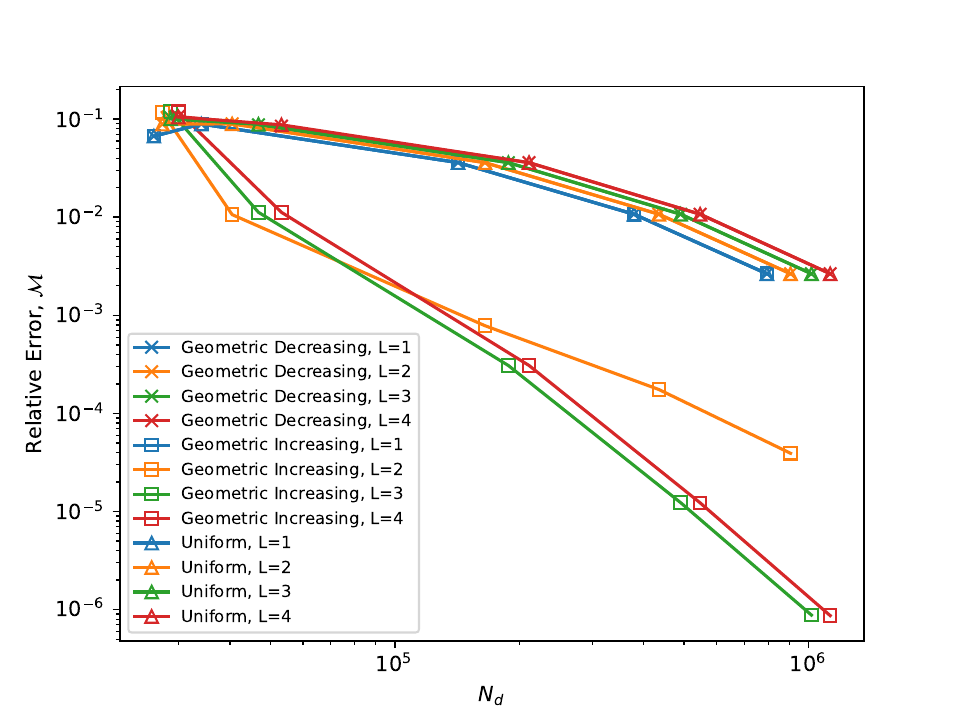} &
\includegraphics[width=0.45\textwidth]{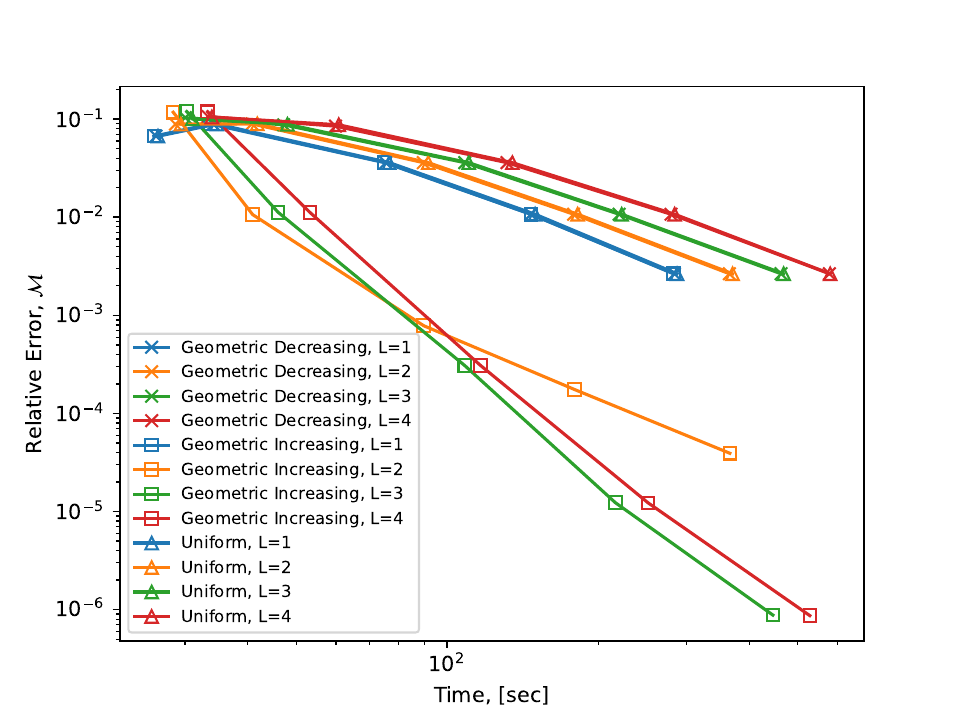} \\
(c) & (d)\\
\includegraphics[width=0.45\textwidth]{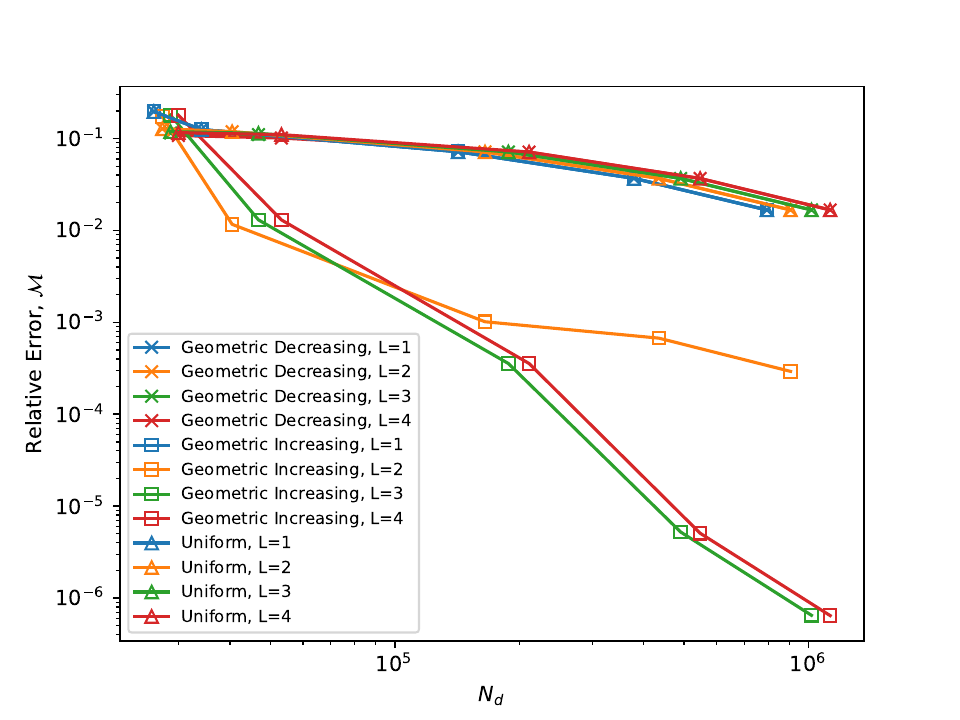} &
\includegraphics[width=0.45\textwidth]{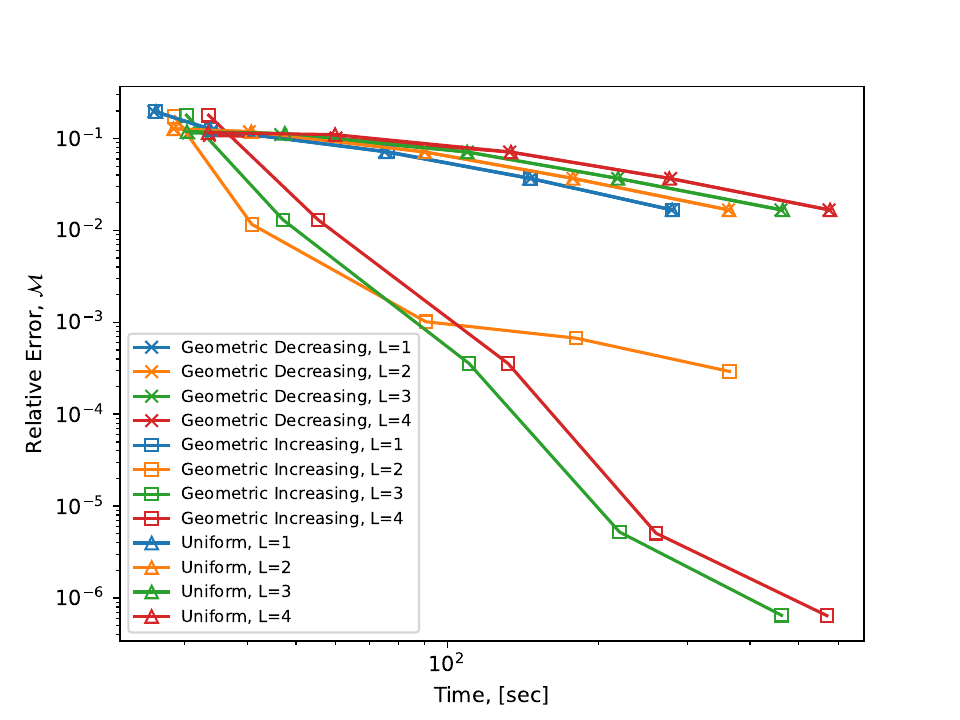} \\
(e) & (f) 
\end{array}$
\caption{Magnetic conducting sphere:  Showing the effect of $p$--refinement on  $E = \lVert \mathcal{M}^{hp} - \mathcal{M}\rVert_F / \lVert\mathcal{M}\rVert_F$ for the different prismatic layer strategies with respect to number of degrees of freedom (left column) and computational time (right column) for  $\mu_r =1 \, (a, b), \,  16 \, (c,d),$ and $64\, (e,f)$.}
\label{fig:boundary_layer_tests}
\end{figure} 

In the case of $\mu_r=1$, and the aforementioned parameters, the resulting skin depth can be well resolved for all three schemes and $p$-refinement. Furthermore, with the exception of $L=3$ and the geometric increasing strategy, all schemes lead to similar convergence curves in terms of both $N_d$ and time. A small benefit is observed for $L=2$ using the geometric increasing scheme over the other schemes. Indeed, for this case, using the initial tetrahedral mesh alone is already able to achieve exponential convergence when $p$--refinement is applied.

For $\mu_r>1$ the uniform and geometric decreasing strategies are seen to produce similar results for all values of $L$, while there is  a considerable benefit in convergence rate, and by extension accuracy, by using the geometric increasing scheme with $L\ge 2$ both with respect to $N_d$ and computational time. As $\mu_r$ increases, further benefits in accuracy with respect to {both} $N_d$ and computational time are offered by using $L\ge 3$ and, while not shown, exponential convergence with respect to $N_d^{1/3}$ is obtained if $\log E$ is plotted against  $N_d^{1/3}$ for sufficiently large $L$. 
However, while combining $p$--refinement with $L\ge 3$ can achieve $E < 10^{-6}$,  a relative error of $E=1\times 10^{-3}$ is sufficient for classification purposes given the precision to which MPT coefficients can currently be measured, and the ability to which materials and geometry of an object to be characterised are known. This level of accuracy can be achieved using $L=2$, the geometric increasing scheme and $p$--refinement and, therefore, in the later practical computations,  this is what {is employed}. These findings are also consistent with the theory of~\cite{schawbsuri}, which would suggest a first layer of thickness $O((p+1)\tau)$ if their findings for their one-dimensional problem are extrapolated to our three-dimensional problem.

\begin{remark}
For highly magnetic conducting spheres with the aforementioned properties and $\mu_r=100,200,400, 800$, and applying geometric increasing refinement with $L=2$ layers leads to a relative error of $E=1\times 10^{-3}$ using order $p=3$ elements. The strategy has also been applied to problems with even smaller skin depths, such using the same materials with an object of size $\alpha=0.01$ m, and similar convergence behaviour is obtained.
\end{remark}

\begin{remark}
We briefly recall that the MPT characterisation of an object under size scaling $B_\alpha \to sB_\alpha$~\cite{Wilson2021}[Lemma 3] can be done exactly and at no additional computational cost by
\begin{equation*}\label{eqn:scaling}
\left(\mathcal{M}\left[s\alpha B, \omega, \sigma_*, \mu_r\right]\right)_{ij} = s^3 \left(\mathcal{M}\left[\alpha B, s^2 \omega, \sigma_*, \mu_r\right]\right)_{ij}, 
\end{equation*}
with $\left(\mathcal{M}\left[\alpha B, s^2 \omega, \sigma_*, \mu_r\right]\right)_{ij}$ denoting the original tensor coefficients at angular frequency $s^2\omega$ provided that the same discretisation is used in both cases. 
For consideration with the prismatic boundary elements; since the discretisation is applied to $B$, then we should choose the size of the prismatic elements according to $\tau[s^2\omega, \alpha/s]=\tau[\omega, \alpha]$ if $B_\alpha \to sB_\alpha$ and $\omega$ is the target frequency for the $B_\alpha$ case.
\end{remark}

\section{Numerical Examples} \label{sect:numexp}
In this section, a range of numerical examples to illustrate the improvements in accuracy and speedup for  calculation of the {MPT spectral signature} using PODP, the use of adaption to choose new frequency snapshots, and the geometric increasing recipe proposed for the construction of prismatic layers in Section~\ref{sect:boundarylayers} {are considered.}

\subsection{Conducting Sphere}\label{sect:sphere}
The conducing sphere described in Section~\ref{sect:boundarylayers} for the particular case where $\mu_r=32$ {is considered, and the computational domain $\Omega$ is again discretised by 21\,151 unstructured tetrahedra  and $L=2$ layers of prismatic elements following the geometric increasing strategy, resulting in  1\,275 prisms.  A total of $N=13$ full order solutions at logarithmically spaced snapshot (SS) frequencies $1\times 10^1\le \omega \le1\times 10^8$ rad/s are computed using order $ p=3$ elements and, using $TOL_\Sigma = 1\times 10^{-6}$, this leads to  $M=9$ TSVD modes}. 
The MPT spectral signature obtained using IM, FMM, and MM are compared in Figure~\ref{fig:magnetic_sphere_tensors} and excellent agreement is obtained for all frequencies, with lines indistinguishable on this scale. Note that due to object symmetries,  $({\mathcal M})_{11}=({\mathcal M})_{22}= ({\mathcal M})_{33}$ and $({\mathcal M})_{12}= ({\mathcal M})_{13}= ({\mathcal M})_{23}=0$, and, thus, only one on--diagonal and one off--diagonal coefficient are shown.
Similar agreement can be found for spheres using other values of $\mu_r$.
Timings performed using the workstation described in Section~\ref{sect:compres} for the IM, FMM, and MM methods indicate that with $P=2$ cores, and the use of multi-threading, the computation time is reduced from 13\,821 seconds to 237 seconds to 34 seconds for the IM, FMM, and MM approaches respectively. This gives an overall time of just 539 seconds for the entire process, including obtaining the full order solution snapshots, when using the MM approach. The breakdown in computational time is illustrated in Figure~\ref{fig:magnetic_sphere_times}, where memory usage is plotted against computation time for the cases of ($a$) IM, ($b$) FMM, and ($c$) MM methods. The computation time is further broken down into common sections for easier comparison. While the off-line stage of the POD has the same computational cost for all three methods, there are significant differences in the on--line stage of the POD. For the IM method, which is the naive implementation, the memory usage is low, however a disproportionate amount of the overall computational time is spent computing the tensor coefficients. In the FMM and MM cases, the efficient post-processing significantly decreases the overall computation time. In the FMM case, large sparse matrices, including interior degrees of freedom, are used throughout, leading to a larger memory footprint. The MM method is the fastest and uses less memory than the FMM approach since the matrices $\mathbf{K}_{ij}^M$, $\mathbf{C}_{ij}^M$, and $\mathbf{C}_{i}^{(2), M}$ are precomputed and reduced in size. Further memory savings are made by computing the real and imaginary parts of the tensor coefficients separately to avoid having $\mathbf{K}$ and $\mathbf{C}$ stored in memory at the same time. For the naive IM approach, where large matrices are not required, this is not a concern.

\begin{figure}
\centering
$\begin{array}{c c}
\includegraphics[width=0.45\textwidth]{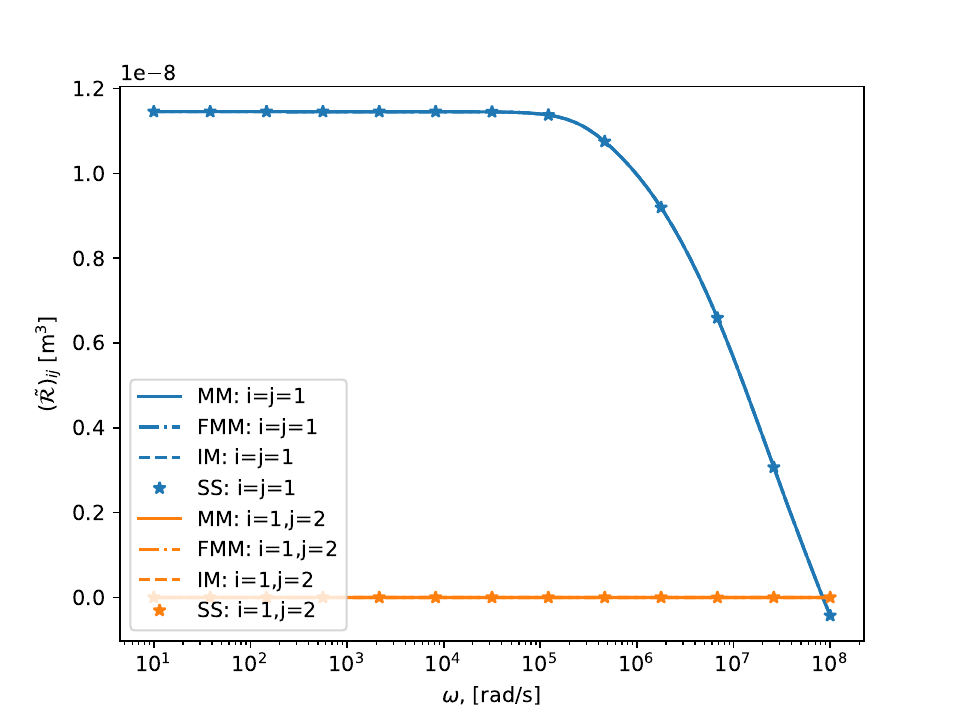} & \includegraphics[width=0.45\textwidth]{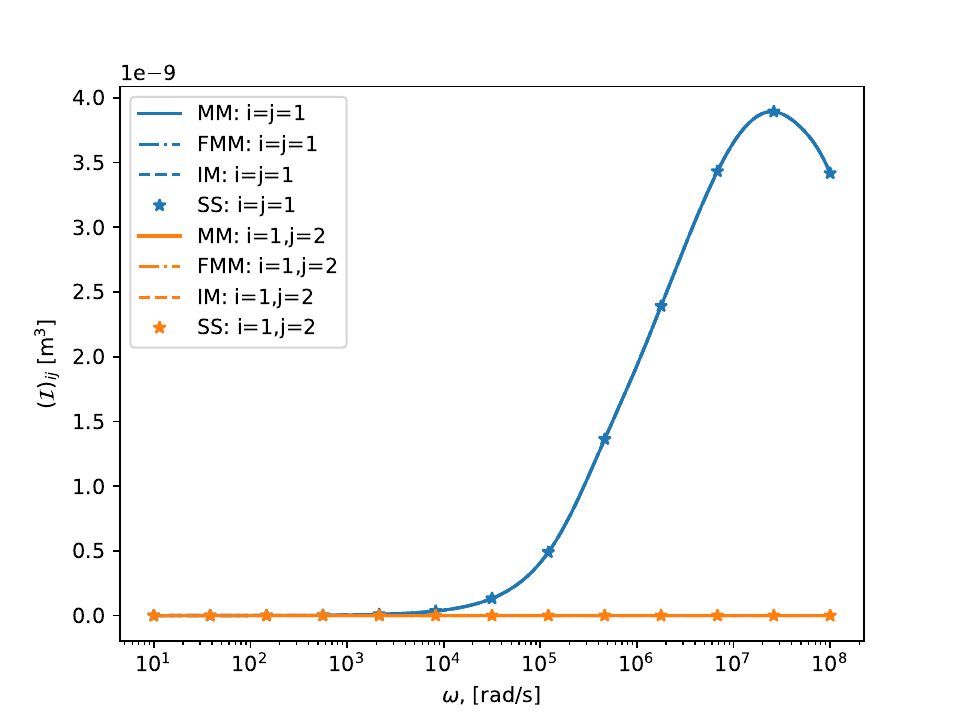} \\
(a) & (b) 
\end{array}$
\caption{Magnetic conducting sphere: 
Showing a comparison between the original IM and the new faster FMM and MM approaches for the calculation of the  MPT spectral signature  using PODP $(a)$ $(\tilde{\mathcal R})_{ij}$ and $(b)$ $({\mathcal I})_{ij}$.}

\label{fig:magnetic_sphere_tensors}
\end{figure}

\begin{figure}[h!]
\centering
$\begin{array}{c}
\includegraphics[width=0.45\textwidth]{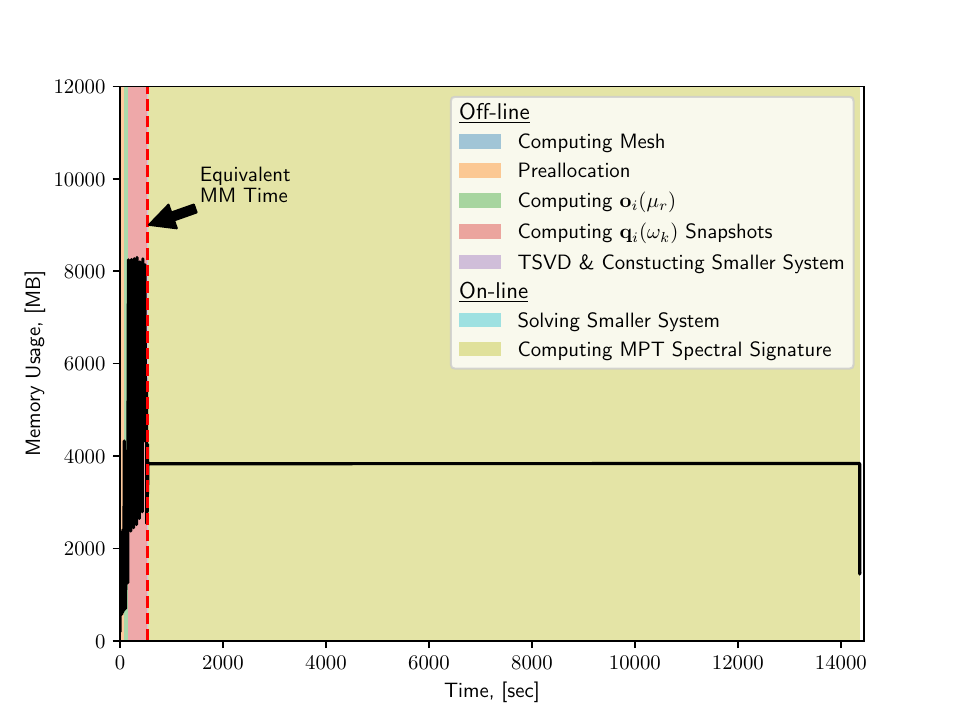} \\
(a)
\end{array}$
$\begin{array}{c c}
\includegraphics[width=0.45\textwidth]{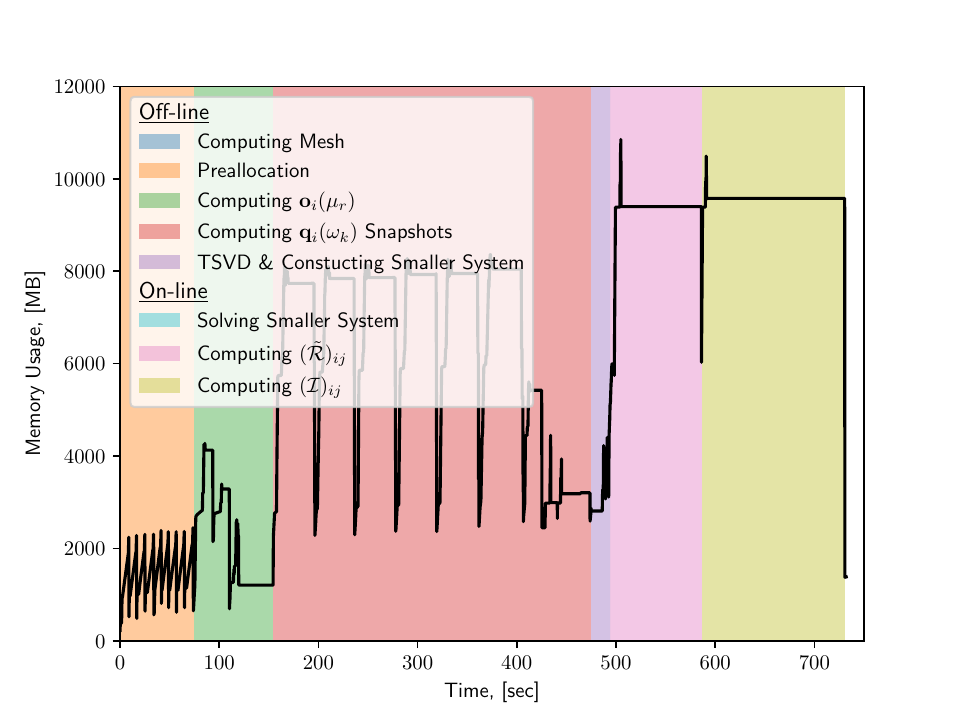} & 
\includegraphics[width=0.45\textwidth]{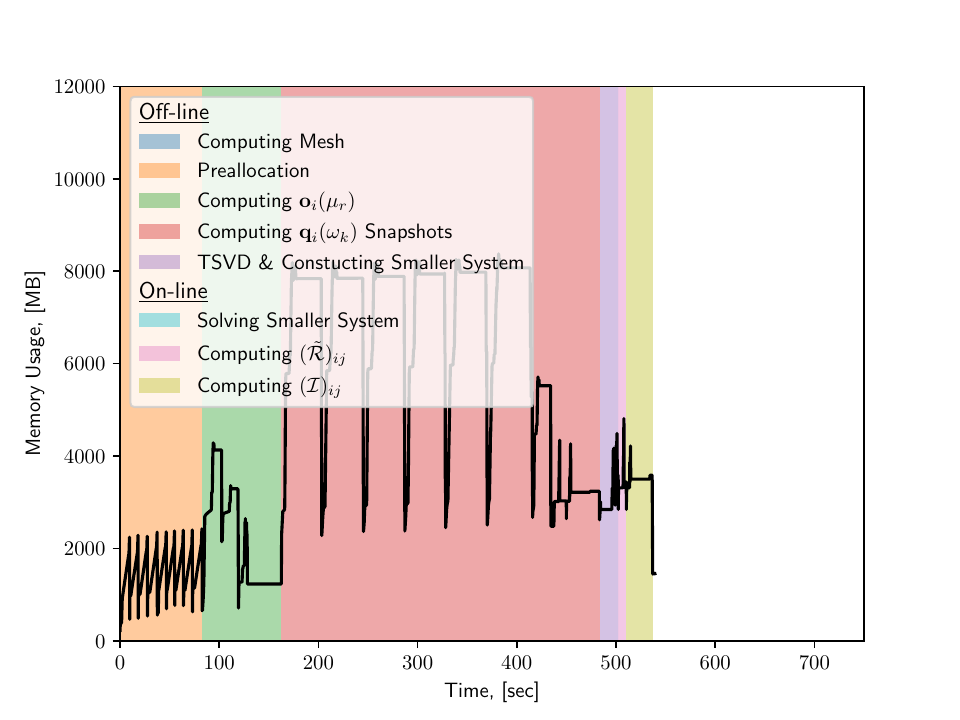} \\
(b) & (c)
\end{array}$

\caption{Magnetic conducting sphere:
A comparison between time and memory usage for ($a$) IM, ($b$) FMM, and ($c$) MM methods showing substantial speed up for the on--line stage of POD. Timings are further broken down into common parts of the problem to correlate memory usage with specific tasks.}
\label{fig:magnetic_sphere_times}
\end{figure}


The adaptive procedure outlined in Algorithm \ref{alg:adapt} is demonstrated for the same discretisation in Figure~\ref{fig:Adaptive_Sphere}, which shows the MPT  spectral signature for $(\tilde{\mathcal R}[\alpha B, \omega,\sigma_*,\mu_r ])_{ij}$, that we subsequently refer to as $(\tilde{\mathcal R})_{ij}$, including the a-posteriori error certificates $(\tilde{\mathcal R}\pm \Delta)_{ij}$ obtained at different iterations where $(\Delta)_{ij}$
reduces as new snapshots are adaptively chosen. In this example $TOL_\Delta = 1\times 10^{-3}$ {and results for the iterations $k=1,2,3,4$ are shown, resulting in $N=13,15,17,19$ snapshots, respectively}. 
Importantly, while the effectivity indices of the error certificates are large, {the additional cost  required to compute them during the on--line stage is small} and, as the figures show, {the algorithm} provides an effective way to choose new snapshots to reduce $(\Delta)_{ij}$.
Only the  behaviour for $(\tilde{\mathcal R})_{ij}$  {is} shown here, since the error certificates are the same for both $(\tilde{\mathcal R})_{ij}$ and $({\mathcal I})_{ij}$.  
Similarly to Figure~\ref{fig:magnetic_sphere_tensors}, only one on--diagonal and one off--diagonal coefficient are shown since the other tensor coefficients and error certificates are indistinguishable on this scale.
{Figure~\ref{fig:Adaptive_Sphere_Convergence} shows the maximum error, $\Lambda$, against $N$ to }illustrate the performance of the adaptive POD, compared to non-adaptive logarithmically spaced {snapshots}. The figure shows the improvement associated with the adaptive POD compared to the non-adaptive scheme. Nevertheless, using {logarithmically spaced} snapshots with $N=13$ still provides  a very good starting point for an initial choice of frequencies from which adaption can then be performed.

\begin{figure}[!h]
\centering
$\begin{array}{cc}
\includegraphics[width=0.45\textwidth]{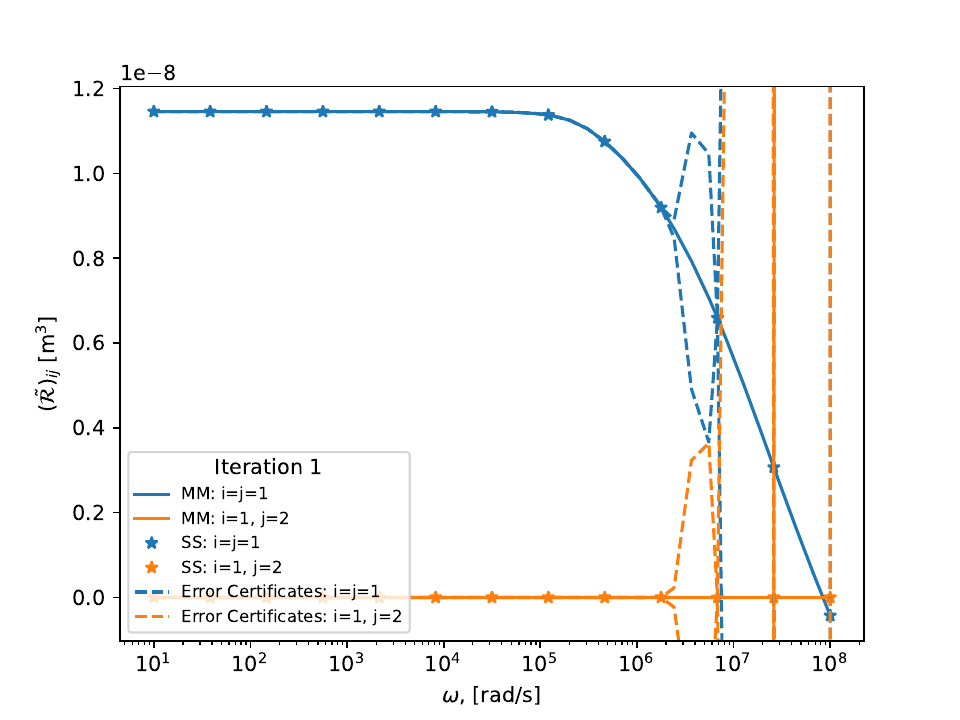} & \includegraphics[width=0.45\textwidth]{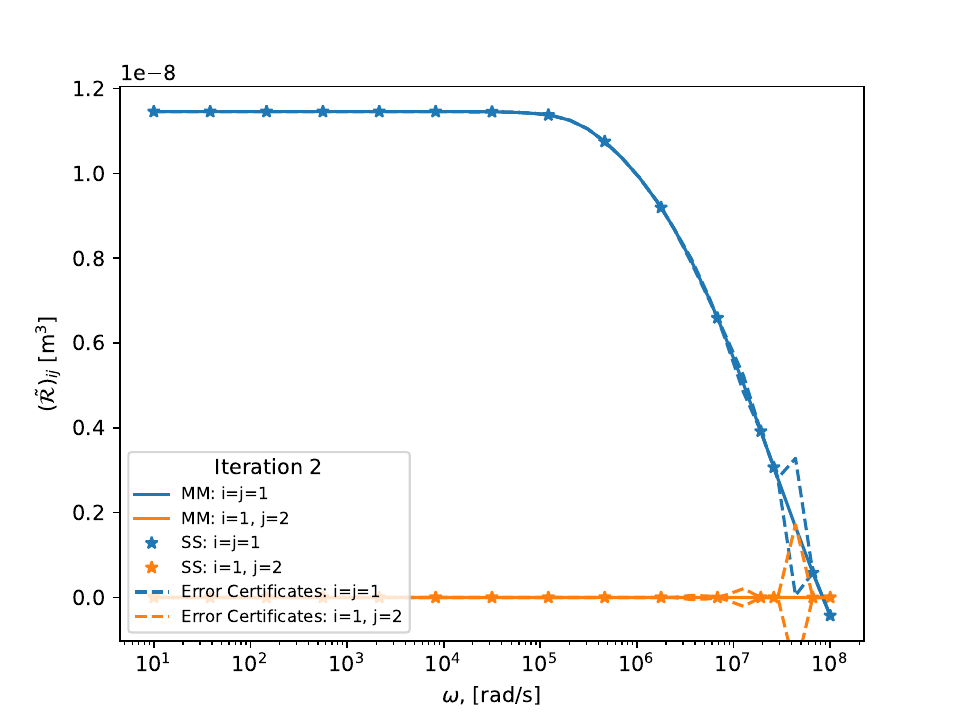} \\
(a) & (b) \\
\includegraphics[width=0.45\textwidth]{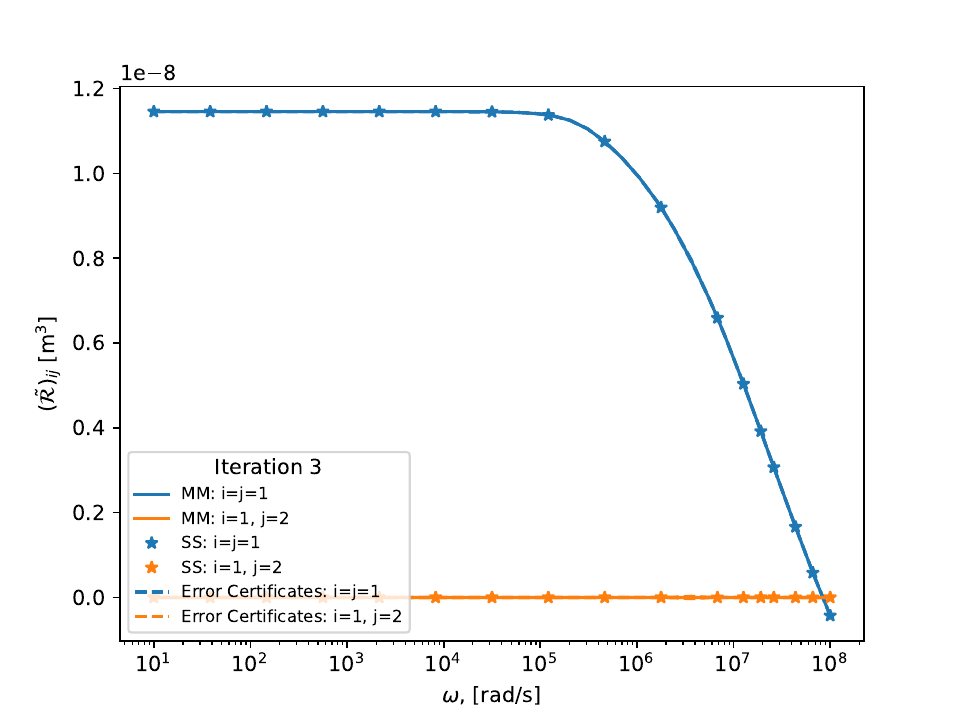} & \includegraphics[width=0.45\textwidth]{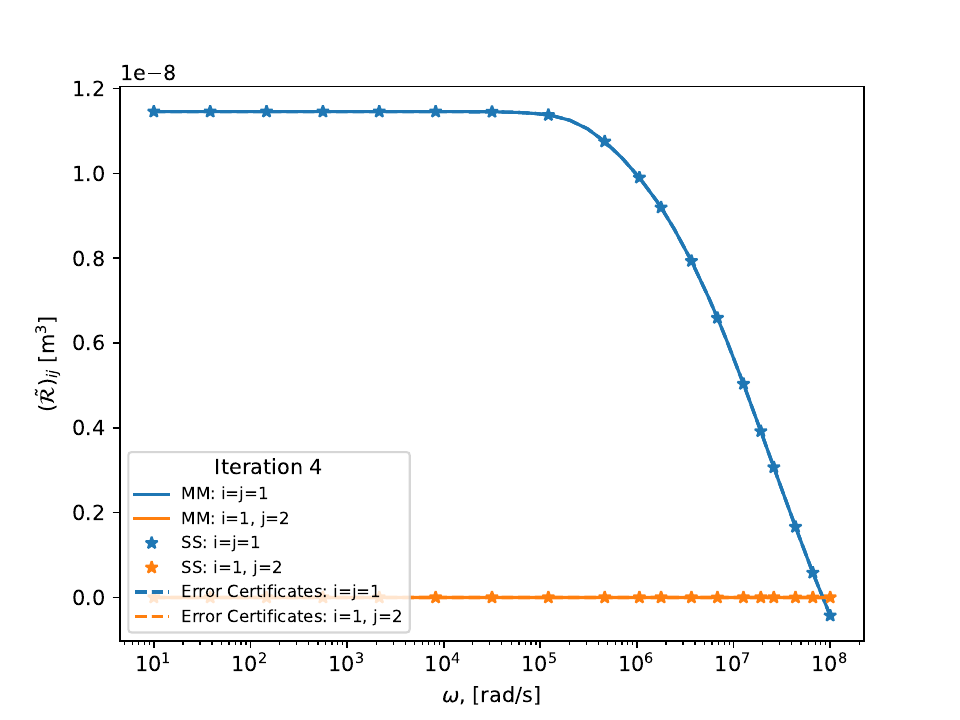} \\
(c) & (d)
\end{array}$
\caption{Magnetic conducting sphere:
{Showing  $(\tilde{\mathcal R})_{ij}$ and $(\tilde{\mathcal R}\pm \Delta)_{ij}$ obtained by applying adaptive PODP in Algorithm~\ref{alg:adapt} at iterations $(a)$ $k=1$, $(b)$ $k=2$, $(c)$  $k=3$, and $(d)$  $k=4$. }}

\label{fig:Adaptive_Sphere}
\end{figure}

\begin{figure}[!h]
\centering
\includegraphics[width=0.45\textwidth]{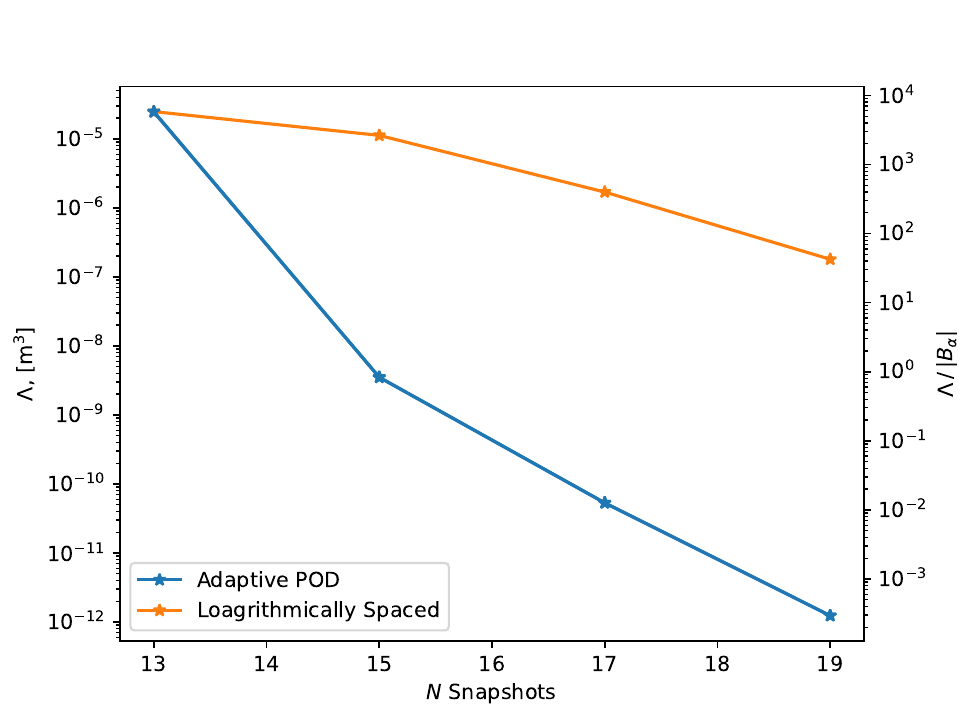}
\caption{Magnetic conducting sphere: { Showing $\Lambda$ and $\Lambda/|B_\alpha|$ obtained in the adaptive PODP in Algorithm~\ref{alg:adapt}
against $N$ compared with logarithmically spaced snapshot frequencies.}}

\label{fig:Adaptive_Sphere_Convergence}
\end{figure}

\subsection{Conducting and Magnetic Disks} \label{sect:magnetic_disk}
{This section considers the MPT characterisation of thin conducting and magnetic disks with their circular face in the $x_1-x_3$ plane. A disk with physical dimensions radius $r=1 \times 10^{-2}$ m and thickness $h=1 \times10^{-3}$ m and  with  $\mu_r = 32$ and $\sigma_*=1\times 10^6$ S/m is  initially considered}. The computational domain $\Omega$ consists of a dimensionless disk $B$ of radius $r/\alpha$, and thickness $h/\alpha$ with $\alpha =1\times 10^{-3}$ m placed centrally in the box $[-1000,1000]^3$ units.  The  geometric increasing methodology from Section~\ref{sect:boundarylayers} is applied to construct prismatic boundary layers for different values of $\mu_r$ at a target frequency of $\omega = 1 \times 10^8$ rad/s. 
This process results in a discretisation of 24\,748 tetrahedra and 2\,995 prisms with $p=3$ giving converged solutions at the $N=13$ snapshot  frequencies in the range $1\times 10^1\le \omega \le1\times 10^8$ rad/s, which reduces to $M=11$ TSVD modes using  $TOL_\Sigma=1 \times 10^{-6}$.
Due to  mirror symmetries of the disk, and its rotational symmetry about ${\boldsymbol e}_2$, the non-zero independent tensor coefficients associated with the object reduce to  $({\mathcal M})_{11} = ({\mathcal M})_{33}$ and $({\mathcal M})_{22}$. For this reason,  we show only $(\tilde{\mathcal R})_{11} = (\tilde{\mathcal{R}})_{33}$, $(\tilde{\mathcal R})_{22}$, $(\tilde{\mathcal R})_{12} = (\tilde{\mathcal R})_{23} = (\tilde{\mathcal R})_{13} = 0$ and $({\mathcal I})_{11} = (\mathcal{I})_{33}$, $({\mathcal I})_{22}$, $({\mathcal I})_{12} = ({\mathcal I})_{23} = ({\mathcal I})_{13} = 0$ in Figure~\ref{fig:thin_disk_magnetic_tensors}. The figures show excellent agreement between the IM, FMM, and MM methodologies, {with curves indistinguishable on this scale}.
\begin{figure}[!h]
\centering
$\begin{array}{cc}
\includegraphics[width=0.45\textwidth]{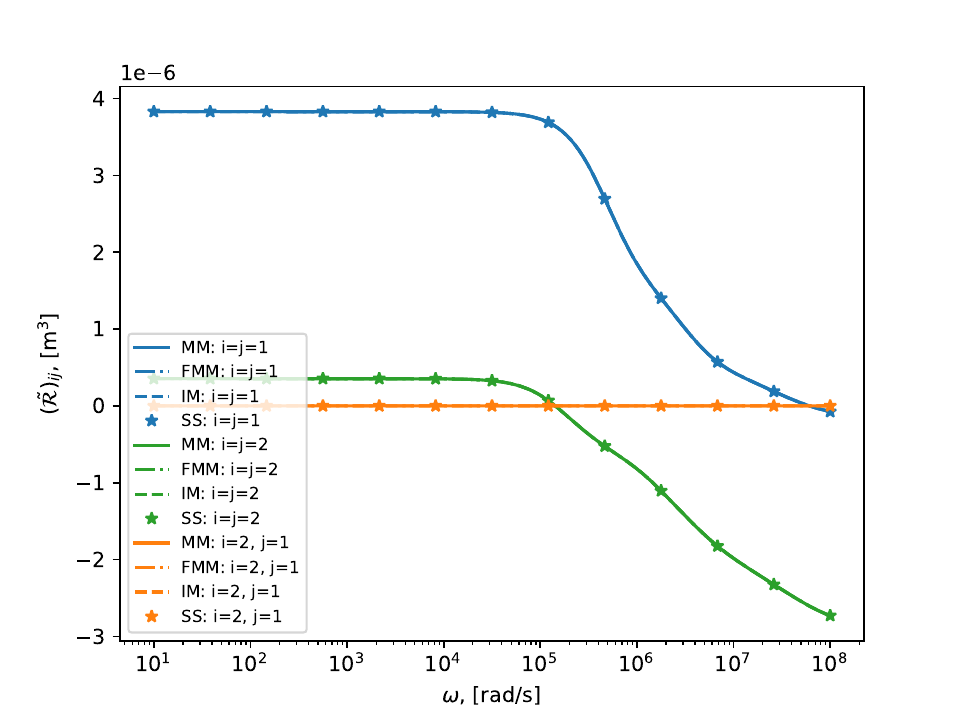} & \includegraphics[width=0.45\textwidth]{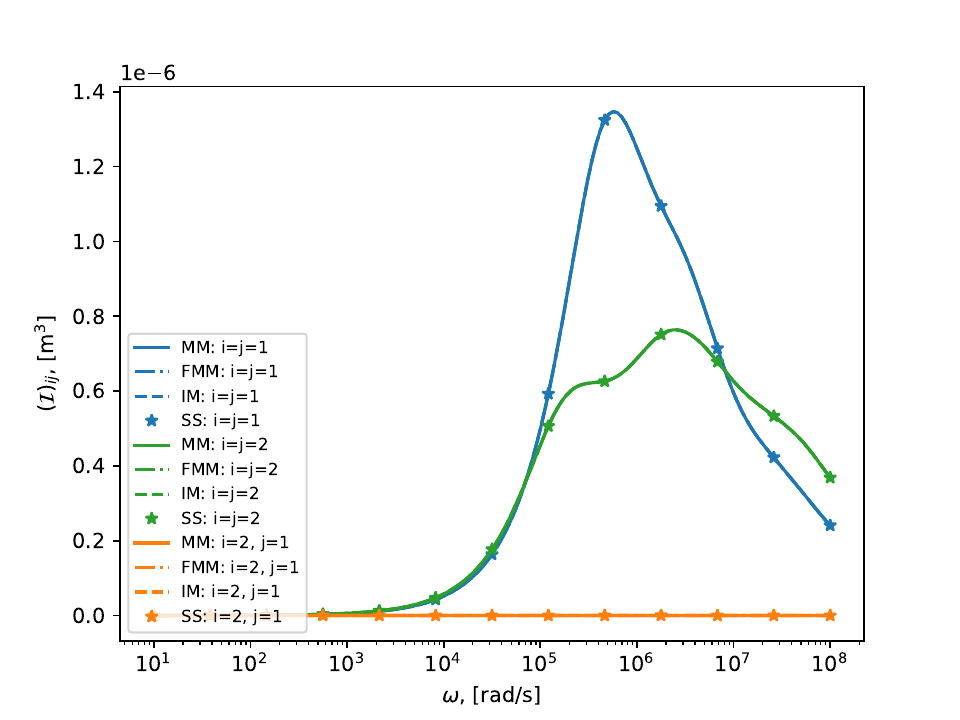}\\
(a) & (b) 
\end{array}$
\caption{{Thin conducting magnetic disk:  Showing a comparison between the original IM and the new faster FMM and MM approaches for the calculation of the MPT spectral signature using PODP $(a)$ $(\tilde{\mathcal R})_{ij}$ and $(b)$ $({\mathcal I})_{ij}$. }}

\label{fig:thin_disk_magnetic_tensors}
\end{figure}

Similarly to the magnetic sphere example, timings were performed comparing the MM, FMM, and IM methods for computing tensor coefficients for 160 output frequencies over the frequency range. As with the sphere example, a substantial improvement in speed is seen when going from the IM to FMM to MM methods, with each taking a total time of 19\,497, 1\,338, and 1\,038 seconds respectively with the post-processing time reducing from 18\,474 seconds to 385 seconds to 70 seconds. Timing plots are shown in figure~\ref{fig:magnetic_disk_times} and formatted in the same way as Figure~\ref{fig:magnetic_sphere_times}. The earlier observations about both computational time and memory usage again apply.
\begin{figure}[h!]
\centering
$\begin{array}{c}
\includegraphics[width=0.45\textwidth]{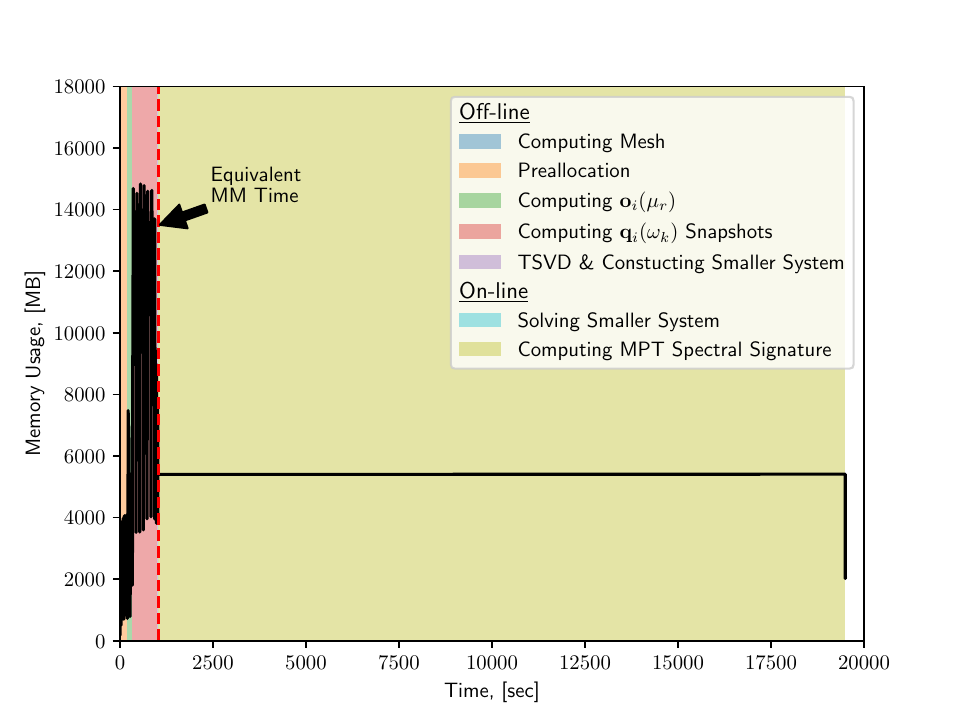} \\
(a)
\end{array}$
$\begin{array}{c c}
\includegraphics[width=0.45\textwidth]{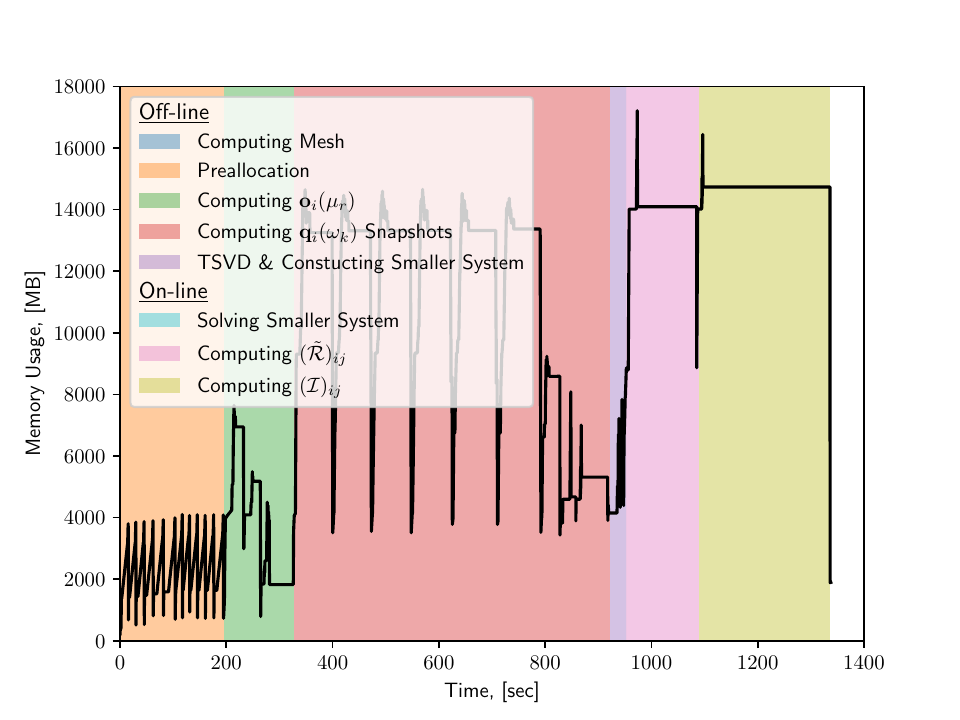} & 
\includegraphics[width=0.45\textwidth]{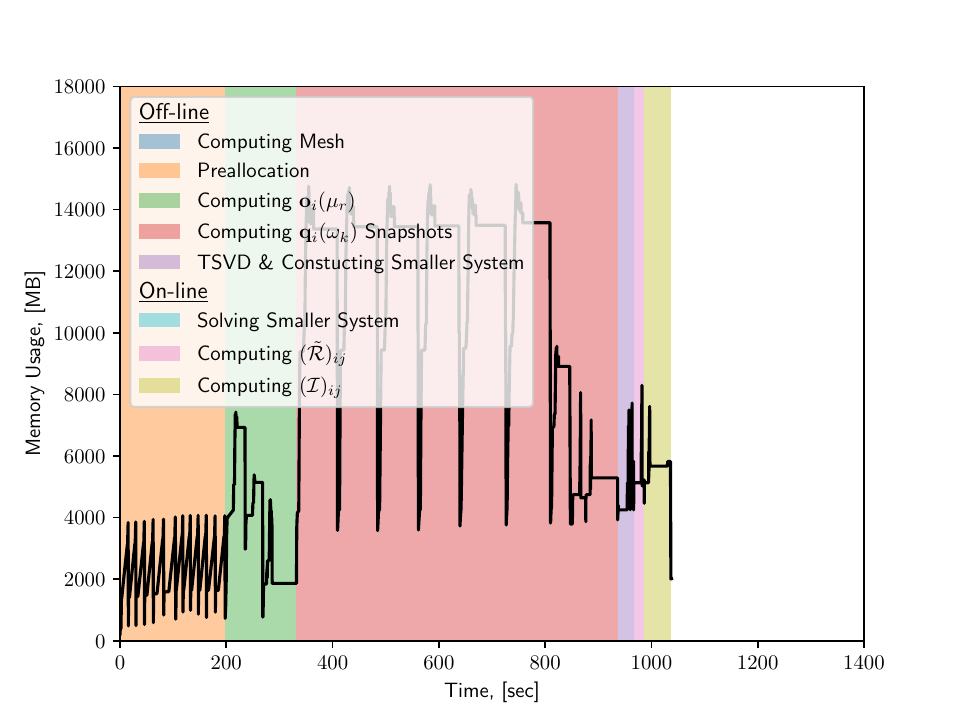} \\
(b) & (c)
\end{array}$
\caption{Thin conducting magnetic disk:
A comparison between time and memory usage for ($a$) IM, ($b$) FMM, and ($c$) MM methods showing substantial speed up for the on--line stage of POD. Timings are further broken down into common parts of the problem to correlate memory usage with specific tasks.}
\label{fig:magnetic_disk_times}
\end{figure}

To illustrate the adaptive procedure described by Algorithm~\ref{alg:adapt}, {the spectral signature for $(\tilde{\mathcal R})_{ij}$,  including the a-posteriori error certificates $(\tilde{\mathcal R}\pm \Delta)_{ij}$ obtained at different iterations, is shown  for the same magnetic disk  in Figure~\ref{error_iteration_examplere}}. Starting with the setup used in Figure~\ref{fig:thin_disk_magnetic_tensors}, and a stopping tolerance of $TOL_\Delta = 1\times 10^{-3}$, the figure shows
{iterations $k=1,2,3,4$ 
of the adaptive algorithm beginning with $N=13$ logarithmically spaced snapshots and resulting in $N=  15, 17,  19$ non-logarithmically} spaced snapshots, respectively, in the subsequent {three} iterations.  The convergence behaviour for certificates for $({\mathcal I} +\Delta)_{ij}$ is {identical} and, therefore, not shown. 
\begin{figure}[h!]
$\begin{array}{c c}
\includegraphics[width=0.45\textwidth]{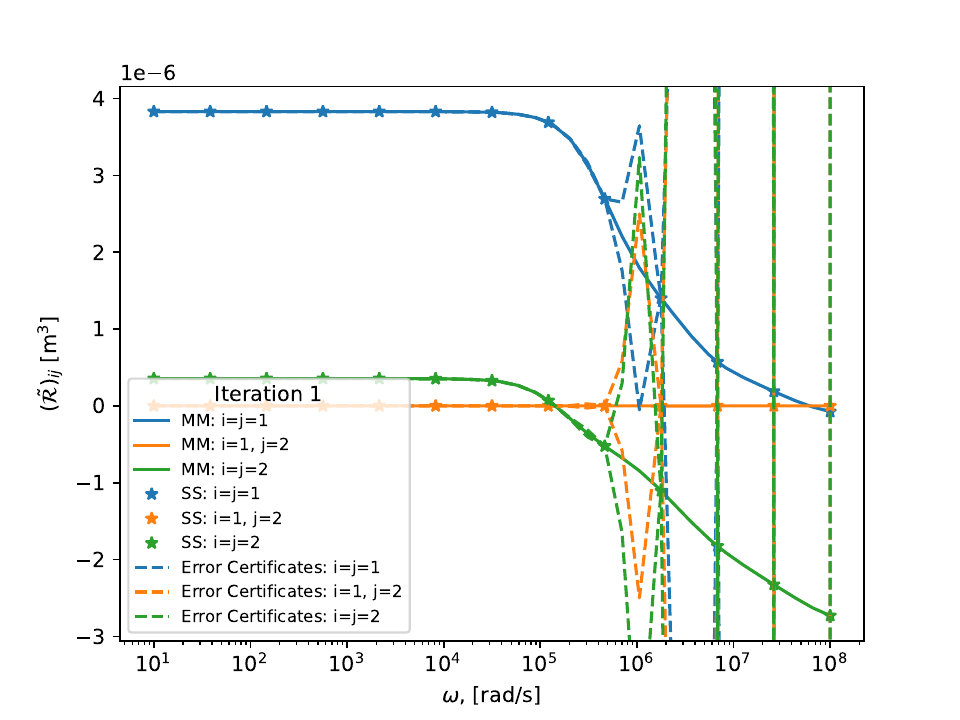} & \includegraphics[width=0.45\textwidth]{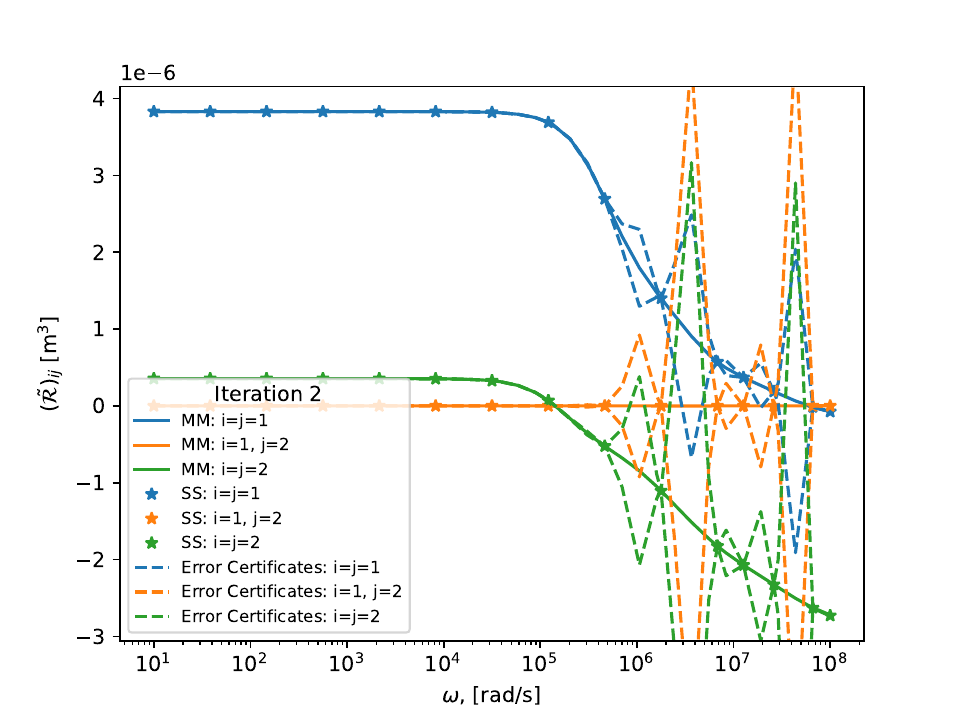} \\
(a) & (b) \\
\includegraphics[width=0.45\textwidth]{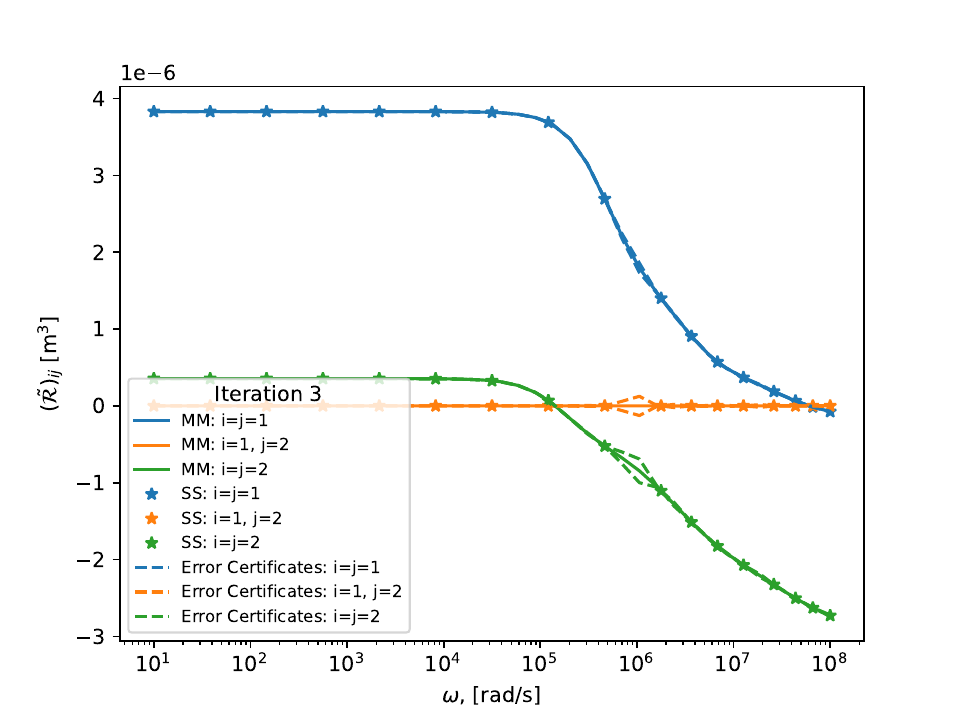} & \includegraphics[width=0.45\textwidth]{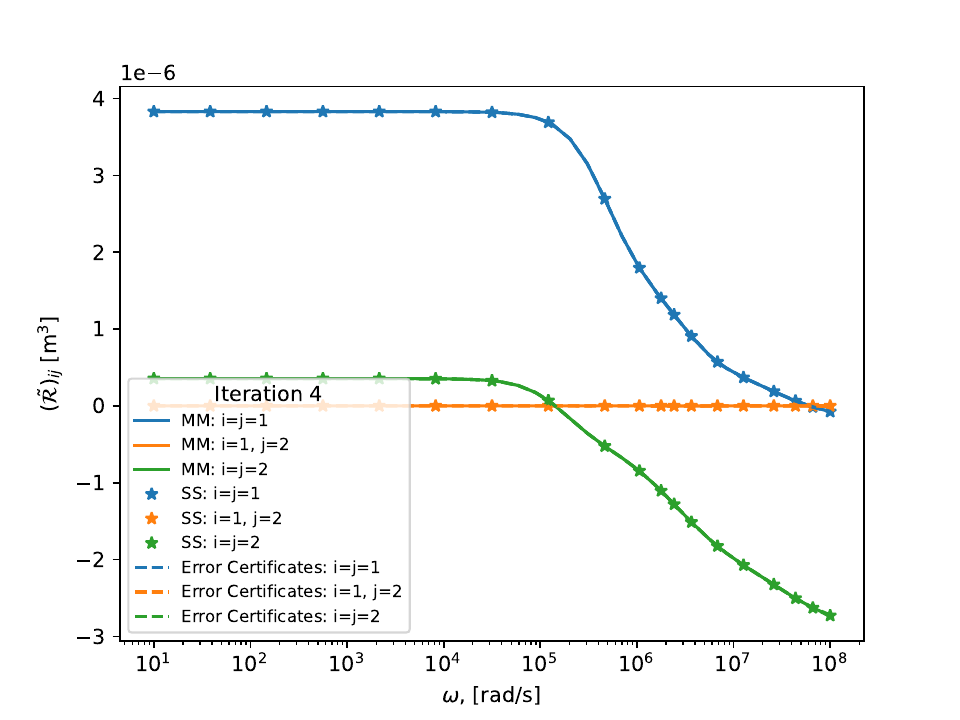} \\
(c) & (d)
\end{array}$
\caption{{Thin conducting magnetic disk:  Showing $(\tilde{\mathcal R})_{ij}$ and $(\tilde{\mathcal R} + \Delta)_{ij}$ obtained by applying the adaptive PODP Algorithm~\ref{alg:adapt} at iterations  $(a)$ $k=1$,  $(b)$ $k=2$,  $(c)$  $k=3$, and $(d)$  $k=4$.}}

\label{error_iteration_examplere}
\end{figure}

To illustrate the performance of the adaptive POD, compared to non-adaptive logarithmically spaced snapshots, Figure~\ref{fig:error_decay} shows the maximum error, $\Lambda $, against $N$, which, like the earlier sphere example, shows the benefits of the adaptive scheme over logarithmically spaced snapshots.

\begin{figure}[h!]
\centering
\includegraphics[width=0.45\textwidth]{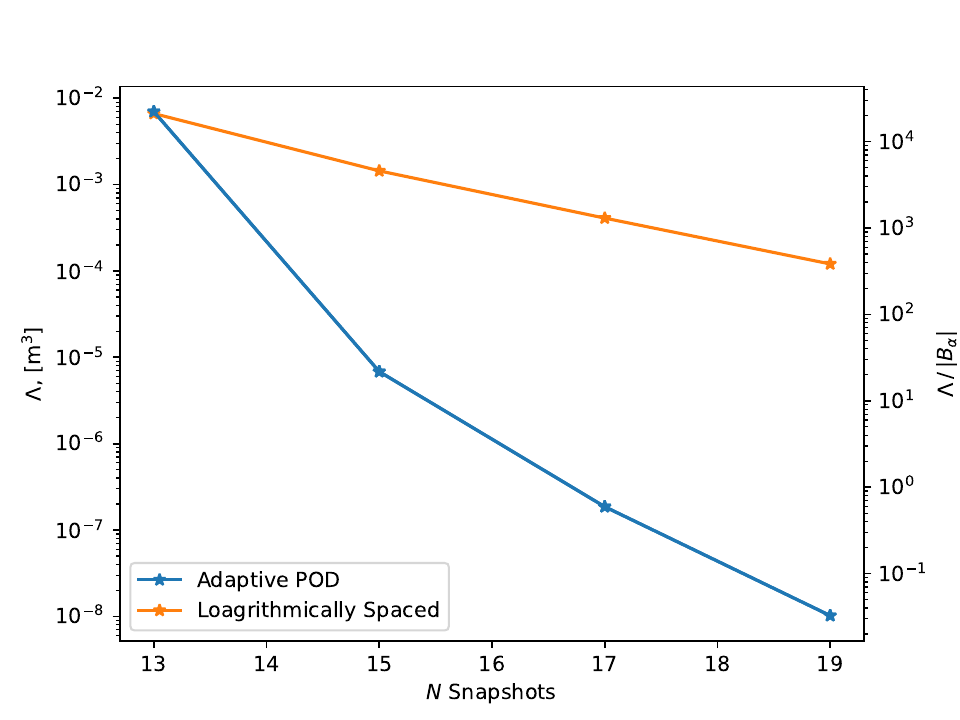}
\caption{{Thin conducting magnetic disk:  Showing  $\Lambda$ and $\Lambda / |B_\alpha|$ obtained by applying the adaptive PODP Algorithm~\ref{alg:adapt}   against $N$ compared with logarithmically spaced snapshot frequencies.}}
\label{fig:error_decay}
\end{figure}
Considering the limiting case of an infinitely thin conducting non-magnetic  disk  in the $x_1-x_3$ plane  and the corresponding limiting case when the  disk is magnetic,    {based on measurement observations for finitely thick disks~\cite{abdel, marsh2013}} (which is rotated to a disk in the $x_1-x_2$ plane for our situation), it has been proposed that the form of the MPT (when its coefficients are displayed as a matrix)  undergo the transition 
\begin{equation*}
\mathcal{M} [\alpha B, \omega, \sigma_*,1]= \begin{bmatrix}
0 & 0 & 0 \\ 0 & c & 0 \\ 0 & 0 & 0
\end{bmatrix}
 \to \mathcal{M} [\alpha B, \omega, \sigma_*,\mu_r \to \infty ] = \begin{bmatrix}
m & 0 & 0 \\ 0 & 0 & 0 \\ 0 & 0 & m
\end{bmatrix}
\end{equation*}
as $\mu_r$ becomes large.
{This proposal is investigated numerically in Figure~\ref{fig:thin_disk_mur}. for the aforementioned finitely thick disk geometry, but now with $\mu_r =1, 8,16, 64$, in turn.} In each case, different thickness prismatic layers were constructed according to the recipe proposed in Section~\ref{sect:boundarylayers} for a target value of $\omega =1 \times 10^8$ rad/s and then the MPT coefficients {are} obtained by using {the} PODP method for $N=13$ snapshot solutions at logarithmically spaced frequencies in the range $1\times 10^1\le \omega \le1\times 10^8$ rad/s and $TOL_\Sigma=1\times 10^{-6}$. Interestingly, despite the relative simplicity of the geometry and its homogeneous materials, {the presence of multiple local maxima in $({\mathcal I})_{ij}$  for the cases of $\mu_r =16$ and $\mu_r=64$ are observed}. This can be explained by the spectral theory of MPT spectral signatures~\cite{LedgerLionheart2019} where, in this case, the first term in the expansions in Lemma 8.5 is no longer dominant and multiple terms play an important role.

\begin{figure}
$\begin{array}{c c}
\includegraphics[width=0.45\textwidth]{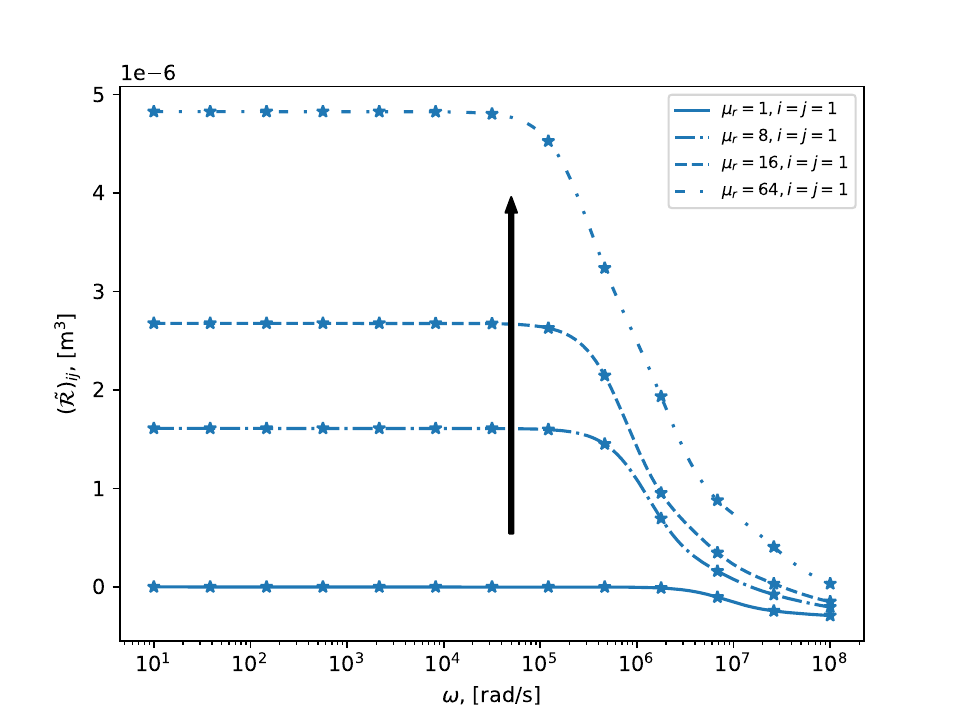} & 
\includegraphics[width=0.45\textwidth]{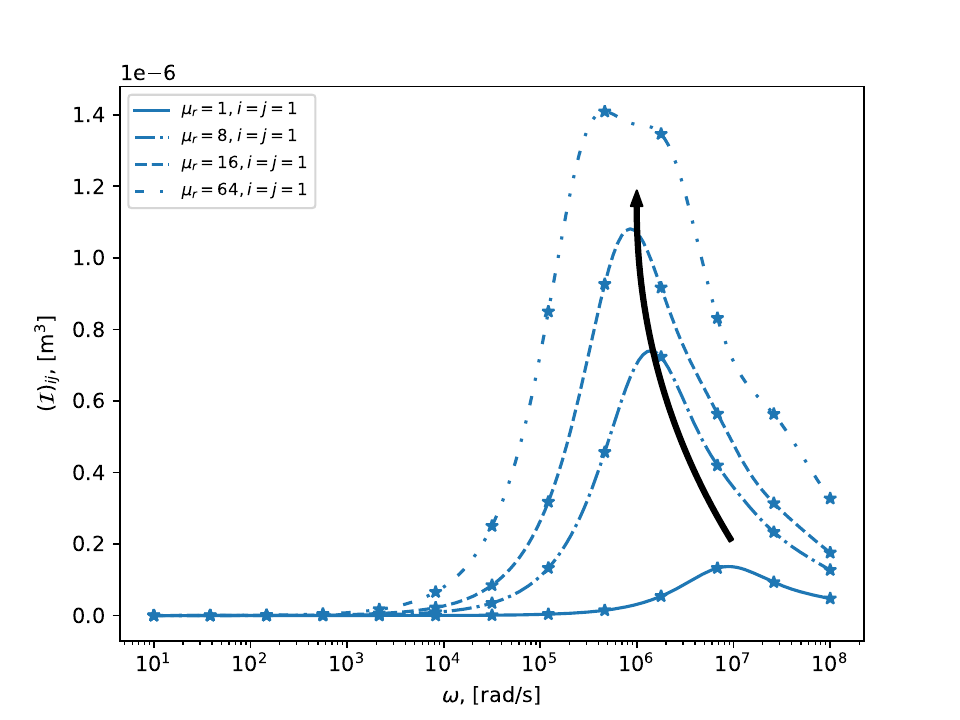} \\
(a) & (b) \\
\includegraphics[width=0.45\textwidth]{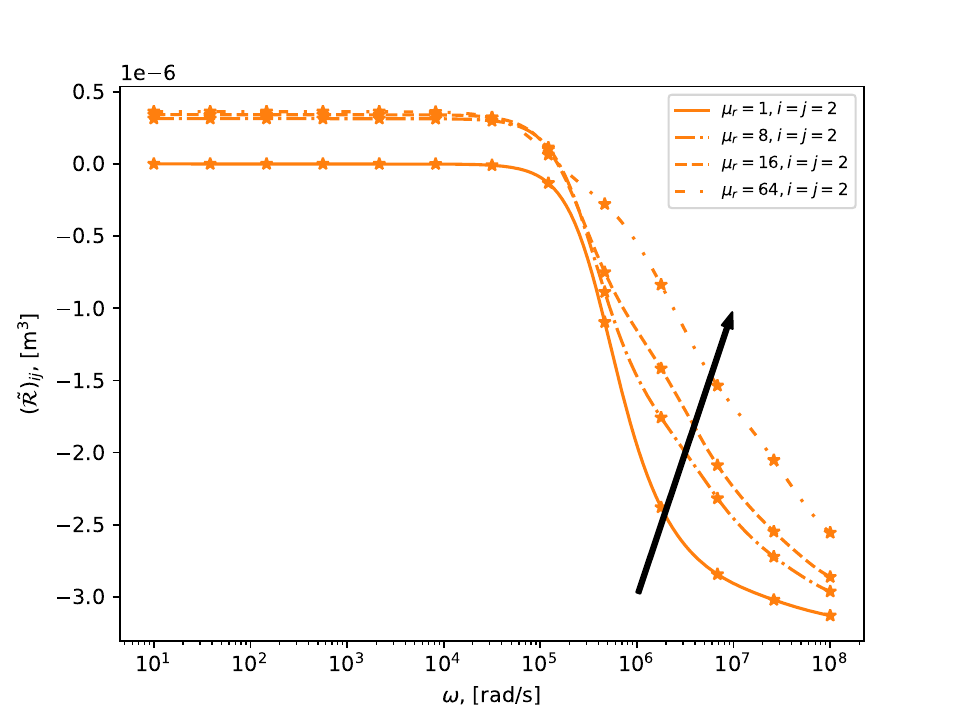} &
\includegraphics[width=0.45\textwidth]{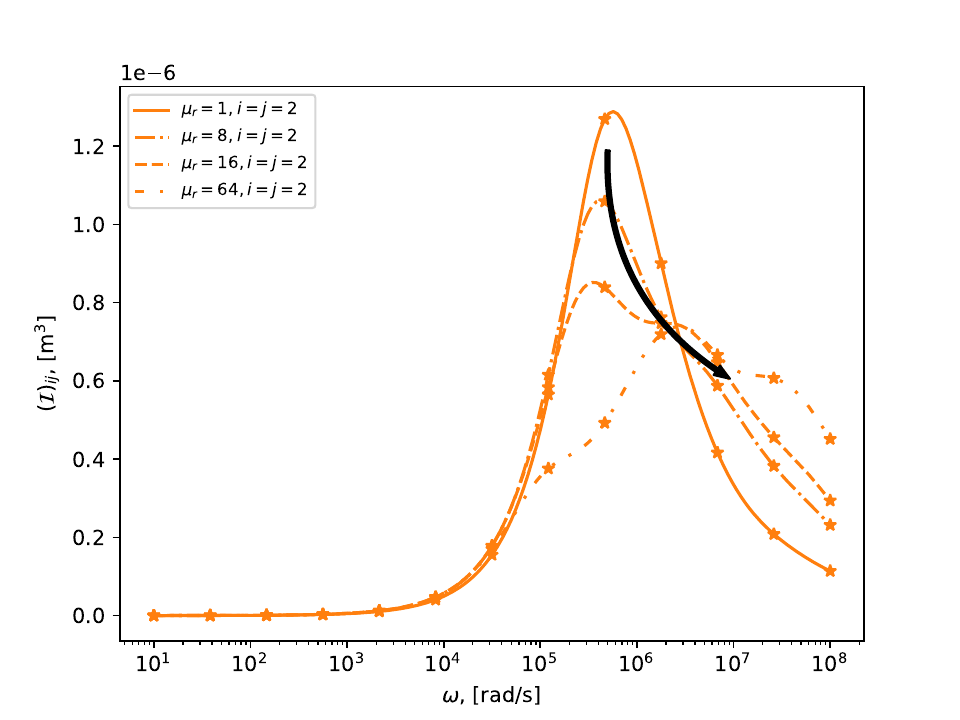} \\
(c) & (d) \\
\end{array}$
\caption{ {Thin conducting magnetic disk: Showing the effect of increasing $\mu_r$ with arrows indicating the evolution of the spectral signature with} ($a$) $(\tilde{\mathcal{R}})_{11}=(\tilde{\mathcal{R}})_{33}$, ($b$) $(\mathcal{I})_{11}=({\mathcal{I}})_{33}$, ($c$) $(\tilde{\mathcal{R}})_{22}$, and ($d$) $(\mathcal{I})_{22}$.}

\label{fig:thin_disk_mur}
\end{figure}

\section{Realistic Metallic Cleaver} \label{sect:cleaver}

One potential practical application is the use of metal detection to identify threat objects, such as knifes and firearms, in security screening applications, such as at event venues and transport hubs. To assist with this process, it is important to differentiate between concealed weaponry and other metallic clutter and an accurate  MPT characterisation of common objects is potentially of value. Therefore, as a more challenging example, we consider the case of a large metallic kitchen cleaver of the sort found in many homes and restaurants. In practice a cleaver often consists of a large thin blade made of either stainless or high-carbon steel and a plastic or wooden handle, often with the addition of small cylindrical rivets. In this work, we model the cleaver as only the metallic components; modelling only the steel blade and three copper rivets along the tang of the blade, with an illustration provided in figure~\ref{fig:cleaver_mesh}. {In the figure, we show ($a$) a surface triangulation of the full object $B$ and ($b$) a cut through of the full mesh with tetrahedra inside the object shown in green, prisms positioned at the surface of the object in blue, and tetrahedra outside the object in red.

\begin{figure}[!h]
\centering
$\begin{array}{c c}
\includegraphics[align=c, width=0.3\textwidth]{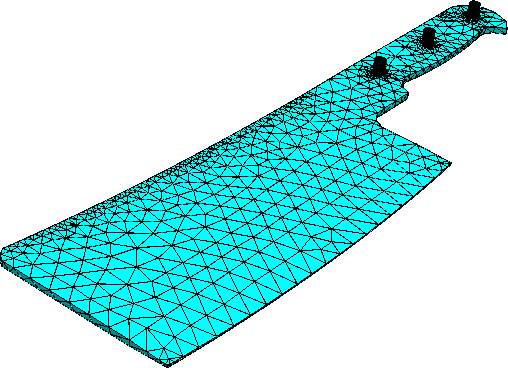} &
\includegraphics[align=c, width=0.3\textwidth]{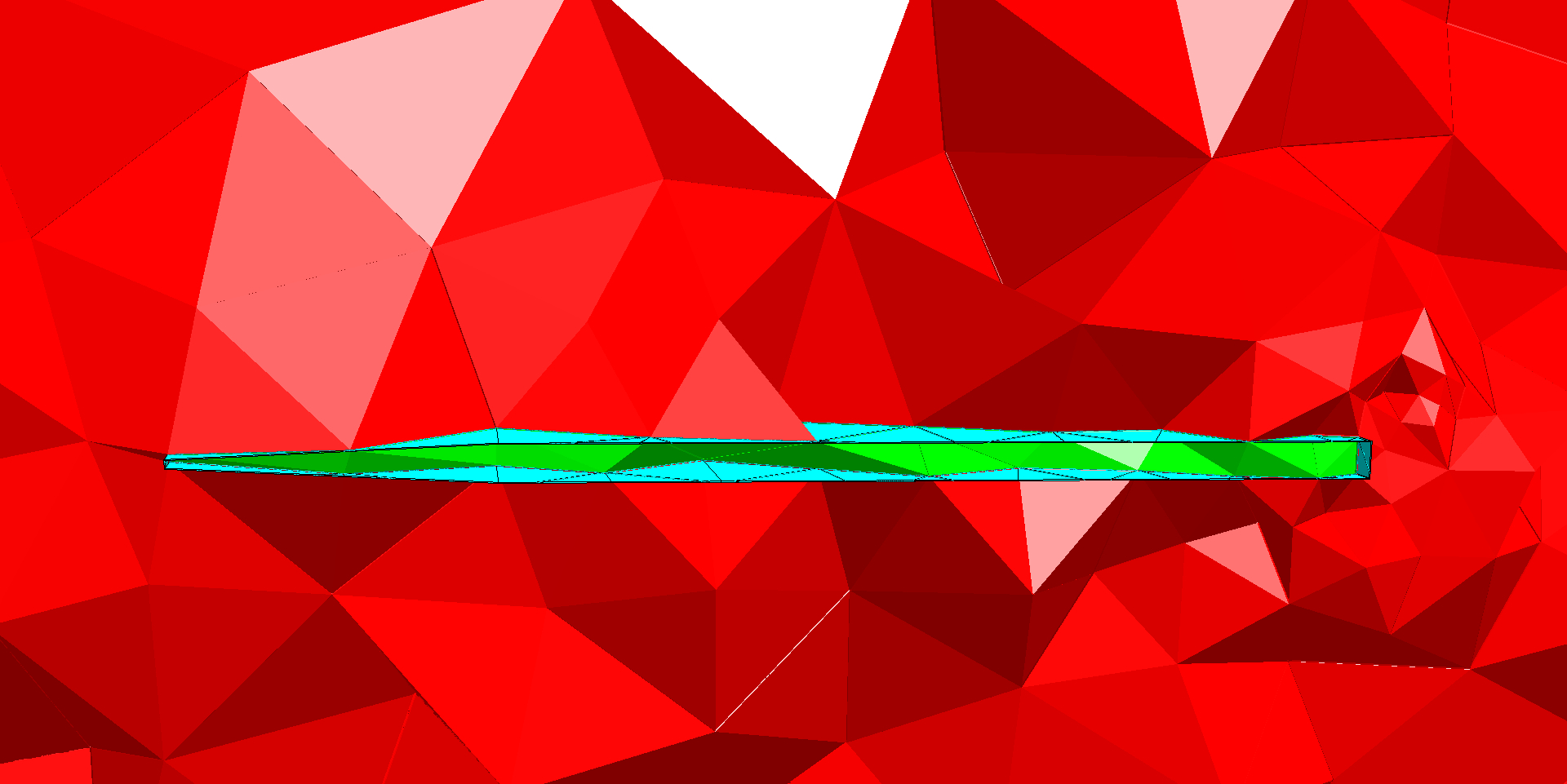}\\
(a) & (b) 
\end{array}$
\caption{Conducting magnetic cleaver: Mesh illustration of ($a$) the object $B$ showing both the steel blade and copper rivets with the prisms shown in blue and ($b$) a cut through of the mesh showing the interior tetrahedra, thin layer of prisms, and exterior tetrahedra. The mesh consists of 32\,023 tetrahedra and 6\,671 prisms.}
\label{fig:cleaver_mesh}
\end{figure}

For this example, the steel blade realistically modelled using $\sigma_* = 6\times 10^6$ S/m and $\mu_r = 100$ to model highly magnetic cast steel~\cite{conductivity}. Three copper rivets are introduced and are realistically modelled with conductivity $5.8\times 10^7$ S/m and $\mu_r = 1$, corresponding to 100$\%$ IACS, and the object is placed in a $[-1000, 1000]^3$ units non-conducting box, in the same way as the previous examples. Prismatic boundary layer elements are added to the steel blade section of the mesh following the geometric increasing strategy (and not to the copper rivets due to the much larger skin depths in the non-magnetic material), using a mesh consisting of 32\,023 unstructured tetrahedra and 6\,671 prisms which is converged when using $p=4$ for the $N=13$ POD snapshot frequencies over the range $1\times 10^1 \le \omega \le 1\times 10^8$ rad/s.

The object is configured such that there is a reflection symmetry in the $\boldsymbol{e}_1$ direction and the object has no additional rotational symmetries. This means that, $(\mathcal{M})_{13} = (\mathcal{M})_{23} = (\mathcal{M})_{31} = (\mathcal{M})_{32} = 0$, $(\mathcal{M})_{12} = (\mathcal{M})_{21}$, and $(\mathcal{M})_{11} \neq (\mathcal{M})_{22} \neq (\mathcal{M})_{33}$. For this reason, all three on-diagonal entries are shown and the non-zero off-diagonal entry are shown in Figure~\ref{fig:cleaver_tensors}, where excellent agreement between the MM, FMM, and IM methods are exhibited.

\begin{figure}[!h]
\centering
$\begin{array}{c c}
\includegraphics[width=0.45\textwidth]{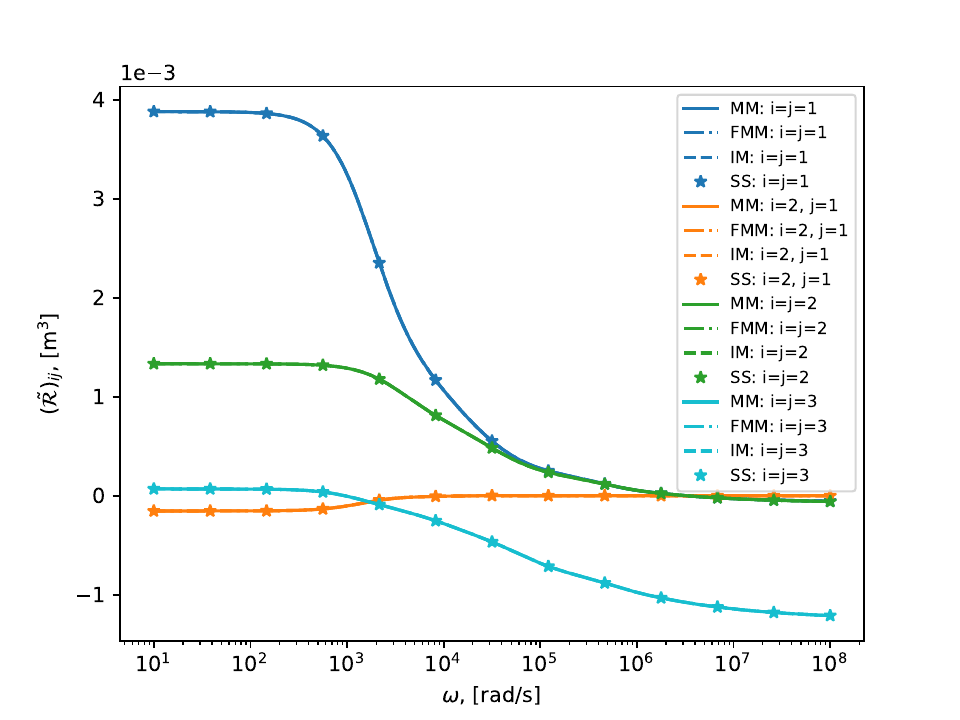} & \includegraphics[width=0.45\textwidth]{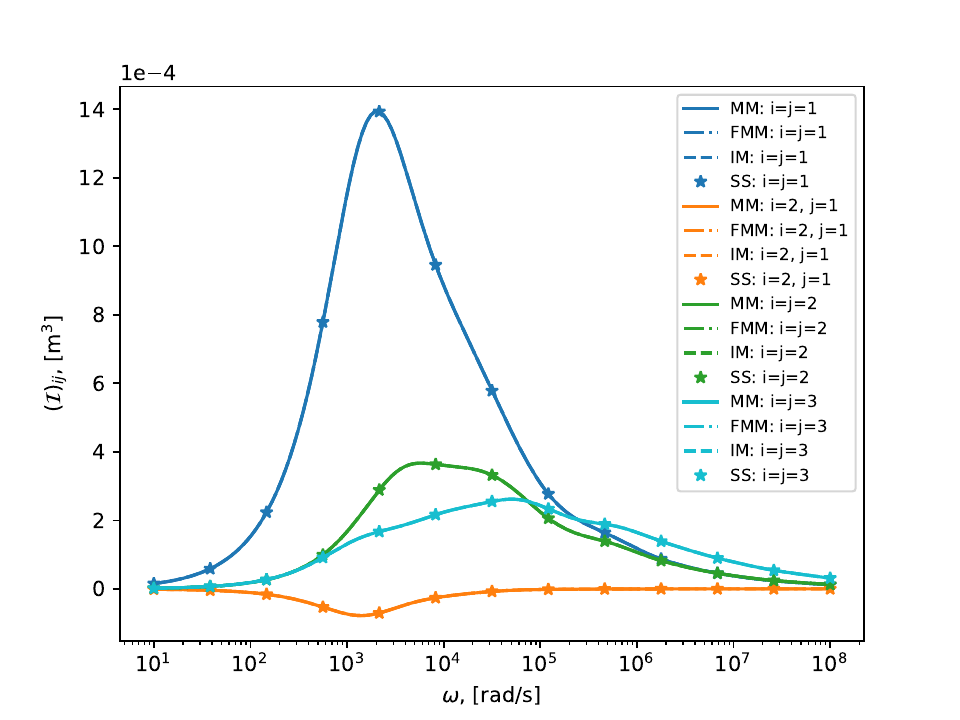}\\
(a) & (b) 
\end{array}$
\caption{Conducting magnetic cleaver: Showing a comparison between the original IM and the new faster FMM, and MM approaches for the calculation of the MPT spectral signature using PODP $(a)$ $(\tilde{\mathcal R})_{ij}$ and $(b)$ $({\mathcal I})_{ij}$.}
\label{fig:cleaver_tensors}
\end{figure}

Similarly to Figures~\ref{fig:magnetic_sphere_times} and~\ref{fig:magnetic_disk_times}, the MM approach is shown, in Figure~\ref{fig:cleaver_times}, to be significantly faster than both the IM and FMM methods. The total computation time, including both off-line and on-line stages for the IM, FMM, and MM methods is 50\,846 seconds, 8\,214 seconds, and 6\,603 seconds for the IM, FMM, and MM methods respectively, with computing tensor coefficients contributing 44\,415 seconds, 1\,837 seconds, and 208 seconds respectively, with the earlier observations about both time and memory usage also applying in this case. 

\begin{figure}[h!]
\centering
$\begin{array}{c}
\includegraphics[width=0.45\textwidth]{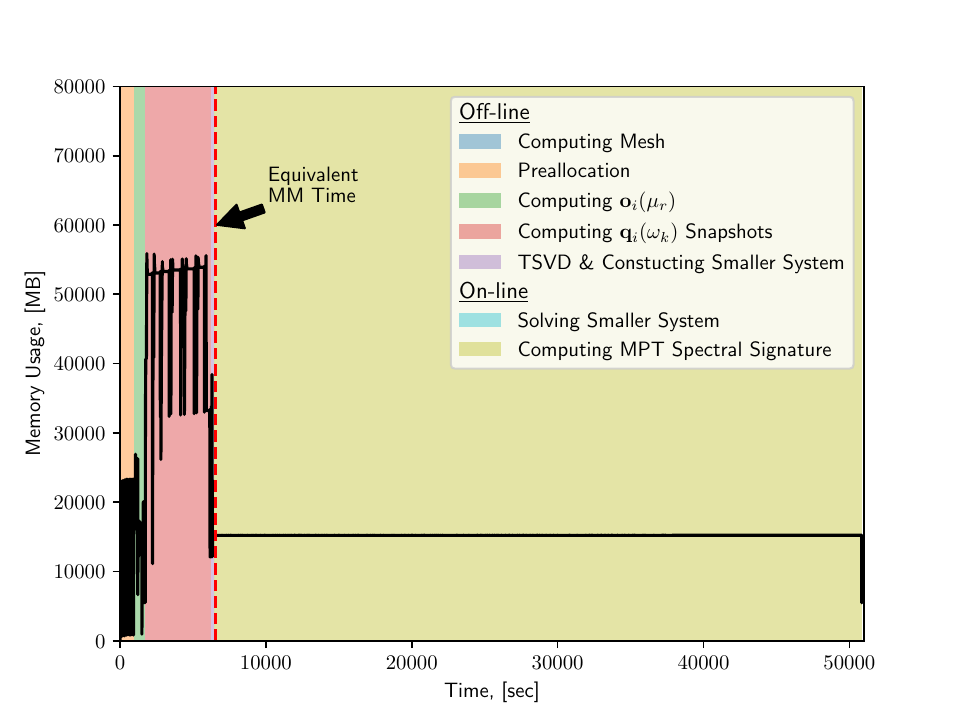} \\
(a)
\end{array}$
$\begin{array}{c c}
\includegraphics[width=0.45\textwidth]{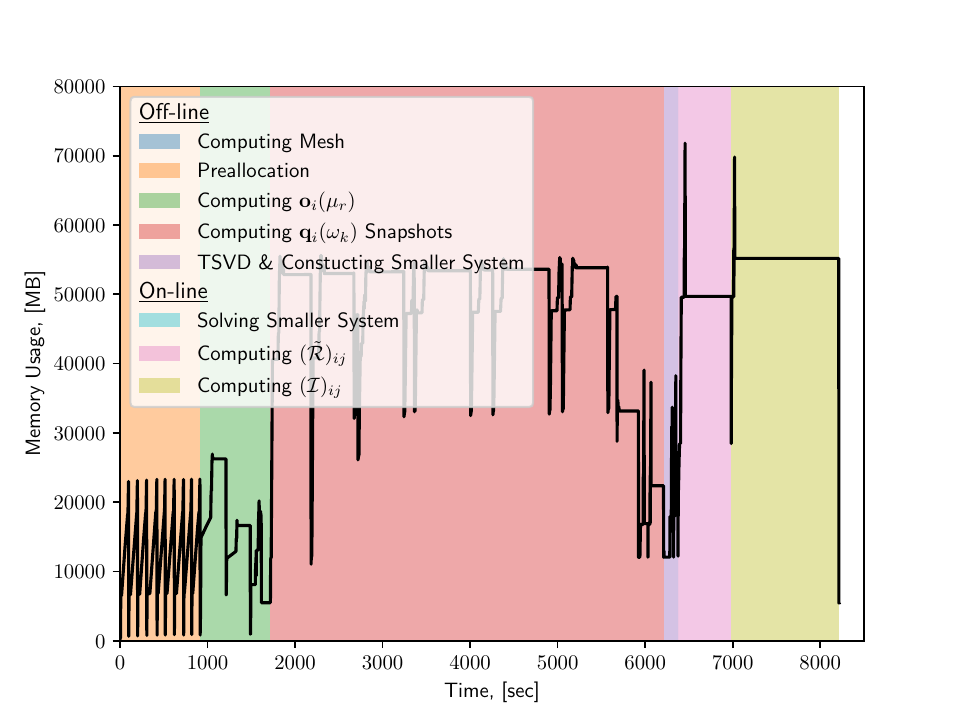} & 
\includegraphics[width=0.45\textwidth]{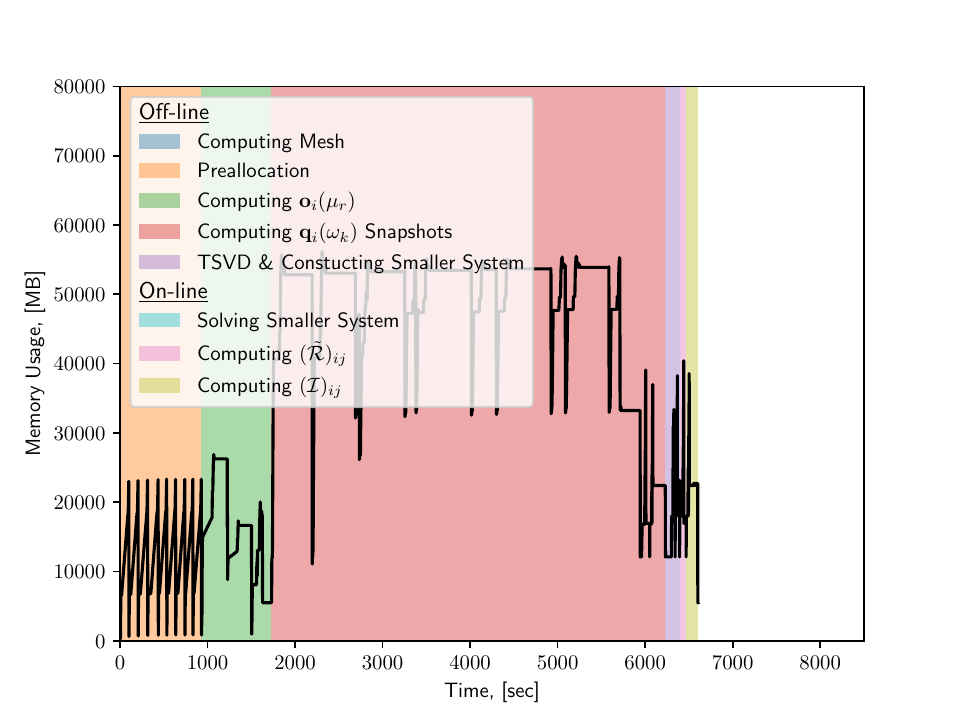} \\
(b) & (c)
\end{array}$

\caption{Conducting magnetic cleaver:
A comparison between time and memory usage for ($a$) IM, ($b$) FMM, and ($c$) MM methods showing substantial speed up for the on--line stage of POD. Timings are further broken down into common parts of the problem to correlate memory usage with specific tasks.}
\label{fig:cleaver_times}
\end{figure}

To illustrate the adaptive POD method described in Algorithm~\ref{alg:adapt} we show, in Figure~\ref{fig:cleaver_errors}, $k=1,2,4,8$ iterations of the adaptive POD, after which the error certificates become indistinguishable from the computed tensor coefficients on this scale. 
As with the previous two examples, we illustrate the adaptive procedure using the spectral signature for $(\tilde{\mathcal R})_{ij}$, including the a-posteriori error certificates $(\tilde{\mathcal R}\pm \Delta)_{ij}$ obtained at different iterations of the adaptive algorithm beginning with $N=13$ logarithmically spaced snapshots and resulting in $N= 13, 15, 21, 29$ non-logarithmically spaced snapshots, respectively, in the subsequent iterations.
In this case, due to the large size of the inhomogeneous object and highly magnetic materials, the off--line stage of the POD computes 13 logarithmically spaced snapshots over the range $1\times 10^1 \le \omega \le 1\times 10^8$ rad/s. However, the adaption and error estimates are applied over the smaller frequency range $1\times 10^1 \le \omega \le 1\times 10^7$ rad/s. This is in order to avoid edge effects associated with the higher frequencies and is not of great concern since the eddy current model will break down at high frequencies in any case.

\begin{figure}[!h]
$\begin{array}{c c}
\includegraphics[width=0.45\textwidth]{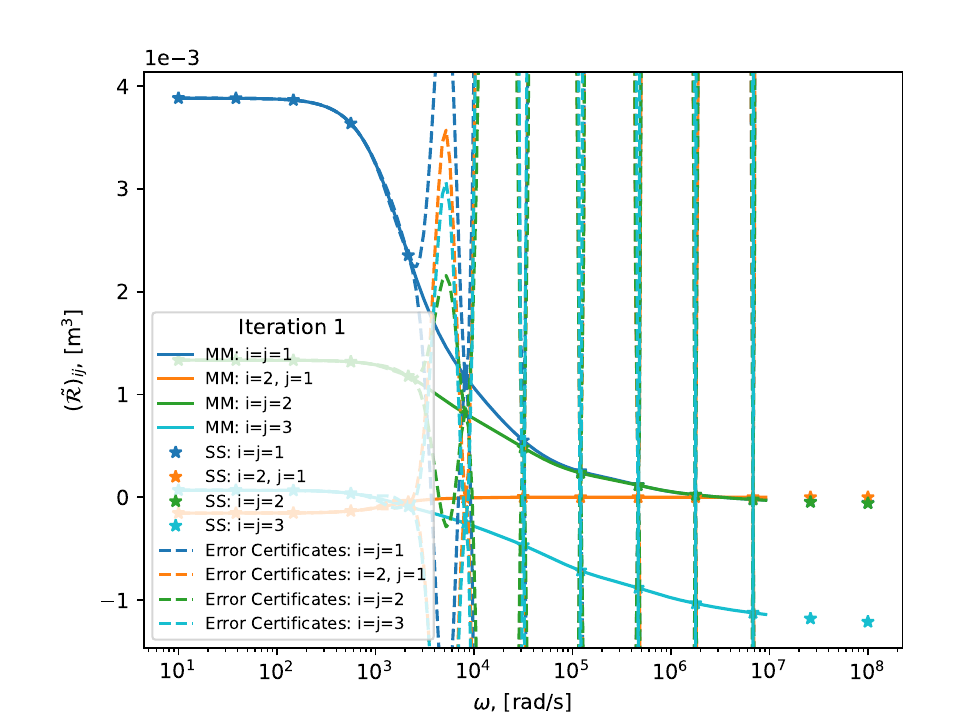} &
\includegraphics[width=0.45\textwidth]{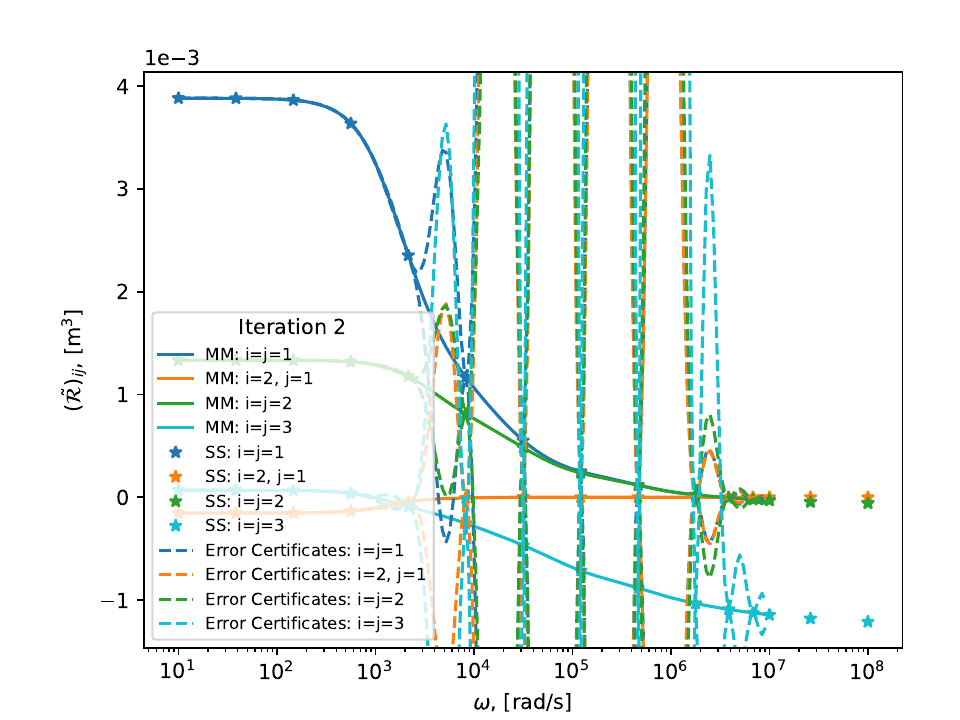} \\
(a) & (b) \\
\includegraphics[width=0.45\textwidth]{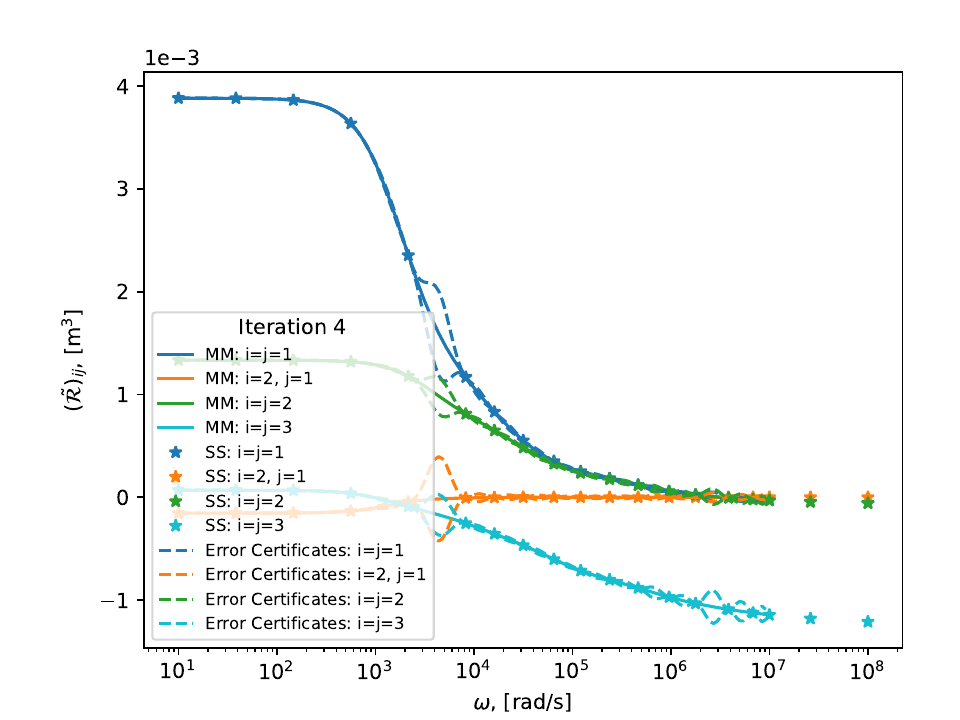} & 
\includegraphics[width=0.45\textwidth]{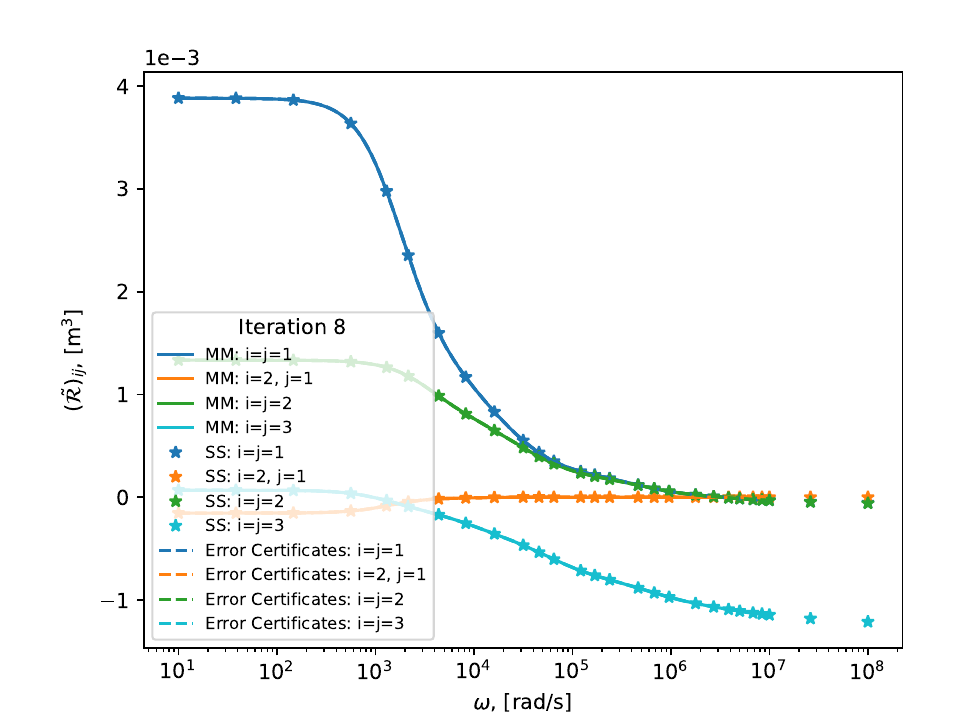} \\
(c) & (d)
\end{array}$
\caption{Conducting magnetic cleaver:  
Showing $(\tilde{\mathcal R})_{ij}$ and $(\tilde{\mathcal R} + \Delta)_{ij}$ obtained by applying the adaptive PODP Algorithm~\ref{alg:adapt} at iterations  $(a)$ $k=1$,  $(b)$ $k=2$,  $(c)$  $k=4$, and $(d)$  $k=8$.}

\label{fig:cleaver_errors}
\end{figure}

Finally, to illustrate the difficulties associated with this object, contour plots of normalised $\lvert {\boldsymbol J}_3^{e,hp} \rvert = \lvert\mathrm{Re}(\mathrm{i}\omega\sigma_* \boldsymbol{\theta}^{(1, hp)}_3)\rvert$ are included to illustrate the decreasing skin depths as $\omega$ is increased from $\omega = 1\times 10^2,$ to $1\times 10^4$, and $1\times 10^6$ rad/s in Figure~\ref{fig:cleaver_contours}. In each case, the results are normalised between 0 and 1 to allow for easier comparison, showing the decreasing skin depths and frequency dependence of the induced eddy currents, which become increasingly concentrated along the surface of the object with increasing frequency. To show the rotation of the induced eddy currents clearly throughout the entire blade, we demonstrate, using the $\omega=1\times 10^4$ rad/s example, streamlines of the induced eddy currents throughout the blade in Figure~\ref{fig:cleaver_currents}. Similarly to Figure~\ref{fig:cleaver_contours}b the eddy currents decay quickly towards the center of the object. 

\begin{figure}[!h]
\centering
$\begin{array}{c c}
\includegraphics[width=0.45\textwidth]{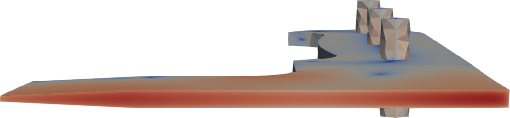} & (a) \\
\includegraphics[width=0.45\textwidth]{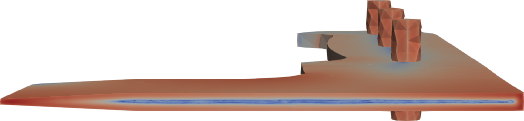} & (b) \\
\includegraphics[width=0.45\textwidth]{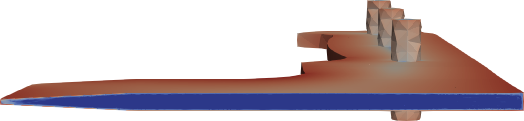} & (c) \\
\includegraphics[width=0.45\textwidth]{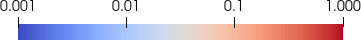}
\end{array}$
\caption{Conducting magnetic cleaver: Contour plots of $\lvert {\boldsymbol J}_3^{e,hp} \rvert = \lvert\mathrm{Re}(\mathrm{i}\omega\sigma_* \boldsymbol{\theta}^{(1, hp)}_3)\rvert$ for the frequencies ($a$) $\omega = 1\times 10^2$ rad/s, ($b$) $\omega = 1\times 10^4$ rad/s, and ($c$) $\omega = 1\times 10^6$ rad/s. The figure shows both the decreasing as $\omega$ increases and demonstrates the necessity of using the thin layers of prisms advocated for in this work.}
\label{fig:cleaver_contours}
\end{figure}

\begin{figure}[!h]
\centering
$\begin{array}{c}
\includegraphics[width=0.45\textwidth]{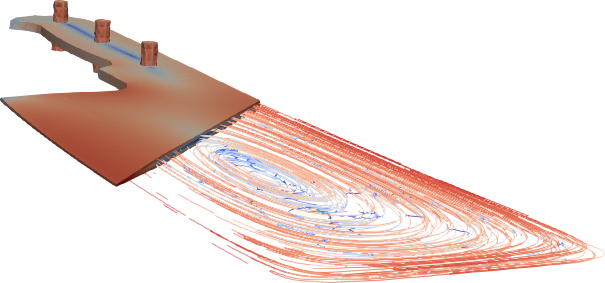}\\
\includegraphics[width=0.45\textwidth]{Cleaver_colourbar}
\end{array}$

\caption{Conducting magnetic cleaver: Contour plot of $\lvert {\boldsymbol J}_3^{e,hp} \rvert = \lvert\mathrm{Re}(\mathrm{i}\omega\sigma_* \boldsymbol{\theta}^{(1, hp)}_3)\rvert$ for the frequency $\omega = 1\times 10^4$ rad/s. The induced eddy currents are illustrated using streamlines, showing both the rotation of the currents and the decay in magnitude towards the center of the object.}
\label{fig:cleaver_currents}
\end{figure}

\section{Conclusions}\label{sect:Conclusions}

In this work, a new computational formulation for efficiently computing MPT coefficients from POD predictions has been developed. This, in turn, 
has led to  improvements in computational performance associated with obtaining MPT spectral signature characterisations of complex and highly magnetic conducting objects including magnetic disks of varying magnetic permeability and a realistic inhomogeneous magnetic cleaver. Three contrasting approaches IM, FMM, and MM for computing the MPT coefficients with POD are presented, which, up to rounding error, all produce identical results. While IM is acceptable for small problems, it scales poorly with large numbers of output frequencies and for large problems. FMM offers an improvement over IM by avoiding having to recompute integrals for each output frequency and MM is the most superior as the computation of MPT coefficients for new frequencies becomes independent of the underlying discretisation. These performance gains have been obtained by appealing to the explicit expressions for the MPT coefficients and writing them in terms of finite element matrices and solution vectors. 


The paper has also proposed an enhancement to our previous POD scheme by incorporating  adaptivity  to choose off-line snapshot frequencies based on an a--posteriori error estimate to produce an accurate on-line prediction of the MPT spectral signature.
This addresses an important shortcoming of our previous work by using a-posteriori error estimates in order to choose snapshot frequencies and obtain accurate results to some desired accuracy. The results included illustrate that the adaptive method produces results that are four orders of magnitude more accurate compared to using the same number of logarithmically spaced frequency snapshots. 

In addition, a simple recipe for choosing  the number and thicknesses of prismatic boundary layers has been proposed which provides a concrete way of constructing discretisations for modelling magnetic objects when combined with $p$--refinement. By considering a magnetic sphere, our results show that by using {two} prismatic layers of suitable thicknesses, and employing $p$--refinement, an exponential rate of convergence can be obtained for the MPT coefficients. In addition, our strategy achieves  a relative error of $E<1\times 10^{-3}$ for the MPT coefficients over a wide range of materials and frequency excitations. Results have also been included for a magnetic disk and realistic cleaver with inhomogeneous magnetic materials.

To show the benefits of our improved methodology we also demonstrate both timings and the adaptive POD using a large magnetic multi-material cleaver. This model, and the included contour plots, show the importance of using prismatic layers to capture the thin electromagnetic skin depths, and the importance of our new MPT calculation formulation.


The procedures presented in this work are expected  to be invaluable for constructing large dictionaries of MPT characterisations of complex in-homogeneous real-world metallic objects. An open dataset of these characterisations will be released in the near future. 

\section{Acknowledgements}\label{sect:Ack}
The authors would like to thank Prof. Peyton, Prof. Lionheart and Dr Davidson for their helpful discussions and comments on polarizability tensors.
The authors are grateful for the financial support received from the Engineering and Physical Science Research Council (EPSRC, U.K.) through the research grant EP/V009028/1.

\bibliographystyle{elsarticle-num}

\bibliography{paperbib_short}

\end{document}